\documentclass[twoside]{article}
\usepackage{amssymb, floatflt}
\input enum.mac
\renewcommand{\copyrightheading}[1]{}

\runninghead{Enumerating the Prime Alternating Knots, Part I }
{Enumerating the Prime Alternating Knots, Part I }
\begin{document}
\copyrightheading{}		    

\fpage{1}
\centerline{\bf Enumerating the Prime Alternating Knots, Part I}
\centerline{\footnotesize Stuart Rankin, John Schermann, Ortho Smith\fnm{$\dagger$}}
\fnt{$\dagger$}{Partially supported by the Applied Mathematics Department at the University of
Western Ontario, Wolfgang Struss, Pei Yu, Turab Lookman, Ghia Flint and Dean Harrison.}
\centerline{\footnotesize Department of Mathematics, University of Western Ontario}
\centerline{\footnotesize srankin@uwo.ca, johns@csd.uwo.ca, dsmith6@uwo.ca}
\vskip20pt
\def\ArnoldEtAl{1}
\def\calvo{2}
\def\conway{3}
\def\DasbachHougardy{4}
\def\FrayMendez{5}
\def\DowkerThistlethwaite{6}
\def\ErnstSumners{7}
\def\HosteThistlethwaiteWeeks{8}
\def\kirkman{9}
\def\kirkmantwo{10}
\def\gomez{11}
\def\lickorish{12}
\def\lickorishtwo{13}
\def\little{14}
\def\MenascoThistlethwaite{15}
\def\MenascoThistlethwaitetwo{16}
\def\perko{17}
\def\perkotwo{18}
\def\tait{19}
\def\thistlethwaite{20}
\def\SundbergThistlethwaite{21}
\def\welsh{22}

\section{Abstract}
 The enumeration of prime knots has a long and storied history, beginning
 with the work of T. P. Kirkman [\kirkman,\kirkmantwo], C. N. Little [\little],
 and P. G. Tait [19] in the late 1800's, and continuing through to the
 present day, with significant progress and related results provided along
 the way by J. H. Conway [\conway], K. A. Perko [\perko, \perkotwo], M. B.
 Thistlethwaite [\DowkerThistlethwaite, \HosteThistlethwaiteWeeks,
 \MenascoThistlethwaite, \MenascoThistlethwaitetwo, \thistlethwaite],
 C. B. Dowker [\DowkerThistlethwaite],  J. Hoste [\ArnoldEtAl,
 \HosteThistlethwaiteWeeks], J. Calvo [\calvo], W. Menasco
 [\MenascoThistlethwaite, \MenascoThistlethwaitetwo], W. B. R. Lickorish
 [\lickorish, \lickorishtwo], J. Weeks [\HosteThistlethwaiteWeeks] and
 many others. Additionally, there have been many efforts to establish bounds
 on the number of prime knots and links, as described in the works of O.
 Dasbach and S. Hougardy [\DasbachHougardy], D. J. A. Welsh [\welsh], C.
 Ernst and D. W. Sumners [\ErnstSumners], and C. Sundberg and M.
 Thistlethwaite [\SundbergThistlethwaite] and others.

 In this paper, we provide a solution to part of the enumeration problem,
 in that we describe an efficient inductive scheme which uses a total of
 four operators to generate all prime alternating knots of a given
 minimal crossing size, and we prove that the procedure does in fact
 produce them all. The process proceeds in two steps,
 where in the first step, two of the four operators are applied to the 
 prime alternating knots of minimal crossing size $n$ to produce
 approximately 98\% of the prime alternating knots of minimal crossing
 size $n+1$, while in the second step, the remaining two operators are
 applied to these newly constructed knots, thereby producing the remaining
  prime alternating knots of crossing size $n+1$. The process begins with
 the prime alternating knot of four crossings, the figure eight knot. 

 In the sequel, we provide an actual implementation of our procedure,
 wherein we spend considerable effort to make the procedure efficient.
 One very important aspect of the implementation is a new way of
 encoding a knot. We are able to assign an integer array (called the master array) to a
 prime alternating knot in such a way that each regular projection, or
 plane configuration, of the knot can be constructed from the data in the
 array, and moreover, two knots are equivalent if and only if their
 master arrays are identical. A fringe benefit of this scheme is a
 candidate for the so-called ideal configuration of a prime alternating
 knot.
 
 We have used this generation scheme to enumerate the prime alternating
 knots up to and including those of 19 crossings. The knots up to and
 including 17 crossings produced by our generation scheme concurred with
 those found by M. Thistlethwaite, J. Hoste and J. Weeks (see
 [\HosteThistlethwaiteWeeks]).

 The current implementation of the algorithms involved in the generation
 scheme allowed us to produce the $1,769,979$ prime alternating knots of
 17 crossings on a five node beowulf cluster in approximately $2.3$ hours,
 while the time to produce the prime alternating knots up to and including
 those of 16 crossings totalled approximately 45 minutes. The prime
 alternating knots at 18 and 19 crossings were enumerated using the 48 node
 Compaq ES-40 beowulf cluster at the University of Western Ontario (we
 also received generous support from Compaq at the SC 99 conference). The
 cluster was shared with other users and so an accurate
 estimate of the running time is not available, but the generation of the
 $8,400,285$ knots at 18 crossings was completed in 17 hours, and the
 generation of the $40,619,385$ prime alternating knots at 19 crossings
 took approximately 72 hours. With the improvements that are described in
 the sequel, we anticipate that the knots at 19 crossings will be
 generated in not more than 10 hours on a current Pentium III personal
 computer equipped with 256 megabytes of main memory.

\section{Introduction}
  A knot is a smooth embedding of the unit sphere $S^1$ into $\mathbb{R}^3$,
  and a knot configuration is a projection of a knot into a plane with
  preimages of size at most two for which a point has preimage of size two only
  if the point is the image of a crossing. Two knots are said to be
  equivalent if there exists a homeomorphism of $\mathbb{R}^3$ onto itself
  taking one of the knots to the other. It is conventional to indicate
  on a knot configuration the additional information required to construct a
  knot equivalent to the one from which the diagram was made. This
  information takes the form of over/under passes, indicated by a break
  in the curve in a neighborhood of the point in the preimage which
  is furthest from the projection plane. This paper provides an effective
  method for the construction (and hence enumeration) of all prime
  alternating knots. We rely on the fact that two prime alternating knots
  are equivalent if and only if they are flype equivalent (the Tait flyping
  conjecture, proven by W. W. Menasco and M. B. Thistlethwaite [15]).

  Up to the present time, the enumeration of knots has primarily relied on
  a process of representing each configuration of a knot by a sequence of
  integers. The most efficient of these representation schemes is
  the Dowker-Thistlethwaite (see [\DowkerThistlethwaite]) encoding
  scheme, whereby each configuration of an alternating knot of $n$ crossings
  is represented by a permutation of the even integers from $2$ to $2n$.
  The standard approach to prime alternating knot enumeration is essentially
  to produce permutations of the even integers from $2$ to $2n$, where $n$
  is the crossing size, then eliminate those permutations that are not the
  Dowker-Thistlethwaite code of any prime alternating knot, and finally,
  identify amongst the remaining permutations that represented the same
  alternating knot.

  In this work, we introduce an alternative approach which utilizes a new
  knot configuration encoding, which we call the group code for the knot
  configuration. 

  We then introduce an inductive scheme which uses a total of four operators
  to generate the prime alternating knots of minimal crossing size $n+1$
  from those of minimal crossing size $n$. Each stage requires two steps.
  In the first step, we apply two of the four operators to the prime
  alternating knots of minimal crossing size $n$ to produce a subset of the
  set of prime alternating knots of minimal crossing size $n+1$. In
  our computations so far (up to the knots of 19 crossings), we have found
  that over $98\%$ of the total number of knots at a given crossing size
  have been constructed by the first two operators. In
  the second step, we apply the remaining two operators to these newly
  constructed knots to obtain the prime alternating knots of crossing size
  $n+1$. The process begins with the prime alternating knot of four
  crossings, the figure eight knot.

  Two of the operators, those that we call $D$ and $ROTS$, are simply
  specific instances of the general splice operation (see Calvo [\calvo]
  for an extensive discussion of the splicing operation). A form of $D$
  was also used by H. de Fraysseix and P. Ossona de Mendez (see
  [\FrayMendez]) in their work to characterize Gauss codes. Another
  one of the operators, $OTS$, made an appearance in Conway's seminal
  paper [\conway]. The fourth operator seems not to have put in an
  appearance in the literature to date, but it is conceptually related to
  Conway's notion of tangle insertion in templates.
 
\par
\section{Knot encodings}
 As remarked above, the DT encoding scheme has featured prominently in
 recent knot enumeration activity. However, our knot operators could not
 be conveniently represented on the data provided by the DT encoding scheme.
 Consequently, we introduce a new knot encoding scheme, which we call a
 group code for the knot configuration. 

 We may regard a knot configuration as a planar graph by considering each
 crossing as a vertex of the graph and the portion of the curve between two
 consecutive crossings as an edge between the two vertices. Such a planar
 graph is always 4-regular.
 
\begin{definition}
 Given the planar graph for a knot configuration, a $k$-group,
 or a group of $k$, is a maximal length sequence
 $c_1,c_2,\ldots,c_k$ of crossings, subject to the
 requirement that $k\ge2$ and for each $i$ from 1 to $k-1$, $c_i$ and
 $c_{i+1}$ are joined by two edges. The two
 edges that are incident to $c_1$ but not incident to $c_2$ are called
 the edges at one end of the group, as are the two edges incident to
 $c_k$ but not to $c_{k-1}$. If $c_1,c_2,\ldots,c_k$ is a group, any
 subsequence of the form $c_i,c_2,\ldots,c_j$, where $1\le i<j\le k$ is called
 a subgroup of the group $c_1,c_2,\ldots,c_k$. More precisely, if $t=j-i+1$, then
 $c_i,c_2,\ldots,c_j$ is called a $t$-subgroup of the $k$-group $c_1,c_2,\ldots,c_k$.
 The two edges incident to $c_i$ but not to $c_{i+1}$ are called the edges
 at one end of the subgroup, as are the two edges that are incidednt to
 $c_j$ but not to $c_{j-1}$.

 Schematically, we shall illustrate a group by

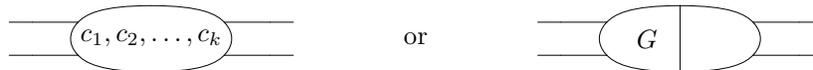
\sarcenter{\xy /r25pt/:,
 (2,.25)="a",
 {\hloop\hloop-}="a",
 (2,-.25)*{\hbox{$c_1,c_2,\ldots,c_k$}},
 (3.1,0)="a2",(3.8,0)="b2","a2";"b2"**\dir{-},
 (3.1,-.5)="a3",(3.8,-.5)="b3","a3";"b3"**\dir{-},
 (10,.25)="q",
 {\hloop\hloop-}="q",
 (9.5,-.25)*{\hbox{$G$}},
 (10,.25)="p1",(10,-.75)="o1","p1";"o1"**\dir{-},
 (6,-.25)*{\hbox{or}},
 (11.1,-.5)="q1",(11.8,-.5)="r1","q1";"r1"**\dir{-},
 (11.1,0)="q2",(11.8,0)="r2","q2";"r2"**\dir{-},
 (.9,0)="a5",(.2,0)="b5","a5";"b5"**\dir{-},
 (.9,-.5)="a4",(.2,-.5)="b4","a4";"b4"**\dir{-},
 (8.9,-.5)="q3",(8.2,-.5)="r3","q3";"r3"**\dir{-},
 (8.9,0)="q4",(8.2,0)="r4","q4";"r4"**\dir{-},
 (4.15,0)="b21","b21";"b2"**\dir{-},
 (4.15,-.5)="b31","b31";"b3"**\dir{-},
 (-.15,0)="b51","b51";"b5"**\dir{-},
 (-.15,-.5)="b41","b41";"b4"**\dir{-},
 (12.15,0)="r21","r21";"r2"**\dir{-},
 (7.85,-.5)="r31","r31";"r3"**\dir{-},
 (12.15,-.5)="r11","r11";"r1"**\dir{-},
 (7.85,0)="r41","r41";"r4"**\dir{-},
\endxy}

\noindent Each group of $k$ is referred to as a positive or negative
 group according to the following convention: if, upon traversing the knot
 in either of the two possible directions of travel, it is found that the
 group is either being entered or exited on both arcs at one end of the
 group, then the group is said to be positive, otherwise the group is said
 to be negative.

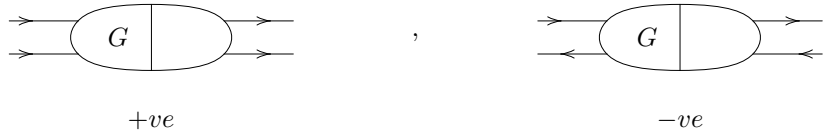
\sarcenter{\xy /r25pt/:,
 (2,.25)="a",
 {\hloop\hloop-}="a",
 (1.5,-.25)*{\hbox{$G$}},
 (2,.25)="a1",(2,-.75)="b1","a1";"b1"**\dir{-},
 (3.1,0)="a2",(3.8,0)="b2","a2";"b2"**\dir{-},"b2"*\dir{>},
 (3.1,-.5)="a3",(3.8,-.5)="b3","a3";"b3"**\dir{-},"b3"*\dir{>},
 (10,.25)="q",
 {\hloop\hloop-}="q",
 (9.5,-.25)*{\hbox{$G$}},
 (10,.25)="p1",(10,-.75)="o1","p1";"o1"**\dir{-},
 (6,-.25)*{\hbox{,}},
 (11.1,-.5)="q1",(11.8,-.5)="r1","q1";"r1"**\dir{-},
 (11.1,0)="q2",(11.8,0)="r2","q2";"r2"**\dir{-},"r2"*\dir{>},
 (.9,0)="a5",(.2,0)="b5","a5";"b5"**\dir{-},
 (.9,-.5)="a4",(.2,-.5)="b4","a4";"b4"**\dir{-},
 (8.9,-.5)="q3",(8.2,-.5)="r3","q3";"r3"**\dir{-},"r3"*\dir{>},
 (8.9,0)="q4",(8.2,0)="r4","q4";"r4"**\dir{-},
 (4.15,0)="b21","b21";"b2"**\dir{-},
 (4.15,-.5)="b31","b31";"b3"**\dir{-},
 (-.15,0)="b51","b51";"b5"**\dir{-},"b5"*\dir{>},
 (-.15,-.5)="b41","b41";"b4"**\dir{-},"b4"*\dir{>},
 (12.15,0)="r21","r21";"r2"**\dir{-},
 (7.85,-.5)="r31","r31";"r3"**\dir{-},
 (12.15,-.5)="r11","r11";"r1"**\dir{-},"r1"*\dir{>},
 (7.85,0)="r41","r41";"r4"**\dir{-},"r4"*\dir{>},
 (2,-1.5)*{\hbox{$+ve$}},
 (10,-1.5)*{\hbox{$-ve$}},
\endxy}
 
 
\noindent Any crossing $c$ that is not an element of a group of $k$ for any
 $k$ shall be called a group of 1, or a loner. The schematic we shall use
 for a loner is

  \sarcenter
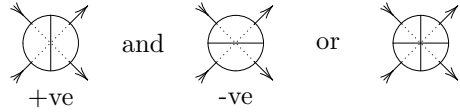
{\xy /r14pt/:,
   (0,-.25)="a11";"a11";{\ellipse(.75){}},
   (2.5,-.25)*{\hbox{and}},
   (5,-.25)="a12";"a12";{\ellipse(.75){}},
   (10,-.25)="a13";"a13";{\ellipse(.75){}},
   (7.5,-.25)*{\hbox{or}},
   (0,-1.75)*{\hbox{+ve}},(5,-1.75)*{\hbox{-ve}},
   (0,-1)="a1",(0,0.5)="b1","a1";"b1"**\dir{-},
   (4.25,-.25)="a2",(5.75,-.25)="b2","a2";"b2"**\dir{-},
   (9.25,-.25)="a3",(10.75,-.25)="b3","a3";"b3"**\dir{-},
   (10,.5)="a4",(10,-1)="b4","a4";"b4"**\dir{-},
   (-.5,-.75)="a5",(.5,.25)="b5",
   (4.5,-.75)="a6",(5.5,.25)="b6",
   (9.5,-.75)="a7",(10.5,.25)="b7",
   "a5";"b5"**\dir{.},
   "a6";"b6"**\dir{.},
   "a7";"b7"**\dir{.},
   (1,.75)="c5","b5";"c5"**\dir{-},"c5"*\dir{>},
   (6,.75)="c6","b6";"c6"**\dir{-},"c6"*\dir{>},
   (11,.75)="c7","b7";"c7"**\dir{-},"c7"*\dir{>},
   (-.5,.25)="a8",(.5,-.75)="b8",
   (4.5,.25)="a9",(5.5,-.75)="b9",
   (9.5,.25)="a10",(10.5,-.75)="b10",
   "a8";"b8"**\dir{.},
   "a9";"b9"**\dir{.},
   "a10";"b10"**\dir{.},
   (1,-1.25)="c8",
   (1,-1.25)="d8","b8";"c8"**\dir{-},"c8"*\dir{>},
   (6,-1.25)="c9","b9";"c9"**\dir{-},"c9"*\dir{>},
   (11,-1.25)="c10","b10";"c10"**\dir{-},"c10"*\dir{>},
   (-.75,.5)="d8",(-1,.75)="e8",
   "e8";"d8"**\dir{-},"d8"*\dir{>},"a8";"d8"**\dir{-},
   (4.25,.5)="d9",(4,.75)="e9",
   "e9";"d9"**\dir{-},"d9"*\dir{>},"a9";"d9"**\dir{-},
   (9.25,.5)="d10",(9,.75)="e10",
   "e10";"d10"**\dir{-},"d10"*\dir{>},"a10";"d10"**\dir{-},
   (-.75,-1)="g8",(-1,-1.25)="f8",
   "f8";"g8"**\dir{-},"g8"*\dir{>},"a5";"g8"**\dir{-},
   (4.25,-1)="g9",(4,-1.25)="f9",
   "f9";"g9"**\dir{-},"g9"*\dir{>},"a6";"g9"**\dir{-},
   (9.25,-1)="g10",(9,-1.25)="f10",
   "f10";"g10"**\dir{-},"g10"*\dir{>},"a7";"g10"**\dir{-},
  \endxy}

\noindent In certain circumstances, we shall refer to a loner as a positive
 group, while in other circumstances the same loner may be considered as a
 negative group. When a loner is being considered as a positive group, then
 the two adjacent outbound edges (as determined by any traversal of the knot)
 are the two edges at one end of the group, and the two inbound
 edges are the two edges at the other and of the group, while if the loner
 is being considered as a negative group, then an inbound edge and the
 adjacent outbound edge are considered to be the two edges at one end of
 the group, with the remaining pair of inbound and outbound edges being the
 two edges at the other end of the group.

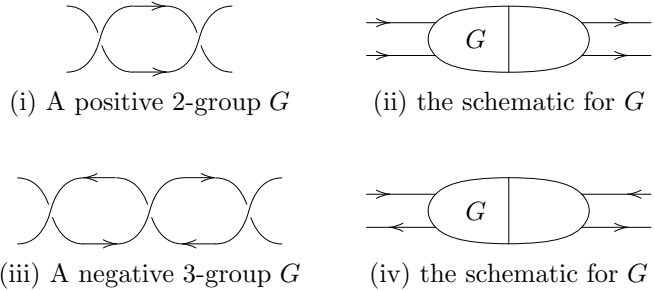
\begin{figure}[ht]
 \centering
 \begin{tabular}{c@{\hskip25pt}c}
  $\vcenter{\xy /r25pt/:,
 (2,.3)="c1",
 {\hcross}="c1",
 (3.5,.3)="c2",
 {\hcross}="c2",
 (3.01,.3)="k1",(3.485,.3)="j1","k1";"j1"**\dir{-},"j1"*\dir{>},
 (3,-.7)="k2",(3.5,-.7)="j2","k2";"j2"**\dir{-},"j2"*\dir{>},
 \endxy}$
 &
  $\vcenter{\xy /r25pt/:,
 (3,.25)="a",
 {\hloop\hloop-}="a",
 (2.5,-.25)*{\hbox{$G$}},
 (3,.25)="a1",(3,-.75)="b1","a1";"b1"**\dir{-},
 (4.1,0)="a2",(4.8,0)="b2","a2";"b2"**\dir{-},"b2"*\dir{>},
 (4.1,-.5)="a3",(4.8,-.5)="b3","a3";"b3"**\dir{-},"b3"*\dir{>},
 (1.9,0)="a5",(1.2,0)="b5","a5";"b5"**\dir{-},
 (1.9,-.5)="a4",(1.2,-.5)="b4","a4";"b4"**\dir{-},
 (5.15,0)="b21","b21";"b2"**\dir{-},
 (5.15,-.5)="b31","b31";"b3"**\dir{-},
 (.85,0)="b51","b51";"b5"**\dir{-},"b5"*\dir{>},
 (.85,-.5)="b41","b41";"b4"**\dir{-},"b4"*\dir{>},
 \endxy}$ \\
 \noalign{\vskip6pt}
 (i) A positive 2-group $G$ & (ii) the schematic for $G$\\
\noalign{\vskip10pt}\\
  $\vcenter{\xy /r25pt/:,
 (2,.3)="c1",
 {\hcross}="c1",
 (3.5,.3)="c2",
 {\hcross}="c2",
 (5,.3)="c3",
 {\hcross}="c3",
 (3.01,.3)="k1",(3.485,.3)="j1","j1";"k1"**\dir{-},"k1"*\dir{>},
 (3,-.7)="k2",(3.5,-.7)="j2","k2";"j2"**\dir{-},"j2"*\dir{>},
 (4.51,.3)="k3",(4.985,.3)="j3","k3";"j3"**\dir{-},"j3"*\dir{>},
 (4.5,-.7)="k4",(5,-.7)="j4","j4";"k4"**\dir{-},"k4"*\dir{>},
 \endxy}$
 &
 $\vcenter{\xy /r25pt/:,
 (3,.25)="q",
 {\hloop\hloop-}="q",
 (2.5,-.25)*{\hbox{$G$}},
 (3,.25)="p1",(3,-.75)="o1","p1";"o1"**\dir{-},
 (4.1,-.5)="q1",(4.8,-.5)="r1","q1";"r1"**\dir{-},"r1"*\dir{>},
 (4.1,0)="q2",(4.8,0)="r2","q2";"r2"**\dir{-},
 (1.9,-.5)="q3",(1.2,-.5)="r3","q3";"r3"**\dir{-},"r3"*\dir{>},
 (1.9,0)="q4",(1.2,0)="r4","q4";"r4"**\dir{-},
 (5.15,0)="r21","r21";"r2"**\dir{-},"r2"*\dir{>},
 (.85,-.5)="r31","r31";"r3"**\dir{-},
 (5.15,-.5)="r11","r11";"r1"**\dir{-},
 (.85,0)="r41","r41";"r4"**\dir{-},"r4"*\dir{>},
 \endxy}$ 
\\
 \noalign{\vskip6pt}
 (iii) A negative 3-group $G$ & (iv) the schematic for $G$\\
\end{tabular} 
\end{figure}

  \noindent Finally, by a group of the knot configuration, we mean any group
 of $k$ for $k\ge1$.
\end{definition}

 It is clear that each crossing of a knot configuration belongs to exactly
 one group of the knot configuration, and so the groups provide us with a
 partition of the set of crossings of the knot configuration. 

 We are now ready to define the group code for a knot configuration.

\begin{definition}\label{group code def}
 Let $C(K)$ denote a configuration of the knot $K$. For each $k\ge1$ for
 which $K$ has at least one group of $k$, let $i_k$ be the number of groups
 of $k$. Label the groups of $k$ in an arbitrary fashion with the labels
 $k_1,k_2,\ldots,k_{i_k}$. Assign an orientation to $K$ and select
 arbitrarily a group to start at. Select an outbound edge at one end of
 this group and, starting with this edge, traverse the knot. Form the
 sequence of group labels in the order that the groups are encountered in
 this traversal. Each label will appear exactly twice in this sequence. If
 a group is negative, then arbitrarily choose one occurrence of the group's
 label in the sequence and preface it with a minus sign. We adopt the
 convention of labelling any loner as a positive group. The resulting
 sequence is called a group code for $C(K)$.
\end{definition}

One consequence of the fact that the labels $k_1,k_2,\ldots,k_{i_k}$ are
assigned arbitrarily is that any cyclic permutation of a group code is
again a group code for the same knot configuration. This is in keeping with
the cyclic structure of a knot, and it is certainly beneficial to view
a group code as a cyclic object rather than a linear array.

  \sarcenter{\xy /r20pt/:,
  {\vloop|(.9)\khole|(.05)>}="a",
  {\vtwistneg\vover-}="b",
  {\xbendr~{"b"+(-.7,3)}{"b"+(0,3)}{"b"+(-2.2,0)}{"b"}}="c","c"+(0,1)="d",
  {\xbendl~{"d"+(0,3)}{"d"+(.7,3)}{"d"}{"d"+(2.2,0)}},
  "d"+(1.07,1.5)="e",
  "b"+(-1.07,1.5)="f",
  {\vcap~{"f"+(0,2.3)}{"e"+(0,2.3)}{"f"}{"e"}|(.4)\khole}?(.4)="x",
  "x"+(-.9 ,0)*!<3pt,-4pt>{\hbox{$\ssize 2$}},
  "x"+(.15,0)*!<-2pt,-4pt>{\hbox{$\ssize 1$}},
  "x"+(-.4,-1)*!<1pt,-3pt>{\hbox{$\ssize 3$}},
  "x"+(-.4,-2.5)*!<1pt,-6pt>{\hbox{$\ssize 4$}},
  (9,1)*{\hbox{\small$+$}},
  (9.5,1)*{\hbox{\small$-$}},(10,1)*{\hbox{\small$+$}},(10.5,1)*{\hbox{\small$-$}},(11,1)*{\hbox{\small$+$}},(11.5,1)*{\hbox{\small$-$}},(12,1)*{\hbox{\small$+$}},(12.5,1)*{\hbox{\small$-$}},
  (7,.5)*{\hbox{\small Gauss code:}},
  (9,.5)*{\hbox{\small$1,$}},
  (9.5,.5)*{\hbox{\small$2,$}},(10,.5)*{\hbox{\small$3,$}},(10.5,.5)*{\hbox{\small$4,$}},(11,.5)*{\hbox{\small$2,$}},(11.5,.5)*{\hbox{\small$1,$}},(12,.5)*{\hbox{\small$4,$}},(12.5,.5)*{\hbox{\small$3,$}},
  (7,-1.7)*{\hbox{\small Group code:}},
  (9,-.5)*{\hbox{\small$1,$}},
  (9.5,-.5)*{\hbox{\small$2,$}},(10,-.5)*{\hbox{\small$3,$}},(10.5,-.5)*{\hbox{\small$4,$}},(11,-.5)*{\hbox{\small$2,$}},(11.5,-.5)*{\hbox{\small$1,$}},(12,-.5)*{\hbox{\small$4,$}},(12.5,-.5)*{\hbox{\small$3,$}},
  (8.8,-.75)="a1",(8.8,-.95)="b1","a1";"b1"**\dir{-},
  (9.65,-.75)="a2",(9.65,-.95)="b2","a2";"b2"**\dir{-},
  (9.25,-.95)="a3",(9.25,-1.25)="b3","a3";"b3"**\dir{-},"b3"*\dir{>},
  (9.85,-.75)="a4",(9.85,-.95)="b4","a4";"b4"**\dir{-},
  (10.65,-.75)="a5",(10.65,-.95)="b5","a5";"b5"**\dir{-},
  (10.25,-.95)="a6",(10.25,-1.25)="b6","a6";"b6"**\dir{-},"b6"*\dir{>},
  (10.85,-.75)="a7",(10.85,-.95)="b7","a7";"b7"**\dir{-},
  (11.65,-.75)="a8",(11.65,-.95)="b8","a8";"b8"**\dir{-},
  (11.25,-.95)="a9",(11.25,-1.25)="b9","a9";"b9"**\dir{-},"b9"*\dir{>},
  (11.85,-.75)="a10",(11.85,-.95)="b10","a10";"b10"**\dir{-},
  (12.65,-.75)="a11",(12.65,-.95)="b11","a11";"b11"**\dir{-},
  (12.25,-.95)="a12",(12.25,-1.25)="b12","a12";"b12"**\dir{-},"b12"*\dir{>},
  "b1";"b2"**\dir{-},
  "b4";"b5"**\dir{-},
  "b7";"b8"**\dir{-},"b10";"b11"**\dir{-},
  (9.25,-1.7)*{\hbox{\small$2_1,$}},
  (10.25,-1.7)*{\hbox{\small$2_2,$}},(11.25,-1.7)*{\hbox{\small$-2_1,$}},(12.25,-1.7)*{\hbox{\small$-2_2,$}},
 \endxy}

\par
\section{The Orbit of a Group}

 We shall need to have a graph-theoretic description
 of an $m$-tangle. A non-trivial $m$-tangle contains at least one crossing,
 although it may consist of an $m-2$ tangle together with an edge that is
 not incident to any crossing in the tangle. We shall refer to an edge that
 passes through a tangle without being incident to any vertex of the tangle
 as a pass-through edge for the tangle. We shall be largely interested
 in $4$-tangles, (in fact, the word tangle shall mean $4$-tangle), and a
 $4$-tangle that has an edge that is not incident to any vertex of the tangle
 is either trivial, or else consists of a $2$-tangle and this edge. In a
 prime knot, there are no $2$-tangles, so the smallest $m$ for which there
 exists an $m$-tangle is $m=4$. Thus in a prime knot, a non-trivial
 $4$-tangle must have each edge that enters the $4$-tangle be incident to a
 vertex of the tangle. Moreover, since a prime knot has no $2$-tangles,
 each non-trivial $4$-tangle of a prime knot configuration is a connected
 graph. In the study of prime knots, we shall only have need of connected
 $m$-tangles, and a connected $m$-tangle with at least
 one crossing can be characterized as a connected induced subgraph of
 the planar graph, together with all edges incident to exactly one vertex
 in the vertex set of the subgraph (by induced subgraph, we mean a subgraph
 formed by choosing a set of vertices and then taking all edges that have
 both endpoints belonging to the selected set of vertices). An edge that is
 incident to a vertex of an $m$-tangle $T$ but that is not an edge of $T$ is
 said to be incident to $T$. Thus an $m$-tangle has $m$ incident edges.

 The following two facts about tangles in a prime knot will prove to be quite
 useful.

\begin{lemma}\label{tangle intersection}
 Let $C(K)$ be a configuration of a prime knot $K$, and let $T_1$ and
 $T_2$ be tangles in $C(K)$ such that $T_1\cap T_2\ne\emptyset$, $T_1-T_2
 \ne\emptyset$ and $T_2-T_1\ne \emptyset$. Then there are at least
 two arcs of $T_1$ each having exactly one endpoint in $T_2$, and at least
 two arcs of $T_2$ each having exactly one endpoint in $T_1$.
\end{lemma}

\begin{proof}
 Since $T_2$ is connected, there must be at least one arc with one
 endpoint in $T_2-T_1$ and the other endpoint in $T_2\cap T_1$. Suppose
 there is exactly one such arc. Then the other 3 arcs
 incident to $T_1$ must have their other endpoints outside of $T_2$.
 We identify four cases.

 \noindent Case 1: none of these three arcs have endpoints in $T_1\cap T_2$.
 Then since $T_1\cap T_2$ must have at least 4 incident arcs, there must
 be at least 3 arcs joining vertices in $T_1\cap T_2$ to vertices in
 $T_1-T_2$, whence three of the four arcs incident to $T_2$ have already
 been accounted for. Now there must be at least four arcs incident to
 $T_2-T_1$, and only one of these can have its other endpoint in $T_1\cap
 T_2$. Thus there are at least three arcs incident to $T_2-T_1$ whose
 other endpoints are not in $T_1$. But then $T_2$ has at least six
 incident arcs, which is not possible.

 \noindent Case 2: exactly one of these three arcs has an endpoint in
 $T_1\cap T_2$. Then the remaining two arcs are incident to $T_1-T_2$,
 with their other endpoint not in $T_2$.
 Since $T_1-T_2$ must have at least four incident arcs, there must be
 at least two arcs joining vertices in $T_1-T_2$ to vertices in $T_1\cap
 T_2$. But then we have accounted for at least three arcs incident to
 $T_2$, and all three are actually incident to $T_1\cap T_2$. Since
 $T_2-T_1$ must have at least four incident arcs, of which exactly one
 is incident to $T_1\cap T_2$, there must be at least three arcs joining
 vertices in $T_2-T_1$ to vertices neither in $T_2$ nor in $T_1$. But then
 we have accounted for at least six arcs incident to $T_2$, which is not
 possible.

 \noindent Case 3: exactly two of these three arcs have an endpoint in
 $T_1\cap T_2$. Then the remaining arc is incident to $T_1-T_2$ with
 its other endpoint not in $T_2$. Since $T_1-T_2$ must have at least
 four incident arcs, there must be at least three arcs joining vertices
 in $T_1-T_2$ to vertices in $T_1\cap T_2$. We have now accounted for
 five arcs incident to $T_2$, which is not possible.

 \noindent Case 4: all three of the arcs are actually incident to
 $T_1\cap T_2$. Since there are at least four arcs incident to $T_1-T_2$,
 none of which may have their other endpoints outside of $T_1$, there
 must be at least four arcs joining vertices in $T_1-T_2$ to vertices
 of $T_1\cap T_2$. But then we have accounted for at least seven arcs
 incident to $T_2$, which is not possible.

 Since every case leads to a contradiction, we conclude that there must
 be at least two arcs having one endpoint in $T_2-T_1$ and the other
 endpoint in $T_2\cap T_1$. Since $T_1$ and $T_2$ have symmetric roles,
 it follows as well that there must be at least two arcs having one
 endpoint in $T_1-T_2$ and the other endpoint in $T_2\cap T_1$.
\end{proof}

\begin{proposition}\label{tangle replacement}
 Let $C(K)$ be a configuration of an alternating prime knot, and let
 $T_1$ be a non-trivial tangle of $K$. Cut the four incident
 arcs to $T_1$ to allow $T_1$ to be removed, and replace $T_1$ by
 a non-trivial tangle $T_2$ such that the result is a configuration
 $C(K_1)$ of an alternating knot $K_1$. If $T_2$ does not contain
 any 2-tangle, then $K_1$ is prime.
\end{proposition}

\begin{proof}
  Suppose that $K_1$ is not prime, and let $T$ be a 2-tangle of
  $C(K_1)$. If
  $T\cap T_2=\emptyset$, then $T$ is a 2-tangle in $C(K)$, which is
  not possible. Thus $T\cap T_2\ne\emptyset$. If $T\subseteq T_2$,
  then we are done, so suppose that $T-T_2\ne \emptyset$. If $T_2
  \subseteq T$, then we may undo the surgery to form $T'$ from $T$
  by replacing $T_2$ with $T_1$. But then $T'$ is a 2-tangle in $C(K)$,
  which is not possible. Thus $T_2-T\ne\emptyset$.

  There are at most two arcs incident to $T\cap T_2$ that are
  actually incident to $T$. There are three cases to consider.

  \noindent Case 1: there are exactly two arcs incident to $T\cap T_2$
  that are actually incident to $T$. Since $T-T_2$ is an $m$-tangle in
  $C(K)$, $m\ge 4$, there are at least four edges incident to $T-T_2$
  each with their other endpoint in $T\cap T_2$. But then they account
  for all four edges incident to $T_2$, in which case $T_2-T$ has
  no incident edges. This is not possible.

  \noindent Case 2: there is exactly one arc incident to $T\cap T_2$
  that is incident to $T$. Then there is one arc incident to $T$ that
  is actually incident to $T-T_2$. Since $T-T_2$ is an $m$-tangle in
  $C(K)$, $m\ge4$, there are at least four edges incident to $T-T_2$,
  exactly one of which has an endpoint outside $T$, whence there are
  at least three edges incident to $T-T_2$ each with their other
  endpoint in $T\cap T_2$. But this accounts for at least three edges
  incident to $T_2$. Thus there is at most one edge incident to $T_2$
  that has an endpoint in $T_2-T$, which means that $T_2-T$ is at
  most a 2-tangle. But then $T_2-T$ must be a 2-tangle, and it is
  contained in $T_2$, as required.

  \noindent Case 3: there are no arcs incident to $T\cap T_2$ that
  are incident to $T$. Since $T-T_2$ is an $m$-tangle in $C(K)$,
  $m\ge4$, and the two arcs incident to $T$ are actually incident
  to $T-T_2$, there must be at least two arcs from $T-T_2$ to $T\cap
  T_2$. But this accounts for at least two of the arcs incident to
  $T_2$, whence there are at most two arcs incident to $T_2$ that
  are actually incident to $T_2-T$. Any additional arcs incident to
  $T_2-T$ would have to have their endpoints in $T\cap T_2$, but these
  would then be arcs incident to $T\cap T_2$ and incident to $T$, of
  which there are none. Thus $T_2-T$ is a 2-tangle contained in $T_2$,
  as required.

  This completes the proof that if $K_1$ is not prime, then $T_2$ contains
  a 2-tangle.
\end{proof}

We remark that in the setting of the preceding proposition, it is
possible that certain tangles $T_2$ might lead to a link or a
non-alternating knot. We are only interested in those situations
as described in the conditions of the proposition.
 
\begin{proposition}\label{prep for orbit}
  Let $C(K)$ be a configuration of a prime knot $K$, and let $e$ and $f$
  be arcs of $C(K)$. Let $a$ and $b$ be endpoints of $e$ and $f$,
  respectively. If $T_1$ and $T_2$ are tangles containing $a$ and $b$ with
  $e$ and $f$ being arcs incident to both $T_1$ and $T_2$, then either
  $T_1\subseteq T_2$ or else $T_2\subseteq T_1$.
\end{proposition} 

\begin{proof}
  Suppose not. Then $T_1-T_2$ is an $m$-tangle with $m\ge4$ since
  $K$ is prime. Thus there are at least four arcs incident to
  $T_1 - T_2$. Since $a, b\in T_2$, these arcs are all different from $e$
  and $f$. But $T_1$ has four incident arcs, two of which are $e$ and $f$,
  so at most two of these $m$ arcs are incident to $T_1$, with the having
  one end point in $T_1-T_2$ and the other endpoint in $T_1\cap T_2$. Thus
  at least two of these $m$ arcs are incident to $T_2$. Since $T_2$ is a
  4-tangle, and $e$ and $f$ are incident to $T_2$, it must be that there
  are exactly two arcs that have one endpoint in $T_1\cap T_2$, and the
  other endpoint in $T_1 - T_2$. The same reasoning applied to $T_2 - T_1$
  establishes that there are exactly two arcs that have an endpoint in
  $T_1\cap T_2$ and the other endpoint in $T_2 - T_1$. But then the four
  arcs incident to $T_1$ are $e$, $f$ and the two arcs going from
  $T_1\cap T_2$ to $T_2 - T_1$, whence $T_1 - T_2$ has only two incident
  arcs. This contradicts the fact that $m\ge4$, so the assumption that
  neither $T_1\subseteq T_2$ nor $T_2\subseteq T_1$ holds is false.
\end{proof}

\begin{corollary}\label{prep for orbit 2}
  If $e$ and $f$ are edges in a configuration $C(K)$ of a prime knot $K$,
  and $a$ and $b$ are endpoints of $e$ and $f$, respectively, such that
  there is a tangle $T$ with $e$ and $f$ incident to $T$ and $a, b\in T$,
  then there is a minimum such tangle.
\end{corollary} 

\begin{proof}
 Of all tangles with this property, let $T$ denote one with fewest vertices.
 By Proposition \ref{prep for orbit}, if $T_1$ is any tangle with this
 property, then either $T\subseteq T_1$ or else $T_1\subseteq T$. If
 $T_1\subseteq T$, then by the minimality of $T$, we have $T_1 = T$.
 Otherwise, $T\subseteq T_1$.
\end{proof}

 We are now ready to present the definition of the orbit (or flype
 circuit) of a group in a prime knot. 

\begin{definition}
  Let $C(K)$ be a configuration of the prime alternating knot $K$, and
  let $G$ be a group of $C(K)$. If $K$ is a torus knot, then $C(K)$ is

  \sarcenter
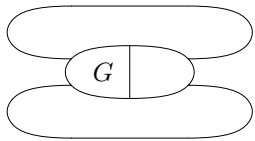
{\xy /r20pt/:,
    (2,.25)="a",
    {\hloop\hloop-}="a",
    (8,-.25)="a12",
    (1.5,-.25)*{\hbox{$G$}},
    (2,.25)="a1",(2,-.75)="b1","a1";"b1"**\dir{-},
    (3.1,0)="a2",
    (3.1,-.5)="a3",
    (.9,1)="c1",
    {\hloop-}="c1",
    (.9,-.5)="c2",
    {\hloop-}="c2",
    (3.1,1)="c3",
    {\hloop}="c3",
    (3.1,-.5)="c4",
    {\hloop}="c4",
    (.9,-1.5)="a6",(3.1,-1.5)="b6","a6";"b6"**\dir{-},
    (.9,1)="a7",(3.1,1)="b7","a7";"b7"**\dir{-},
   \endxy}

\noindent and we define the orbit of $G$ to be the sequence $(G)$.
Otherwise, $K$ is not a torus knot, and $C(K)$ is

  \sarcenter
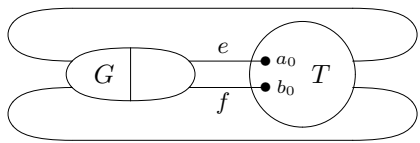
{\xy /r20pt/:,
    (5.25,-.25)="a11";"a11";{\ellipse(1){}},
    (5.6,-.25)*{\hbox{$T$}},
    (2,.25)="a",
    {\hloop\hloop-}="a",
    (1.5,-.25)*{\hbox{$G$}},
    (2,.25)="a1",(2,-.75)="b1","a1";"b1"**\dir{-},
    (3.1,0)="a2",(4.55,0)="b2","a2";"b2"**\dir{-},
    (3.1,-.5)="a3",(4.55,-.5)="b3","a3";"b3"**\dir{-},
    (.9,-1.5)="a6",(6.2,-1.5)="b6","a6";"b6"**\dir{-},
    (.9,1)="a7",(6.2,1)="b7","a7";"b7"**\dir{-},
    (.9,1)="c1",
    {\hloop-}="c1",
    (.9,-.5)="c2",
    {\hloop-}="c2",
    (6.2,1)="c3",
    {\hloop}="c3",
    (6.2,-.5)="c4",
    {\hloop}="c4",
    (4.55,0)*{\hbox{\small$\bullet$}},(4.55,-.5)*{\hbox{\small$\bullet$}},
    (4.95,0)*{\hbox{$\ssize a_0$}},(4.95,-.5)*{\hbox{$\ssize b_0$}},
    (3.75,.25)*{\hbox{\small$e$}},(3.75,-.8)*{\hbox{\small$f$}},
   \endxy}

\noindent where $T$ is a non-trivial tangle. Let $\stackrel{*}{G_0} = G_0
 = G$. By Corollary \ref{prep for orbit 2}, there is a tangle $T_0$ which
 contains $a_0$, $b_0$ and has $e$, $f$ as incident edges, and which is
 contained in all other tangles with this property. In particular,
 $T_0\subseteq T$ and $T$ is of the form

  \sarcenter
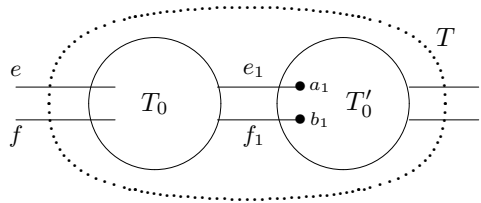
{\xy /r25pt/:,
    (7.1,-.25)*{\hbox{$T_0$}},
    (7.1,-.25)="q";"q";{\ellipse(1){}},
    (5,0)="a4",(6.5,0)="b4","a4";"b4"**\dir{-},
    (5,-.5)="a5",(6.5,-.5)="b5","a5";"b5"**\dir{-},
    (9.95,-.25)="a13";"a13";{\ellipse(1){}},
    (10.2,-.25)*{\hbox{$T'_0$}},
    (8.05,0)="a8",(9.3,0)="b8","a8";"b8"**\dir{-},
    (8.05,-.5)="a9",(9.3,-.5)="b9","a9";"b9"**\dir{-},
    (10.95,0)="l1",(12,0)="p1","l1";"p1"**\dir{-},
    (10.95,-.5)="l2",(12,-.5)="p2","l2";"p2"**\dir{-},
    (5.5,-.25)="t0",(11.5,-.25)="t1",
    "t0";"t1"**\crv{~*=<\jot>{.}(5.5,.5) & (6.25,1.25) & (10.75,1.25)
      & (11.5,.5)},
    "t1";"t0"**\crv{~*=<\jot>{.}(11.5,-1) & (10.75,-1.75) & (6.25,-1.75)
      & (5.5,-1)},
    (9.3,0)*{\hbox{\small$\bullet$}},(9.3,-.5)*{\hbox{\small$\bullet$}},
    (5,.25)*{\hbox{\small$e$}},(5,-.75)*{\hbox{\small$f$}},
    (8.6,.25)*{\hbox{\small$e_1$}},(8.6,-.75)*{\hbox{\small$f_1$}},
    (11.5,.75)*{\hbox{$T$}},
    (9.6,0)*{\hbox{$\ssize a_1$}},(9.6,-.5)*{\hbox{$\ssize b_1$}},
   \endxy}

 \noindent If $T'_0$ is trivial, then the orbit of $G$ is the sequence
 $(\stackrel{*}{G_0},T_0)$. Otherwise, $T'_0$ is non-trivial and we examine
 the arcs $e_1$ and $f_1$.
 
 Case (1): $a_1 = b_1$. Let $G_1$ denote the group containing this crossing. 
  Then $G_1\subseteq T'_0$, and $T'_0$ has the form

  \sarcenter
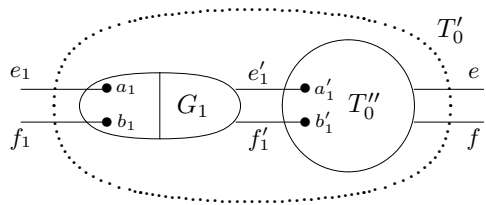
{\xy /r25pt/:,
    (7.6,-.25)*{\hbox{$G_1$}},
    (7.1,.25)="q",
    {\hloop\hloop-}="q",
    (7.1,.25)="q1",(7.1,-.75)="r1","q1";"r1"**\dir{-},
    (5,0)="a4",(6.3,0)="b4","a4";"b4"**\dir{-},
    (5,-.5)="a5",(6.3,-.5)="b5","a5";"b5"**\dir{-},
    (9.95,-.25)="a13";"a13";{\ellipse(1){}},
    (10.2,-.25)*{\hbox{$T''_0$}},
    (8.25,0)="a8",(9.3,0)="b8","a8";"b8"**\dir{-},
    (8.25,-.5)="a9",(9.3,-.5)="b9","a9";"b9"**\dir{-},
    (10.95,0)="l1",(12,0)="p1","l1";"p1"**\dir{-},
    (10.95,-.5)="l2",(12,-.5)="p2","l2";"p2"**\dir{-},
    (5.5,-.25)="t0",(11.5,-.25)="t1",
    "t0";"t1"**\crv{~*=<\jot>{.}(5.5,.5) & (6.25,1.25) & (10.75,1.25)
      & (11.5,.5)},
    "t1";"t0"**\crv{~*=<\jot>{.}(11.5,-1) & (10.75,-1.75) & (6.25,-1.75)
      & (5.5,-1)},
    (6.3,0)*{\hbox{\small$\bullet$}},(6.3,-.5)*{\hbox{\small$\bullet$}},
    (9.3,0)*{\hbox{\small$\bullet$}},(9.3,-.5)*{\hbox{\small$\bullet$}},
    (5,.25)*{\hbox{\small$e_1$}},(5,-.75)*{\hbox{\small$f_1$}},
    (8.6,.25)*{\hbox{\small$e'_1$}},(8.6,-.75)*{\hbox{\small$f'_1$}},
    (11.85,.25)*{\hbox{\small$e$}},(11.85,-.75)*{\hbox{\small$f$}},
    (11.5,.9)*{\hbox{$T'_0$}},
    (6.6,0)*{\hbox{$\ssize a_1$}},(6.6,-.5)*{\hbox{$\ssize b_1$}},
    (9.6,0)*{\hbox{$\ssize a'_1$}},(9.6,-.5)*{\hbox{$\ssize b'_1$}},
   \endxy}

  \noindent where $T''_0$ is necessarily a non-trivial tangle.  By Corollary
  \ref{prep for orbit 2}, there is a tangle $T_1$ which contains $a'_1$,
  $b'_1$ and has $e'_1$, $f'_1$ as incident edges, and which is contained
  in all other tangles with this property. In particular, $T_1\subseteq
  T''_0$ and $T''_0$ is of the form

   \sarcenter
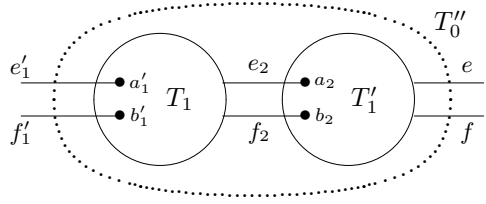
{\xy /r25pt/:,
     (7.4,-.25)*{\hbox{$T_1$}},
     (7.1,-.25)="q";"q";{\ellipse(1){}},
     (5,0)="a4",(6.5,0)="b4","a4";"b4"**\dir{-},
     (5,-.5)="a5",(6.5,-.5)="b5","a5";"b5"**\dir{-},
     (9.95,-.25)="a13";"a13";{\ellipse(1){}},
     (10.2,-.25)*{\hbox{$T'_1$}},
     (8.05,0)="a8",(9.3,0)="b8","a8";"b8"**\dir{-},
     (8.05,-.5)="a9",(9.3,-.5)="b9","a9";"b9"**\dir{-},
     (10.95,0)="l1",(12,0)="p1","l1";"p1"**\dir{-},
     (10.95,-.5)="l2",(12,-.5)="p2","l2";"p2"**\dir{-},
     (5.5,-.25)="t0",(11.5,-.25)="t1",
     "t0";"t1"**\crv{~*=<\jot>{.}(5.5,.5) & (6.25,1.25) & (10.75,1.25) & (11.5,.5)},
     "t1";"t0"**\crv{~*=<\jot>{.}(11.5,-1) & (10.75,-1.75) & (6.25,-1.75) & (5.5,-1)},
     (6.5,0)*{\hbox{\small$\bullet$}},(6.5,-.5)*{\hbox{\small$\bullet$}},
     (9.3,0)*{\hbox{\small$\bullet$}},(9.3,-.5)*{\hbox{\small$\bullet$}},
     (5,.25)*{\hbox{\small$e'_1$}},(5,-.75)*{\hbox{\small$f'_1$}},
     (8.6,.25)*{\hbox{\small$e_2$}},(8.6,-.75)*{\hbox{\small$f_2$}},
     (11.75,.25)*{\hbox{\small$e$}},(11.75,-.75)*{\hbox{\small$f$}},
     (11.5,.9)*{\hbox{$T''_0$}},
     (6.8,0)*{\hbox{$\ssize a'_1$}},(6.8,-.5)*{\hbox{$\ssize b'_1$}},
     (9.6,0)*{\hbox{$\ssize a_2$}},(9.6,-.5)*{\hbox{$\ssize b_2$}},
    \endxy}

  \noindent Let $\stackrel{*}{G_1} = G_1$.

  Case (2): $a_1\ne b_1$.  By Corollary \ref{prep for orbit 2}, there is a
  tangle $T_1$ which contains $a_1$, $b_1$ and has $e_1$, $f_1$ as incident
  edges, and which is contained in all other tangles with this property. In
  particular, $T_1\subseteq T'_0$ and $T'_0$ is of the form

    \sarcenter
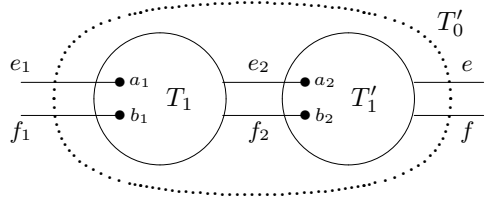
{\xy /r25pt/:,
      (7.4,-.25)*{\hbox{$T_1$}},
      (7.1,-.25)="q";"q";{\ellipse(1){}},
      (5,0)="a4",(6.5,0)="b4","a4";"b4"**\dir{-},
      (5,-.5)="a5",(6.5,-.5)="b5","a5";"b5"**\dir{-},
      (9.95,-.25)="a13";"a13";{\ellipse(1){}},
      (10.2,-.25)*{\hbox{$T'_1$}},
      (8.05,0)="a8",(9.3,0)="b8","a8";"b8"**\dir{-},
      (8.05,-.5)="a9",(9.3,-.5)="b9","a9";"b9"**\dir{-},
      (10.95,0)="l1",(12,0)="p1","l1";"p1"**\dir{-},
      (10.95,-.5)="l2",(12,-.5)="p2","l2";"p2"**\dir{-},
      (5.5,-.25)="t0",(11.5,-.25)="t1",
      "t0";"t1"**\crv{~*=<\jot>{.}(5.5,.5) & (6.25,1.25)
          & (10.75,1.25) & (11.5,.5)},
      "t1";"t0"**\crv{~*=<\jot>{.}(11.5,-1) & (10.75,-1.75)
          & (6.25,-1.75) & (5.5,-1)},
      (6.5,0)*{\hbox{\small$\bullet$}},(6.5,-.5)*{\hbox{\small$\bullet$}},
      (9.3,0)*{\hbox{\small$\bullet$}},(9.3,-.5)*{\hbox{\small$\bullet$}},
      (5,.25)*{\hbox{\small$e_1$}},(5,-.75)*{\hbox{\small$f_1$}},
      (8.6,.25)*{\hbox{\small$e_2$}},(8.6,-.75)*{\hbox{\small$f_2$}},
      (11.75,.25)*{\hbox{\small$e$}},(11.75,-.75)*{\hbox{\small$f$}},
      (11.5,.9)*{\hbox{$T'_0$}},
      (6.8,0)*{\hbox{$\ssize a_1$}},(6.8,-.5)*{\hbox{$\ssize b_1$}},
      (9.6,0)*{\hbox{$\ssize a_2$}},(9.6,-.5)*{\hbox{$\ssize b_2$}},
     \endxy}

  \noindent Let $\stackrel{*}{G_1} = \{ e_1, f_1 \}$. In either case, if
  $T'_1$ is trivial, we define the orbit of $G$ to be the sequence
  $(\stackrel{*}{G_0}, T_0, \stackrel{*}{G_1}, T_1)$. If $T'_1$ is
  non-trivial, we repeat the process. For some $k$, we end up with $T'_k$
  trivial and the process stops. The orbit of $G$ is the sequence
  $(\stackrel{*}{G_0}, T_0, \stackrel{*}{G_1}, T_1,\ldots,\stackrel{*}{G_k},
  T_k)$.

  \noindent $C(K)$ has the form

    \sarcenter
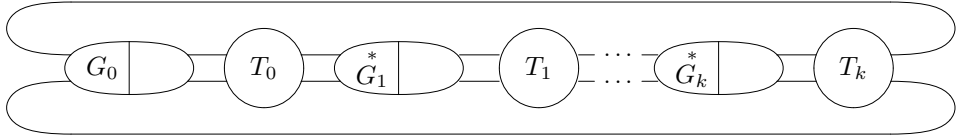
{\xy /r20pt/:,
   (4.55,-.25)="a11",
   ;"a11";{\ellipse(.75){}},
   (4.55,-.25)*{\hbox{$T_0$}},
   (2,.25)="a",
   {\hloop\hloop-}="a",
   (6.6,-.25)*{\hbox{$\stackrel{*}{G_1}$}},
   (7.1,.25)="q",
   {\hloop\hloop-}="q",
   (7.1,.25)="q1",(7.1,-.75)="r1","q1";"r1"**\dir{-},
   (1.5,-.25)*{\hbox{$G_0$}},
   (2,.25)="a1",(2,-.75)="b1","a1";"b1"**\dir{-},
   (3.1,0)="a2",(3.85,0)="b2","a2";"b2"**\dir{-},
   (3.1,-.5)="a3",(3.85,-.5)="b3","a3";"b3"**\dir{-},
   (5.25,0)="a4",(6,0)="b4","a4";"b4"**\dir{-},
   (5.25,-.5)="a5",(6,-.5)="b5","a5";"b5"**\dir{-},
   (.9,-1.5)="a6",(16.4,-1.5)="b6","a6";"b6"**\dir{-},
   (.9,1)="a7",(16.4,1)="b7","a7";"b7"**\dir{-},
   (.9,1)="c1",
   {\hloop-}="c1",
   (.9,-.5)="c2",
   {\hloop-}="c2",
   (9.75,-.25)="a13",
   ;"a13";{\ellipse(.75){}},
   (9.75,-.25)*{\hbox{$T_1$}},
   (8.25,0)="a8",(9,0)="b8","a8";"b8"**\dir{-},
   (8.25,-.5)="a9",(9,-.5)="b9","a9";"b9"**\dir{-},
   (16.4,1)="c3",
   {\hloop}="c3",
   (16.4,-.5)="c4",
   {\hloop}="c4",
   (15.7,-.25)="a14",
   ;"a14";{\ellipse(.75){}},
   (15.7,-.25)*{\hbox{$T_k$}},
   (12.65,-.25)*{\hbox{$\stackrel{*}{G_k}$}},
   (13.15,.25)="u",
   {\hloop\hloop-}="u",
   (13.15,.25)="u1",(13.15,-.75)="v1","u1";"v1"**\dir{-},
   (14.25,0)="a10",(15,0)="b10","a10";"b10"**\dir{-},
   (14.25,-.5)="a11",(15,-.5)="b11","a11";"b11"**\dir{-},
   (11.25,0)*{\hbox{\small$\ldots$}},(11.25,-.5)*{\hbox{\small$\ldots$}},
   (10.5,0)="l1",(10.8,0)="p1","l1";"p1"**\dir{-},
   (10.5,-.5)="l2",(10.8,-.5)="p2","l2";"p2"**\dir{-},
   (11.7,0)="l3",(12,0)="p3","l3";"p3"**\dir{-},
   (11.7,-.5)="l4",(12,-.5)="p4","l4";"p4"**\dir{-},
  \endxy}

  \noindent The tangles $T_0,T_1,\ldots,T_k$ are called the {\em min-tangles} of
   the orbit of $G$, and for each $i$, $0\le i\le k$, $\stackrel{*}{G_i}$
   marks {\em position} $i$ of the orbit. For each min-tangle $T_i$, the
   orbit causes the four arcs incident to $T_i$ to be partitioned into two
   pairs, and each pair shall be referred to as a {\em position pair} for the min-tangle.
   For each $i$, $\stackrel{*}{G_i}$ either
   denotes a group or else a pair of arcs. Let $G_F =
   \{\stackrel{*}{G_i}|\stackrel{*}{G_i}$ is a group\}. $G_F$ is called the
   {\em full group determined by} $G$ and each $\stackrel{*}{G_i}$ that is a group
   is said to be the {\em component in position} $i$ of the full group. If
   $G_F = \{G\}$, then $G_F$ is said to be a {\em full group}. Otherwise,
   we say that $G_F$ is a {\em split group}. If $G_F=\{G\}$, we also say
   that $G$ is a full group, while if $G_F$ is not a full group, then we
   say that $G$ is a member of the split group.
\end{definition}

   We remark that if the two arcs on the other side of group $G$ had been
   chosen to construct the orbit of $G$, the sequence that would result
   would be the present sequence presented in the reverse order. Moreover,
   if for some $i$, $G_i^*$ is a group and we used the two arcs from $G_i^*$
   to $T_i$ to construct the orbit of the group $G_i^*$, the result would
   be the sequence $(G_i^*,T_i,G_{i+1}^*,T_{i+1},\ldots,G_k^*,T_k,G_0^*,T_0,
   \ldots,G_{i-1}^*,T_{i-1})$.

\begin{definition}
   A knot configuration is said to be a split-group configuration
   if there exist groups $G_1$ and $G_2$ and non-trivial tangles $T_1$ and
   $T_2$ such that $C(K)$ is of the form

   \sarcenter
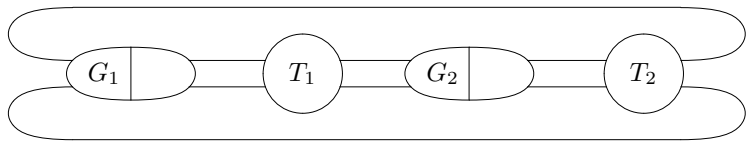
{\xy /r20pt/:,
   (5.25,-.25)="a11",
   ;"a11";{\ellipse(.75){}},
   (5.25,-.25)*{\hbox{$T_1$}},
   (2,.25)="a",
   {\hloop\hloop-}="a",
   (7.9,-.25)*{\hbox{$G_2$}},
   (8.4,.25)="q",
   {\hloop\hloop-}="q",
   (8.4,.25)="q1",(8.4,-.75)="r1","q1";"r1"**\dir{-},
   (1.5,-.25)*{\hbox{$G_1$}},
   (2,.25)="a1",(2,-.75)="b1","a1";"b1"**\dir{-},
   (3.1,0)="a2",(4.55,0)="b2","a2";"b2"**\dir{-},
   (3.1,-.5)="a3",(4.55,-.5)="b3","a3";"b3"**\dir{-},
   (5.95,0)="a4",(7.3,0)="b4","a4";"b4"**\dir{-},
   (5.95,-.5)="a5",(7.3,-.5)="b5","a5";"b5"**\dir{-},
   (.9,-1.5)="a6",(12.4,-1.5)="b6","a6";"b6"**\dir{-},
   (.9,1)="a7",(12.4,1)="b7","a7";"b7"**\dir{-},
   (.9,1)="c1",
   {\hloop-}="c1",
   (.9,-.5)="c2",
   {\hloop-}="c2",
   (11.7,-.25)="a13",
   ;"a13";{\ellipse(.75){}},
   (11.7,-.25)*{\hbox{$T_2$}},
   (9.5,0)="a8",(11,0)="b8","a8";"b8"**\dir{-},
   (9.5,-.5)="a9",(11,-.5)="b9","a9";"b9"**\dir{-},
   (12.4,1)="c3",
   {\hloop}="c3",
   (12.4,-.5)="c4",
   {\hloop}="c4",
  \endxy}
   
  \noindent In this case, we say that $G_1$ and $G_2$ are part of a split
  group, with the splitting achieved by tangles $T_1$ and $T_2$. 
  A knot configuration that is not a split-group configuration is said to
  be a full-group configuration.
\end{definition}  

 Observe that if $G_1$ and $G_2$ are part of a split group, with the
 splitting achieved by tangles $T_1$ and $T_2$, and if $G_1$ is a group of
 $k$ and $G_2$ is a group of $l$, then a sequence of $l$ flype moves can be
 applied to flype each crossing of $G_2$ across $T_1$. In the resulting
 configuration, the groups $G_1$ and $G_2$ have become one group $G$ of
 $k$ + $l$ such that $\flype G,T_1,T_2,$ is a flype scenario. We say that
 $G_2$ has been group flyped (across $T_1$).

\begin{definition}
  A prime knot configuration is said to be a full-group configuration if
  every group is a full group. Otherwise, the knot configuration is said
  to be a split-group configuration.
\end{definition}  

 Our next objective is to establish that if $G$ and $G'$ are groups
 for which $G_F\ne G'_F$, then the orbit of $G$ and the orbit of $G'$
 are in a sense orthogonal to each other. By this we mean that there is
 a min-tangle $T'$ of the orbit of $G'$ such that $C(K)-T'$ is contained
 within a min-tangle $T$ of the orbit of $G$. In order to show this, we
 shall first of all need to identify the possible flype moves for an
 individual crossing. In order to do this, we shall introduce a notion
 of flype scenario for an individual crossing. 
 
\begin{definition}
 Let $C(K)$ be a configuration of a prime knot $K$, and let $c$ be a
 crossing of $C(K)$. A crossing flype scenario for $c$ on edges $e$ and
 $f$ is a sequence $(c,e,f,T_1,T_2)$, where $T_1$ and $T_2$ are tangles,
 $e$ and $f$ are adjacent arcs at $c$, incident to $T_1$, and $C(K)$ has
 the form

  \sarcenter
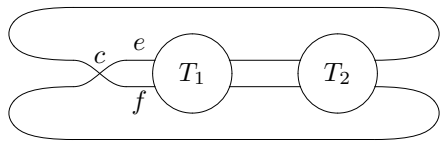
{\xy /r20pt/:,
   (0,1.75)*{\hbox{}},
   (5.25,-.25)="a11",
   ;"a11";{\ellipse(.75){}},
   (5.25,-.25)*{\hbox{$T_1$}},
   (2,.25)="a",
   (8,-.25)="a12",
   ;"a12";{\ellipse(.75){}},
   (8,-.25)*{\hbox{$T_2$}},
   (4,0)="a2",(4.55,0)="b2","a2";"b2"**\dir{-},
   (4,-.5)="a3",(4.55,-.5)="b3","a3";"b3"**\dir{-},
   (5.95,0)="a4",(7.3,0)="b4","a4";"b4"**\dir{-},
   (5.95,-.5)="a5",(7.3,-.5)="b5","a5";"b5"**\dir{-},
   (3,-1.5)="a6",(8.7,-1.5)="b6","a6";"b6"**\dir{-},
   (3,1)="a7",(8.7,1)="b7","a7";"b7"**\dir{-},
   (3,0)="a8";"a3"**\crv{"a8"+(.4,0) & "a3"-(.4,0)},
   (3,-.5)="a9";"a2"**\crv{"a9"+(.4,0) & "a2"-(.4,0)},
   (3.5,-.25)*!<0pt,-6pt>{c},
   (4.25,0)*!<0pt,-6pt>{e},
   (4.25,-.5)*!<0pt,6pt>{f},
   (3,1)="c1",
   {\hloop-}="c1",
   (3,-.5)="c2",
   {\hloop-}="c2",
   (8.7,1)="c3",
   {\hloop}="c3",
   (8.7,-.5)="c4",
   {\hloop}="c4",
  \endxy}
\end{definition}
   
 Note that unlike the situation for flype scenario, it is not implied
 that the endpoints of $e$ and $f$ that belong to $T_1$ be distinct, nor
 that the endpoints of the other two arcs incident to $c$ that belong to
 $T_2$ be distinct. For example, if $K$ is a 3-crossing torus knot, then
 we might take $c$ to be the first crossing in the group of three
 crossings, $T_1$ to be just the second crossing with its incident arcs
 and $T_2$ to be the third crossing with its incident arcs, with $e$ and
 $f$ being the arcs from $c$ to $T_1$.

 It turns out that if we label the edges incident to a crossing
 $c$ in the clockwise direction as $e$, $f$, $e_1$ and $f_1$, and
 there is a crossing flype scenario for $c$ on edges $e$ and $f$, then
 there is no crossing flype scenario for $c$ on edges $f$ and $e_1$.

\begin{proposition}\label{crossing flype scenario}
 Let $C(K)$ be a configuration of a prime knot $K$ and let $c$ be
 a crossing of $C(K)$. Let $e$, $f$ and $e_1$ be edges incident to
 $c$ with $e$ and $f$ adjacent, and $f$ and $e_1$ adjacent. If there
 is a crossing flype scenario for $c$ on edges $e$ and $f$, then
 there is no crossing flype scenario on edges $f$ and $e_1$.
\end{proposition}

\begin{proof}
 Suppose to the contrary that $(c,e,f,T_1,T_2)$ and $(c,f,e_1,T'_1,T'_2)$
 are crossing flype scenarios. Let the endpoints of
 $e$ and $f$ in $T_1$ be denoted by $a$ and $a'$, respectively, and
 let the endpoint of $e_1$ that is in $T'_1$ be denoted by $b'$. Let
 $f_1$ denote the fourth edge incident to $c$, and denote its other
 endpoint by $b$. Then we have $a$ and $a'$ in $T_1$, $b$ and $b'$
 in $T_2$, $a'$ and $b'$ in $T'_1$ and $a$ and $b$ in $T'_2$. Let $h_1$
 and $h_2$ denote the two edges incident to both $T_1$ and $T_2$, and
 let $h'_1$ and $h'_2$ denote the two edges incident to both $T'_1$ and
 $T'_2$. Since $a$ and $a'$ are in $T_1$, which is connected, there is
 a path from $a$ to $a'$ using only edges of $T_1$. But the graph that
 is obtained by deleting $c$ and the two edges $h'_1$ and $h'_2$ from
 $C(K)$ has two connected components, namely $T'_1$ and $T'_2$, and
 $a$ is in $T'_2$, while $a'$ is in $T'_1$. Thus every path in $C(K)$ from
 $a$ to $a'$ that does not go through $c$ must use either edge $h'_1$
 or else $h'_2$.

 \sarcenter
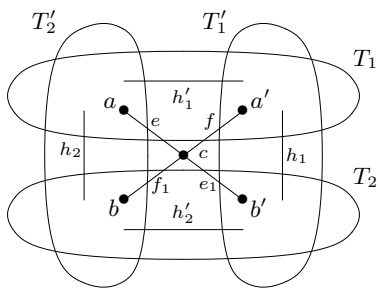
{\xy /r30pt/:(0,.75)::,
   (1.25,1.25)="a1",(2.75,2.75)="b1","a1";"b1"**\dir{-},
   (1.25,2.75)="a2",(2.75,1.25)="b2","a2";"b2"**\dir{-},
   (.75,2.75)="a3",(.75,1.25)="b3","a3";"b3"**\dir{-},
   (1.25,3.25)="a4",(2.75,3.25)="b4","a4";"b4"**\dir{-},
   (3.25,2.75)="a5",(3.25,1.25)="b5","a5";"b5"**\dir{-},
   (1.25,.75)="a6",(2.75,.75)="b6","a6";"b6"**\dir{-},
   (.25,2)="c1",(1.55,2)="c2","c1";"c2"**\crv{(.25,3.75) & (.75,4.25) & (1.25,4.25) & (1.55,3.75)},
   "c2";"c1"**\crv{(1.55,.25) & (1.25,-.25) & (.75,-.25) & (.25,.25)},
   (2.45,2)="c3",(3.75,2)="c4","c3";"c4"**\crv{(2.45,3.75) & (2.75,4.25) & (3.25,4.25) & (3.75,3.75)},
   "c4";"c3"**\crv{(3.75,.25) & (3.25,-.25) & (2.75,-.25) & (2.45,.25)},
   (2,2.25)="c5",(2,3.75)="c6","c5";"c6"**\crv{(.25,2.25) & (-.25,2.75) & (-.25,3.25) & (.25,3.75)},
   "c6";"c5"**\crv{(3.75,3.75) & (4.25,3.25) & (4.25,2.75) & (3.75,2.25)},
   (2,1.75)="c7",(2,.25)="c8","c7";"c8"**\crv{(.25,1.75) & (-.25,1.25) & (-.25,.75) & (.25,.25)},
   "c8";"c7"**\crv{(3.75,.25) & (4.25,.75) & (4.25,1.25) & (3.75,1.75)},
   (1.25,1.25)*{\hbox{\small$\bullet$}},(2.75,2.75)*{\hbox{\small$\bullet$}},
   (1.25,2.75)*{\hbox{\small$\bullet$}},(2.75,1.25)*{\hbox{\small$\bullet$}},
   (2,2)*{\hbox{\small$\bullet$}},
   (2.25,2.025)*{\hbox{$\ssize c$}},
   (4.3,3.625)*{\hbox{\small$T_1$}},
   (2.4,4.25)*{\hbox{\small$T'_1$}},
   (4.3,1.625)*{\hbox{\small$T_2$}},
   (.25,4.25)*{\hbox{\small$T'_2$}},
   (1.075,2.875)*{\hbox{\small$a$}},
   (2.975,2.95)*{\hbox{\small$a'$}},
   (1.125,1.125)*{\hbox{\small$b$}},
   (2.95,1.125)*{\hbox{\small$b'$}},
   (3.5,2)*!<2pt,0pt>{\hbox{$\ssize h_1$}},
   (2,2.85)*!<0pt,-2.5pt>{\hbox{$\ssize h'_1$}},
   (.6,2.1)*!<.5pt,0pt>{\hbox{$\ssize h_2$}},
   (2,1.1)*!<0pt,1.5pt>{\hbox{$\ssize h'_2$}},
   (1.65,2.7)*!<0pt,3pt>{\hbox{$\ssize e$}},
   (2.3,2.7)*!<-1pt,3pt>{\hbox{$\ssize f$}},
   (2.32,1.355)*!<0pt,-3pt>{\hbox{$\ssize e_1$}},
   (1.67,1.355)*!<-1.5pt,-3pt>{\hbox{$\ssize f_1$}},
  \endxy}

 \sarnoindent Since there is a path from $a$ to $a'$ that uses only
 edges in $T_1$, either $h'_1$ or $h'_2$ is an edge in $T_1$. A similar
 argument applied to $b$ and $b'$ shows that either $h'_1$ or $h'_2$ must
 be an edge in $T_2$. Thus exactly one of $h'_1$ and $h'_2$ belongs to
 $T_1$, the other belongs to $T_2$. We may suppose without loss of
 generality that $h'_1$ is in $T_1$ and $h'_2$ is in $T_2$. Similarly,
 we may suppose without loss of generality that $h_1$ is in $T'_1$ and
 $h_2$ is in $T'_2$. Now $a$ belongs to both $T_1$ and $T'_2$, so
 $T_1\cap T'_2$ is an $m$-tangle for some $m\ge4$. Now $h'_1$ is an edge
 of $T_1$ and incident to $T'_2$, so $h'_1$ is incident to $T_1\cap T'_2$.
 Similarly, $h_2$ is incident to $T_1\cap T'_2$. As well, $e$ is
 incident to $T_1\cap T'_2$. There must be at least one more edge $g$
 incident to $T_1\cap T'_2$. The edges incident to $T_1$ are $e$, $f$,
 $h_1$ and $h_2$, and $g$ is not equal to $e$ or $h_2$. Since $f$ is
 incident to $T'_1$ and not to $T'_2$, $g$ is not equal to $f$. Finally,
 $h_1$ is an edge in $T'_1$ and therefore not incident to $T'_2$, whence
 $g$ is not equal to $h_1$. Thus $g$ must be an edge of $T_1$. But now
 a similar argument shows that $g$ must also be an edge of $T'_2$, whence
 $g$ is an edge of $T_1\cap T'_2$. But this contradicts the fact that
 $g$ was incident to $T_1\cap T'_2$, and so it is not possible to have
 crossing flype scenarios $(c,e,f,T_1,T_2)$ and $(c,f,e_1,T'_1,T'_2)$.
\end{proof}

This leads to the next result, which may be thought of as establishing
the independence of orbits, or the non-interference of orbits.

\begin{theorem}\label{non-interference of orbits}
 Let $C(K)$ be a configuration of a prime alternating knot $K$, and let
 $G$ be a group of $K$ and $T_1$ be a min-tangle of the orbit of $G$.
 Let $c$ be a crossing in $T_1$ and suppose that $T$ is a tangle over
 which $c$ may flype. Then either $T\subseteq T_1$ or else $C(K)-T_1
 \subseteq T$.
\end{theorem}

\begin{proof}
  To begin with, we observe that we may assume without loss of generality
  that $T\cap T_1\ne \emptyset$. To see why this is so,
  we note that $c$ may flype over the tangle $T'=
  C(K)-T-\{c\}$, and if $T\cap T_1=\emptyset$, then $T'\cap T_1\ne
  \emptyset$. Now, $T\subseteq T_1$ if and only if $C(K)-T_1\subseteq
  C(K)-T$, and since $c\in T_1$, we know that $c\notin C(K)-T_1$,
  whence $C(K)-T_1\subseteq C(K)-T$ if and only if $C(K)-T_1\subseteq
  C(K)-T-\{c\}$. Thus $T\subseteq T_1$ if and only if $C(K)-T_1\subseteq
  T'$. Finally, since $T=C(K)-T'-\{c\}$, we obtain by symmetry that
  $C(K)-T_1\subseteq T$ if and only if $T'\subseteq T_1$. Consequently,
  $T\subseteq T_1$ or $C(K)-T_1\subseteq T$ if and only if
  $C(K)-T_1\subseteq T'$ or $T'\subseteq T_1$. 

  We have established now that we may assume that $T\cap T_1\ne \emptyset$.
  If $T\subseteq T_1$, then there is nothing to prove, so we consider
  the case when $T$ contains at least one crossing that is not in $T_1$.
  We must show that $C(K)-T_1\subseteq T$. Suppose to the contrary that
  $C(K)-T_1\nsubseteq T$, and let $T_0=C(K)-T_1$. Label
  one of the position pairs of arcs for $T_1$ as $e$ and $f$, with the other
  pair being labelled with $e'$ and $f'$.
  
  \sarcenter
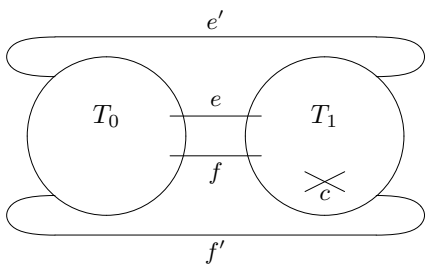
{\xy /r15pt/:,
   (2,2)="a11",
   ;"a11";{\ellipse(2){}},
   (2,2.5)*{\hbox{$T_0$}},
   (7.5,2.5)*{\hbox{$T_1$}},
   (7.5,2)="a12",
   ;"a12";{\ellipse(2){}},
   (3.6,2.5)="a2",(5.9,2.5)="b2","a2";"b2"**\dir{-},
   (3.6,1.5)="a3",(5.9,1.5)="b3","a3";"b3"**\dir{-},
   (.7,-.5)="a6",(8.8,-.5)="b6","a6";"b6"**\dir{-},
   (.7,4.5)="a7",(8.8,4.5)="b7","a7";"b7"**\dir{-},
   (.7,4.5)="c1",
   {\hloop-}="c1",
   (.7,.5)="c2",
   {\hloop-}="c2",
   (13,2)="a13",
   (7,.6)="a8",(8,1.1)="b8","a8";"b8"**\dir{-},
   (8,.6)="a9",(7,1.1)="b9","a9";"b9"**\dir{-},
   (8.8,4.5)="c3",
   {\hloop}="c3",
   (8.8,.5)="c4",
   {\hloop}="c4",
   (4.75,4.8)*!<0pt,-2pt>{\hbox{\small$e'$}},
   (4.75,-.8)*!<0pt,2pt>{\hbox{\small$f'$}},
   (4.75,2.75)*!<0pt,-2pt>{\hbox{\small$e$}},
   (4.75,1.2)*!<0pt,2pt>{\hbox{\small$f$}},
   (7.5,.5)*{\hbox{\small$c$}},
  \endxy}
   
\sarnoindent 
Now, $T_0\nsubseteq T$ implies that $T_0-T\ne\emptyset$, whence
$T_0-T$ is an $m_1$-tangle and $T_0\cap T$ is an $m_2$-tangle for
some even integers $m_1,m_2\ge 4$. If there was at most one edge
joining a vertex of $T_0-T$ to a vertex of $T_0\cap T$, then there
would be least 3 edges each with one endpoint in $T_0-T$ and
the other endpoint not in $T_0$. Since $T_0$ has exactly 4 incident
edges, this would mean that there is exactly one edge with one
endpoint in $T_0\cap T$ and the other endpoint not in $T_0$, in which
case $m_2=2$, which is not the case. Thus there are at least two
edges each with one endpoint in $T_0-T$ and the other endpoint in $T_0\cap
T$. Such arcs belong to $T_0$, so there are at least two arcs of $T_0$
that are incident to $T$. Similarly, since $T_1\cap T\ne\emptyset$ and
$T_1-T\ne\emptyset$, there are at least two arcs of $T_1$ that are
incident to $T$. But $T$ is a tangle, so there are exactly four arcs
incident to $T$, whence there must be exactly two arcs between $T_0-T$
and $T_0\cap T$, and exactly two arcs between $T_1-T$ and $T_1\cap T$.
Thus of the four arcs incident to $T$, two are arcs of $T_0$ and two are
arcs of $T_1$, whence none of the four arcs $e,f,e',f'$ are incident to
$T$. Since $m_1\ge 4$, and there are exactly two arcs from $T_0-T$ to
$T_0\cap T$, at least two of $e,f,e'f'$ are incident to $T_0\cap T$,
hence have both endpoints in $T$. There are two cases to consider:
both arcs of one of the position pairs for $T_1$ belong to $T$, or
exactly one from each of the two position pairs for $T_1$ belongs to
$T$.

\noindent Case 1: both arcs of one of the position pairs for $T_1$ belong
to $T$. Without loss of generality, we may assume that $e$ and
$f$ belong to $T$. But then $T_1\cap T$ is a tangle with $e$ and $f$
incident edges, contradicting the minimality of $T_1$. Thus this case
can't occur.

\noindent Case 2: exactly one from each of the two position pairs for
$T_1$ belongs to $T$. Suppose that the position pairs have been labelled
so that $e$ and $e'$ are the two arcs that belong to $T$. Since two of
the four edges incident to $T$ are edges of $T_1$, and the other two
are edge of $T_0$, we see that $f$ and $f'$ are not in $T$ nor are they
incident to $T$. Thus $e$ and $f$ have different endpoints in $T_0$,
and $e'$ and $f'$ have different endpoints in $T_0$. Thus $e$
and $f$ are not edges incident to the group $G$, nor are $e'$ and $f'$
incident to $G$, whence $T_0$ has a decomposition (obtained from the
orbit of $G$ as a tangle $T_0'$, the group $G$ and a tangle $T_0''$,
so $C(K)$ has the structure

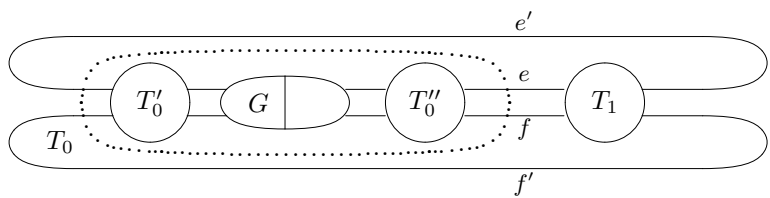
\sarcenter{
\xy /r20pt/:,
   (4.55,-.25)="a11",
   ;"a11";{\ellipse(.75){}},
   (4.55,-.25)*{\hbox{$T'_0$}},
   (6.6,-.25)*{\hbox{$G$}},
   (7.1,.25)="q",
   {\hloop\hloop-}="q",
   (7.1,.25)="q1",(7.1,-.75)="r1","q1";"r1"**\dir{-},
   (3.1,0)="a2",(3.85,0)="b2","a2";"b2"**\dir{-},
   (3.1,-.5)="a3",(3.85,-.5)="b3","a3";"b3"**\dir{-},
   (5.25,0)="a4",(6,0)="b4","a4";"b4"**\dir{-},
   (5.25,-.5)="a5",(6,-.5)="b5","a5";"b5"**\dir{-},
   (3.1,-1.5)="a6",(15,-1.5)="b6","a6";"b6"**\dir{-},
   (3.1,1)="a7",(15,1)="b7","a7";"b7"**\dir{-},
   (3.1,1)="c1",
   {\hloop-}="c1",
   (3.1,-.5)="c2",
   {\hloop-}="c2",
   (9.75,-.25)="a13",
   ;"a13";{\ellipse(.75){}},
   (9.75,-.25)*{\hbox{$T''_0$}},
   (8.25,0)="a8",(9,0)="b8","a8";"b8"**\dir{-},
   (8.25,-.5)="a9",(9,-.5)="b9","a9";"b9"**\dir{-},
   (15,1)="c3",
   {\hloop}="c3",
   (15,-.5)="c4",
   {\hloop}="c4",
  (13.15,-.25)="a14",
   ;"a14";{\ellipse(.75){}},
   (13.15,-.25)*{\hbox{$T_1$}}, 
   (13.85,0)="a10",(15,0)="b10","a10";"b10"**\dir{-},
   (13.85,-.5)="a11",(15,-.5)="b11","a11";"b11"**\dir{-},
   (10.5,0)="l1",(12.375,0)="p1","l1";"p1"**\dir{-}?(.6)
   *!<0pt,-5pt>{\hbox{\small$e$}}*!<0pt,-26pt>{\hbox{\small$e'$}},
   (10.5,-.5)="l2",(12.375,-.5)="p2","l2";"p2"**\dir{-}
   ?(.6)*!<0pt,5pt>{\hbox{\small$f$}}*!<0pt,26pt>{\hbox{\small$f'$}},
  (3.25,-.25)="t1",
(11.35,-.25)="t2",
"t1";"t2"**\crv{~*=<\jot>{.}(3.35,.5) & (3.6,.75) & (11,.75) & (11.25,.5)},
"t2";"t1"**\crv{~*=<\jot>{.}(11.25,-1) & (11,-1.25) & (3.6,-1.25) & (3.35,-1)},
 (2.85,-1)*{\hbox{$T_0$}},  
\endxy}

Now $T\cap T_0''$, $T_0''-T$ and $T-T_0''$ are all nonempty, so by
Lemma \ref{tangle intersection}, there are two arcs of $T_0''$ incident to
$T$. Similarly, there are two arcs of $T_0'$ incident to $T$. As well,
there are two arcs of $T_1$ incident to $T$, whence we have accounted
for six arcs incident to $T$. But this is not possible, so the
assumption that $C(K)-T_1\nsubseteq T$ must be false. We conclude therefore
that $C(K)-T_1\subseteq T$.
\end{proof}   

\begin{corollary}\label{group orbit noninterference}
 Let $C(K)$ be a configuration of a prime knot $K$, and let
 $G$ be a group of $C(K)$. If $T_1$ is any min-tangle from the orbit of
 $G$ and $H$ is any group of $C(K)$ for which $H\cap T_1\ne \emptyset$,
 then $H\subseteq T_1$.
\end{corollary}

\begin{proof}
 Let $H$ be a group such that $H\cap T_1\ne\emptyset$, so there
 is a crossing $c$ of $H$ that is in $T_1$. Suppose that $H$ is
 not contained entirely in $T_1$, so there is a crossing $d$ of $H$
 that is not in $T_1$. Let $T$ denote the connected component of $H-\{c\}$
 that contains $d$. Then $T$ is a tangle over which $c$ may flype.
 By Theorem \ref{non-interference of orbits}, we either have $T\subseteq T_1$
 or else $C(K)-T_1\subseteq T$. Since $d\in T-T_1$, it must be that
 $C(K)-T_1\subseteq T$. But then $G\subseteq T$, whence $d$ is a crossing
 in $G$. By the definition of group, this means that $c$ is in $G$,
 which contradicts the fact that $T_1$ is a min-tangle in the orbit
 of $G$. Since this contradiction follows from the assumption that
 $H$ is not contained entirely in $T_1$, we conclude that $H\subseteq T_1$.
\end{proof}

\begin{definition}
  A {\it flype scenario} $\flype G,T_1,T_2,$ in a knot configuration
  $C(K)$ consists of a group $G$ and two non-trivial tangles $T_1$ and
  $T_2$ such that $C(K)$ is of the form

  \sarcenter
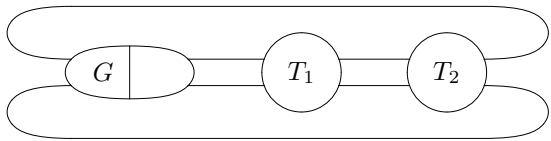
{\xy /r20pt/:,
   (5.25,-.25)="a11";"a11";{\ellipse(.75){}},
   (5.25,-.25)*{\hbox{$T_1$}},
   (2,.25)="a",
   {\hloop\hloop-}="a",
   (8,-.25)="a12";"a12";{\ellipse(.75){}},
   (8,-.25)*{\hbox{$T_2$}},
   (1.5,-.25)*{\hbox{$G$}},
   (2,.25)="a1",(2,-.75)="b1","a1";"b1"**\dir{-},
   (3.1,0)="a2",(4.55,0)="b2","a2";"b2"**\dir{-},
   (3.1,-.5)="a3",(4.55,-.5)="b3","a3";"b3"**\dir{-},
   (5.95,0)="a4",(7.3,0)="b4","a4";"b4"**\dir{-},
   (5.95,-.5)="a5",(7.3,-.5)="b5","a5";"b5"**\dir{-},
   (.9,-1.5)="a6",(8.7,-1.5)="b6","a6";"b6"**\dir{-},
   (.9,1)="a7",(8.7,1)="b7","a7";"b7"**\dir{-},
   (.9,1)="c1",
   {\hloop-}="c1",
   (.9,-.5)="c2",
   {\hloop-}="c2",
   (8.7,1)="c3",
   {\hloop}="c3",
   (8.7,-.5)="c4",
   {\hloop}="c4",
  \endxy}
\end{definition}

 Observe that $T_1$ itself cannot be decomposed in the form

 \sarcenter
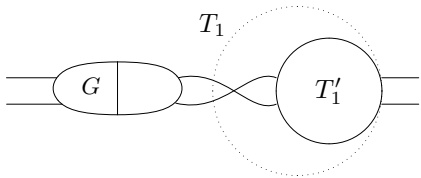
{\xy/r20pt/:,
 (7.5,-.25)="a11",
 ;"a11";{\ellipse(1){}},
 (7.5,-.25)*{\hbox{$T'_1$}},
 (3.5,.3)="a",
 {\hloop\hloop-}="a",
 (3,-.15)*{\hbox{\small$G$}},
 (3.5,.3)="a1",(3.5,-.7)="b1","a1";"b1"**\dir{-},
 (4.65,0)="q1";(6.5,-.5)="q2"**\crv{(5.5,.2) & (6,-.7)},
 (4.6,-.48)="q3";(6.5,0)="q4"**\crv{(5.5,-.7) & (6,.2)},
 (6.9,-.25)="a12",
 ;"a12";{\ellipse(1.6){.}},
 (5.3,1)*{\hbox{$T_1$}},
 (1.4,0)="c1",(2.35,0)="d1","c1";"d1"**\dir{-},
 (1.4,-.5)="c2",(2.45,-.5)="d2","c2";"d2"**\dir{-},
 (8.5,0)="c3",(9.2,0)="d3","c3";"d3"**\dir{-},
 (8.5,-.5)="c4",(9.2,-.5)="d4","c4";"d4"**\dir{-},
\endxy}

\noindent In other words, the two arcs that join the group $G$ to $T_1$
 have distinct endpoints in $T_1$. Similarly, the two arcs that join the
 group $G$ to $T_2$ must have distinct endpoints in $T_2$.

 It is evident that given a flype scenario $\flype G,T_1,T_2,$, each crossing
 of $G$ can be flyped to the two arcs joining $T_1$ and $T_2$. Conversely,
 suppose that $c$ is a crossing for which there is a flype move: that is,
 there is a tangle $T$ across which $c$ can be flyped. Let $G$ denote the
 group that contains $c$. There are several possibilities to be investigated. 

\noindent Case 1: $G\cap T=\emptyset$. 

\noindent Case 1a) $G\cup T=C(K)$. In this case the flype move does not
 change the knot configuration.

\noindent Case 1b) $G\cup T\ne C(K)$. Then $T_1=C(K)-(G\cup T)$ is a
 non-trivial tangle and $G$, $T$, and $T_1$ constitute a flype scenario.

\noindent Case 2: $G\cap T\ne \emptyset$ (by which we mean there is at
 least one crossing that is common to $G$ and to $T$). There are two
 subcases to be addressed.

\noindent Case 2a) $T\subseteq G$. In this case, the flype move does not
 change the knot configuration. 

\noindent Case 2b) $T\not\subseteq G$. Let $T'=T-G$, so that $T'$ is a
 tangle with $G\cap T'=\emptyset$. Furthermore, at least one end of $G$ is
 joined to $T'$. Let $c'$ denote the crossing at an end of $G$ which is
 joined to $T'$. There is a finite sequence of flype moves, each of which
 will leave the knot configuration unchanged, which will flype $c$ to the
 arcs joining $c'$ to $T'$. One final flype move will flype $c$ across $T'$,
 resulting in the knot configuration that is created by flyping $c$ across
 the original tangle $T$. We observe that this is the same configuration
 that would result from flyping $c'$ across $T'$. Accordingly, we replace
 $c$ by $c'$ and $T$ by $T'$ to obtain a situation where $c$ is to be flyped
 across $T$ and $G\cap T=\emptyset$. This situation has been handled in
 Case 1. 

 The preceding discussion makes it clear that the only flype moves of
 importance are those identified by a flype scenario.   

\begin{theorem}\label{full group}
 Given a knot configuration $C(K)$ of an alternating prime knot $K$,
 there is a finite sequence of group flype moves which will transform
 $C(K)$ into a full-group configuration.
\end{theorem}

\begin{proof}  The proof is by induction on the number of groups in the
 configuration, and the base case would be a knot configuration with a
 single group, which is already a full-group configuration. Suppose now
 that $n\ge1$ is an integer such that for any prime alternating knot
 $K$, any configuration of $K$ with $n$ groups can be transformed into
 a full-group configuration of $K$ by a finite sequence of group flype
 moves. Let $C(K)$ be a knot configuration with $n+1$ groups. If $C(K)$
 is a full-group configuration, then there is nothing to do. Suppose
 then that $C(K)$ is a split-group configuration. Let $G_1$ and $G_2$ be
 part of a split group. Then there are non-trivial tangles $T_1$ and
 $T_2$ such that $C(K)$ has the form

   \sarcenter
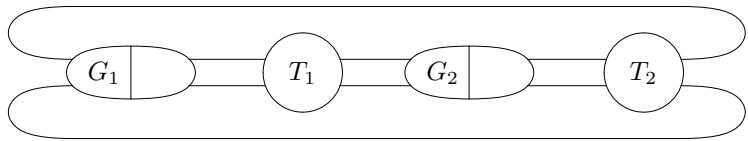
{\xy /r20pt/:,
   (5.25,-.25)="a11",
   ;"a11";{\ellipse(.75){}},
   (5.25,-.25)*{\hbox{$T_1$}},
   (2,.25)="a",
   {\hloop\hloop-}="a",
   (7.9,-.25)*{\hbox{$G_2$}},
   (8.4,.25)="q",
   {\hloop\hloop-}="q",
   (8.4,.25)="q1",(8.4,-.75)="r1","q1";"r1"**\dir{-},
   (1.5,-.25)*{\hbox{$G_1$}},
   (2,.25)="a1",(2,-.75)="b1","a1";"b1"**\dir{-},
   (3.1,0)="a2",(4.55,0)="b2","a2";"b2"**\dir{-},
   (3.1,-.5)="a3",(4.55,-.5)="b3","a3";"b3"**\dir{-},
   (5.95,0)="a4",(7.3,0)="b4","a4";"b4"**\dir{-},
   (5.95,-.5)="a5",(7.3,-.5)="b5","a5";"b5"**\dir{-},
   (.9,-1.5)="a6",(12.4,-1.5)="b6","a6";"b6"**\dir{-},
   (.9,1)="a7",(12.4,1)="b7","a7";"b7"**\dir{-},
   (.9,1)="c1",
   {\hloop-}="c1",
   (.9,-.5)="c2",
   {\hloop-}="c2",
   (11.7,-.25)="a13",
   ;"a13";{\ellipse(.75){}},
   (11.7,-.25)*{\hbox{$T_2$}},
   (9.5,0)="a8",(11,0)="b8","a8";"b8"**\dir{-},
   (9.5,-.5)="a9",(11,-.5)="b9","a9";"b9"**\dir{-},
   (12.4,1)="c3",
   {\hloop}="c3",
   (12.4,-.5)="c4",
   {\hloop}="c4",
  \endxy}
  
\sarnoindent and $G_2$ can be group flyped over $T_1$ (or $T_2$) to merge with
  $G_1$. Let $G'$ denote the group that results when $G_2$ is group
  flyped over $T_1$ to merge with $G_1$. Then the orbit of $G'$ in the
  resulting knot configuration $C_1(K)$ contains one fewer groups than the orbit
  of $G_1$ had in the original configuration $C(K)$. If $H$ is any 
  group of $C(K)$ that does not belong to the orbit of $G_1$, then
  by Corollary \ref{group orbit noninterference}, $H$ is contained within some
  min-tangle of the orbit of $G_1$, whence $H$ is contained entirely within
  either $T_1$ or else $T_2$, in which case $H$ is still a group in $C_1(K)$.
  Thus the configuration $C_1(K)$ has $n$ groups. By the induction
  hypothesis, this configuration can be transformed into a
  full-group configuration by a finite sequence of group flype moves.
  Consequently, we obtain a finite sequence of group flype moves that
  transforms $C(K)$ into a full-group configuration. Now the result follows
  by induction.
\end{proof}

 Note that in a prime knot, every arc incident to a tangle must terminate at
 a vertex of the tangle.

We conclude this section with an observation about the orbit of
a negative group. To begin with, note that once a knot has been
given an orientation, it follows that for any $m$-tangle, exactly
half of the incident arcs are directed into the $m$-tangle, with
the other half of course directed out of the $m$-tangle. In particular,
if $G$ is a group in an oriented prime alternating knot configuration
$C(K)$, then for any min-tangle $T$ in the orbit of $G$, two of the
incident arcs are directed into the tangle and two are directed out of
the tangle. If $G$ is a positive group, then both arcs of each position
pair for $T$ are oriented the same, where as we leave the group to
traverse the knot, the first time that we follow an arc incident to
$T$ establishes that position pair as being directed into $T$, while
the other position pair will be directed outward from $T$. On the
other hand, if $G$ is a negative group, then the two arcs in both
position pairs have opposite orientation.

\begin{proposition}\label{core of neg group}
 Let $C(K)$ be a configuration of a prime alternating knot $K$, and
 let $G$ be a negative group of $C(K)$. There is a unique min-tangle
 $T$ from the orbit of $G$ such that a knot traversal of $K$ begun
 at $G$ travels both arcs of one position pair of $T$ before it
 travels an arc from the other position pair of $T$.
\end{proposition}

\begin{proof}
Let $(\stackrel{*}{G_0}, T_0, \stackrel{*}{G_1}, T_1,\ldots,
\stackrel{*}{G_k},T_k)$ be the orbit of $G$, so that
$T_0, T_1,\ldots,T_k$ are the min-tangles of the orbit of $G$.
Observe that we may perform any necessary group flypes to put
$G_F$ into full group without changing the min-tangles or the
directions assigned to their incident arcs by assigning an orientation
to $C(K)$. We may therefore assume that $G$ is a full group, and that
$\stackrel{*}{G_i}=\{\,e_{i},f_{i}\,\}$ for each $i=1,\ldots,k$.
As well, let the position pair for $T_0$ which connect $G$ to $T_0$
at the beginning of the orbit of $G$ be denoted by $\{\,e_0,f_0\,\}$.
Give $G$ an orientation and suppose that the arcs have been labelled
so that when a knot traversal is performed in the direction determined by
the orientation, starting at $G$, we move initially
towards $T_0$ on $e_0$, and that for each $i=0,1,\ldots,k$, the arc that
is used to move from $T_i$ to $T_{i+1}$ (or to $G$ if $i=k$) is $e_{i+1}$.
The orbit of $G$ is thus of the form

\immediate\write16{** create diagram of orbit of neg group with position
arcs labelled--goes in orbit.tex at Proposition \ref{core of neg group}}

Suppose to begin with that as the knot is traversed, the traversal takes
us into tangle $T_i$ on $e_i$ and leaves via arc $i_{i+1}$ for every
$i=0,1,\ldots,k$. Then we first return to $G$ via arc $e_{k+1}$, and
since the two arcs $e_{k+1}$ and $f_{k+1}$ meet at a crossing of $G$,
the traversal will not take us out of $G$ via $f_{k+1}$. It follows
that the traversal will take us out of $G$ via $e_{0}$, which means
that our knot is a link. Since this is not the case, there must be
some index $i$ such that the knot traversal takes us into $T_i$ via
arc $e_i$ but then leaves $T_i$ next via arc $f_i$, whence $T_i$ is
a min-tangle with the required properties.

The uniqueness follows similarly from the observation that if there
were two such min-tangles, then again we would have to conclude that
$K$ is a link, which it isn't.
\end{proof}

\begin{definition}\label{core def}
 Let $C(K)$ be a configuration of a prime alternating knot $K$, and
 let $G$ be a negative group of $C(K)$. The unique min-tangle
 whose existence was established in Proposition \ref{core of neg group} is called
 the {\em core} of the orbit of $G$.
\end{definition}

\par
\section{The Construction Process}

In this section, we introduce the four operators that will be
utilized in the construction of the prime alternating knots. Two
of the four, $D$ and $ROTS$, take a prime alternating knot configuration with $n$ 
crossings and return a prime alternating knot configuration with 
$n+1$ crossings, while the other two operators do not change the
number of crossings.

\subsection{The $D$ operator}

\begin{figure}[ht]
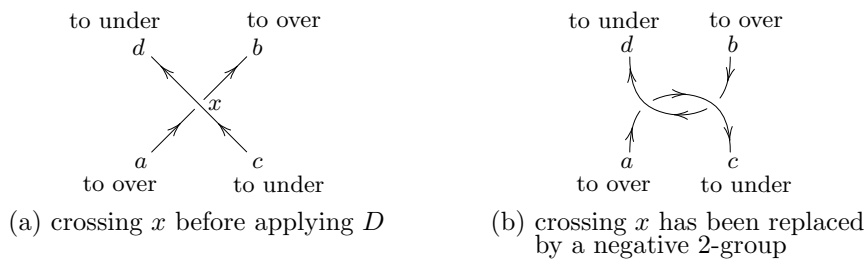

\centering
   \begin{tabular}{c@{\hskip 40pt}c}
   \noalign{\vskip6pt}
$\vcenter{
\xy /r36pt/:,
 (0,0)="a";(1,1)="b"**\dir{}?(.5)="x"?(.05)="u",
 "a";"x"-"u"**\dir{-},
 "x"+"u";"b"**\dir{-},
 "x"*!<-6pt,0pt>{\hbox{\small $x$}},
 "a";"b"**\dir{}*!<-4pt,-4pt>{\hbox{\small $b$}}?(.9)*{\dir{>}}
 ?(.3)*{\dir{>}},
 (1,0)="c";(0,1)="d"**\dir{-}*!<5pt,-4pt>{\hbox{\small $d$}}?(.9)*{\dir{>}}
 ?(.3)*{\dir{>}},
 "b"*!<-14pt,-14pt>{\hbox{\small to over}},
 "d"*!<14pt,-14pt>{\hbox{\small to under}},
 "a"*!<4pt,4pt>{\hbox{\small $a$}},
 "a"*!<12pt,12pt>{\hbox{\small to over}},
 "c"*!<-4pt,4pt>{\hbox{\small $c$}},
 "c"*!<-12pt,12pt>{\hbox{\small to under}},  
\endxy}$
&
$\vcenter{\xy /r38pt/:<0pt,18pt>::,
\vloop <(.05)<<(=.75)\khole <(.48)< <(.92)<,
(0,3),\vloop- <(.05)<  <(.48)< <(=.75)\khole <(.92)<,
(0,1)="a"*!<1pt,4pt>{\hbox{\small $a$}},
(1,1)="c"*!<-1pt,4pt>{\hbox{\small $c$}},
(0,3)="d"*!<1pt,-5pt>{\hbox{\small $d$}},
(1,3)="b"*!<-1pt,-5pt>{\hbox{\small $b$}},
 "b"*!<-6pt,-14pt>{\hbox{\small to over}},
 "d"*!<6pt,-14pt>{\hbox{\small to under}},
 "a"*!<6pt,12pt>{\hbox{\small to over}},
 "c"*!<-6pt,12pt>{\hbox{\small to under}},  
\endxy}$\\
   \noalign{\vskip6pt}
   (a) crossing $x$ before applying $D$  & (b) \vtop{\hsize 1.75in
     \parindent=0pt\leftskip=0pt crossing $x$ has been replaced by
       a negative 2-group}\\
   \noalign{\vskip6pt}
\end{tabular}
\caption{The $D$ operator}\label{d operator}
\end{figure}

The $D$ operator can be applied to any crossing of a knot configuration,
where the crossing is treated as if it were a negative loner, and the end
result is to replace the crossing by a negative group of 2. More precisely,
let $C(K)$ be a configuration of a prime alternating knot $K$, and let
$x$ be a crossing of $C(K)$. Label the arcs incident to $x$ as
shown in Figure \ref{d operator} (a), so that
under a traversal of the knot, arc $a$ is used to visit $x$ for the
first time, and arc $b$ is used to leave $x$ after that visit, while
arc $c$ is used for the second visit, with arc $d$ used to leave $x$
after the second visit. Then the result of applying $D$ to $x$, as shown
in Figure \ref{d operator} (b), is to replace $x$ with a negative 2-group.


For the construction of prime alternating knots, it will suffice to
apply $D$ to any one crossing of any negative group or any positive
2-group, but we prove next that if $C(K)$ is a configuration of a
prime alternating knot and $D$ is applied to any crossing of $C(K)$,
the result is a configuration $C(K')$ of a prime alternating knot $K'$
of one higher crossing size.

First of all, observe that the result of applying $D$ is a knot $K'$,
rather than a link, since in effect, the curve was cut on both the
overpassing strand and the underpassing strand at $x$, then the loose
ends rejoined so that with the orientation as shown in Figure \ref{d
operator}, in following a knot traversal we travel from arc $a$ through
one arc of the 2-group to arc $c$, following the portion of the knot
which is delineated by arcs $b$ and $c$ in the reverse order to that
determined by the orientation of $K$, then through the other arc of the
2-group and out along arc $d$, ultimately to return to our starting
point via arc $a$. Moreover, as is shown in Figure \ref{d operator},
the resulting knot is alternating, whence by Proposition \ref{tangle
replacement}, $K'$ is prime.

It is also important to observe that this operation is fully reversible,
in that if we remove any tangle which consists entirely of two adjacent
crossings in a negative group of a prime alternating knot $K'$ as shown
in Figure \ref{d operator} (b), and replace this tangle with a single crossing,
as shown in Figure \ref{d operator} (a), the result is an alternating knot,
hence prime by Proposition \ref{tangle replacement}. We refer to the
inverse of $D$ as $D^-$.

\subsection{The $ROTS$ operator}

\begin{figure}[ht]
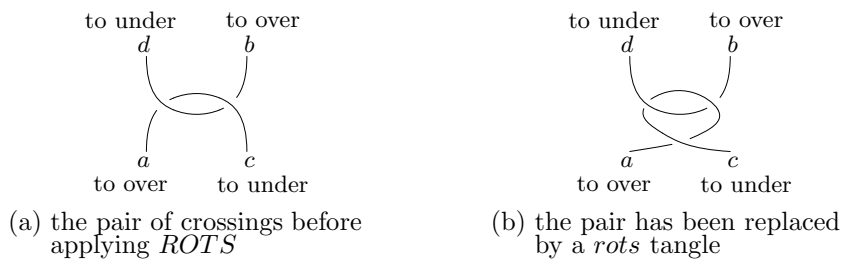

\centering

   \begin{tabular}{c@{\hskip 40pt}c}
   \noalign{\vskip6pt}
$\vcenter{\xy /r38pt/:<0pt,18pt>::,
\vloop <(=.75)\khole,
(0,3),\vloop-  <(=.75)\khole, 
(0,1)="a"*!<1pt,4pt>{\hbox{\small $a$}},
(1,1)="c"*!<-1pt,4pt>{\hbox{\small $c$}},
(0,3)="d"*!<1pt,-5pt>{\hbox{\small $d$}},
(1,3)="b"*!<-1pt,-5pt>{\hbox{\small $b$}},
 "b"*!<-6pt,-14pt>{\hbox{\small to over}},
 "d"*!<6pt,-14pt>{\hbox{\small to under}},
 "a"*!<6pt,12pt>{\hbox{\small to over}},
 "c"*!<-6pt,12pt>{\hbox{\small to under}},  
\endxy}$
&
$\vcenter{\xy /r38pt/:<0pt,18pt>::,
(0,1);(.42,1.15)**\crv{(.3,1.1) },
(.6,1.27);(.21,2.1)**\crv{(.8,1.45) & (1.1,1.8) & (.5, 2.5)},
(.14,1.9);(1,1)**\crv{ (.1,1.7) & (.3,1.4 ) & (.6,1.05) },
(0,3),\vloop-  <(=.75)\khole, 
(0,1)="a"*!<1pt,4pt>{\hbox{\small $a$}},
(1,1)="c"*!<-1pt,4pt>{\hbox{\small $c$}},
(0,3)="d"*!<1pt,-5pt>{\hbox{\small $d$}},
(1,3)="b"*!<-1pt,-5pt>{\hbox{\small $b$}},
 "b"*!<-6pt,-14pt>{\hbox{\small to over}},
 "d"*!<6pt,-14pt>{\hbox{\small to under}},
 "a"*!<6pt,12pt>{\hbox{\small to over}},
 "c"*!<-6pt,12pt>{\hbox{\small to under}},  
\endxy}$
\\
   \noalign{\vskip6pt}
   (a) \vtop{\hsize 1.75in\parindent=0pt\leftskip=0pt
        the pair of crossings before applying $ROTS$}  & (b) \vtop{\hsize 1.75in
       \parindent=0pt\leftskip=0pt the pair has been replaced by
       a $rots$ tangle}\\
   \noalign{\vskip6pt}
\end{tabular}
\caption{The $ROTS$ operator}\label{rots operator}
\end{figure}

This is the second of the two operators that increase the crossing size
by one. $ROTS$ can be applied to any 2-subgroup (that is, any pair of
consecutive crossings) in any
group of two or more. In effect, one of the two arcs joining the two
crossings is cut and the overpassing side of the cut is passed under
the other arc that connects the two crossings, wound around over the arc
and passed back under itself and reconnected to the other cut end. Thus the
result is a knot rather than a link, and, as shown in Figure \ref{rots operator},
the resulting knot is still alternating. Consequently, we can apply
Proposition \ref{tangle replacement} to conclude that we have a configuration of a
prime alternating knot. Since we have replaced a tangle that consisted of
two crossings with one that consists of three crossings, there is a net
increase of one in the crossing size. The tangle shown in Figure
\ref{rots operator} (b) is called a {\em $rots$ tangle}.

Just as for $D$, this process is reversible, and the inverse of the $ROTS$ operator
is denoted by $ROTS^-$.

While the $ROTS$ operator can be applied to any 2-subgroup, the construction of
the prime alternating knots will only use $ROTS$ applied to a negative 2-group,
or to a 2-subgroup of a negative 3-group.
The tangle that results when $ROTS$ is applied to a negative 2-group shall
be called a $ROTS$-tangle.

\subsection{The $OTS$ operator}
\begin{figure}[ht]
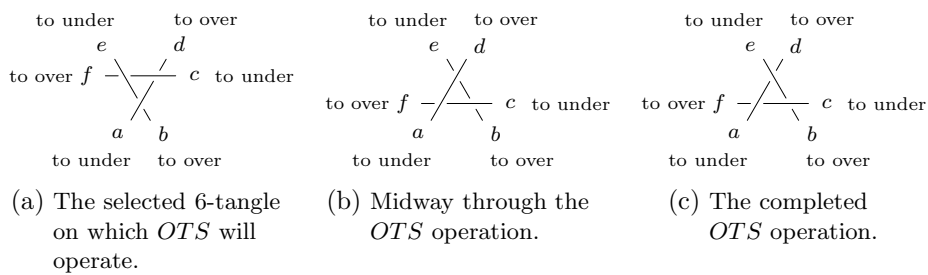

\centering
\begin{tabular}{ccc}
$\vcenter{\xy /l7pt/:,
{\xypolygon3"a"{~>{}~<{}}},
"a1";"a2"**\dir{}?(-.6)="x1"?(-.2)="x2"?(.2)="x3"
    ?(.8)="x4"?(1.2)="x5"?(1.6)="x6",
    "x1";"x4"**\dir{-},
    "x5";"x6"**\dir{-},
"a2";"a3"**\dir{}?(-.6)="y1"?(-.2)="y2"?(.2)="y3"
    ?(.8)="y4"?(1.2)="y5"?(1.6)="y6",
    "y1";"y4"**\dir{-},
    "y5";"y6"**\dir{-},
"a3";"a1"**\dir{}?(-.6)="z1"?(-.2)="z2"?(.2)="z3"
    ?(.8)="z4"?(1.2)="z5"?(1.6)="z6",
    "z1";"z4"**\dir{-},
    "z5";"z6"**\dir{-},
"x1"*!<5pt,5pt>{\hbox{\footnotesize $a$}}*!<15pt,15pt>{\hbox{\scriptsize to under}},
"z6"*!<-5pt,5pt>{\hbox{\footnotesize $b$}}*!<-15pt,15pt>{\hbox{\scriptsize to over}},
"y1"*!<-7pt,0pt>{\hbox{\footnotesize $c$}}*!CL(2){\hbox{\scriptsize to under}},
"x6"*!<-5pt,-6pt>{\hbox{\footnotesize $d$}}*!<-15pt,-15pt>{\hbox{\scriptsize to over}},
"z1"*!<5pt,-5pt>{\hbox{\footnotesize $e$}}*!<15pt,-15pt>{\hbox{\scriptsize to under}},
"y6"*!<7pt,0pt>{\hbox{\footnotesize $f$}}*!CR(2){\hbox{\scriptsize to over}},        
\endxy}$
&
$\vcenter{\xy /r7pt/:,
{\xypolygon3"a"{~>{}~<{}}},
"a1";"a2"**\dir{}?(-.6)="x1"?(-.2)="x2"?(.2)="x3"
    ?(.8)="x4"?(1.2)="x5"?(1.6)="x6",
    "x1";"x6"**\dir{-},
"a2";"a3"**\dir{}?(-.6)="y1"?(-.2)="y2"?(.2)="y3"
    ?(.8)="y4"?(1.2)="y5"?(1.6)="y6",
    "y1";"y2"**\dir{-},
    "y3";"y6"**\dir{-},
"a3";"a1"**\dir{}?(-.6)="z1"?(-.2)="z2"?(.2)="z3"
    ?(.8)="z4"?(1.2)="z5"?(1.6)="z6",
    "z1";"z2"**\dir{-},
    "z3";"z4"**\dir{-},    
    "z5";"z6"**\dir{-},
"x1"*!<-5pt,-5pt>{\hbox{\footnotesize $d$}}*!<-15pt,-15pt>{\hbox{\scriptsize to over}},
"z6"*!<5pt,-5pt>{\hbox{\footnotesize $e$}}*!<15pt,-15pt>{\hbox{\scriptsize to under}},
"y1"*!<7pt,0pt>{\hbox{\footnotesize $f$}}*!CR(2){\hbox{\scriptsize to over}},
"x6"*!<5pt,5pt>{\hbox{\footnotesize $a$}}*!<15pt,15pt>{\hbox{\scriptsize to under}},
"z1"*!<-5pt,5pt>{\hbox{\footnotesize $b$}}*!<-15pt,15pt>{\hbox{\scriptsize to over}},
"y6"*!<-7pt,0pt>{\hbox{\footnotesize $c$}}*!CL(2){\hbox{\scriptsize to under}},            
\endxy}$
&
$\vcenter{\xy /r7pt/:,
{\xypolygon3"a"{~>{}~<{}}},
"a1";"a2"**\dir{}?(-.6)="x1"?(-.2)="x2"?(.2)="x3"
    ?(.8)="x4"?(1.2)="x5"?(1.6)="x6",
    "x1";"x2"**\dir{-},
    "x3";"x6"**\dir{-},    
"a2";"a3"**\dir{}?(-.6)="y1"?(-.2)="y2"?(.2)="y3"
    ?(.8)="y4"?(1.2)="y5"?(1.6)="y6",
    "y1";"y2"**\dir{-},
    "y3";"y6"**\dir{-},
"a3";"a1"**\dir{}?(-.6)="z1"?(-.2)="z2"?(.2)="z3"
    ?(.8)="z4"?(1.2)="z5"?(1.6)="z6",
    "z1";"z2"**\dir{-},
    "z3";"z6"**\dir{-},    
"x1"*!<-5pt,-5pt>{\hbox{\footnotesize $d$}}*!<-15pt,-15pt>{\hbox{\scriptsize to over}},
"z6"*!<5pt,-5pt>{\hbox{\footnotesize $e$}}*!<15pt,-15pt>{\hbox{\scriptsize to under}},
"y1"*!<7pt,0pt>{\hbox{\footnotesize $f$}}*!CR(2){\hbox{\scriptsize to over}},
"x6"*!<5pt,5pt>{\hbox{\footnotesize $a$}}*!<15pt,15pt>{\hbox{\scriptsize to under}},
"z1"*!<-5pt,5pt>{\hbox{\footnotesize $b$}}*!<-15pt,15pt>{\hbox{\scriptsize to over}},
"y6"*!<-7pt,0pt>{\hbox{\footnotesize $c$}}*!CL(2){\hbox{\scriptsize to under}},            
\endxy}$\\
\noalign{\vskip8pt}
(a) \vtop{\hsize=1.25in\leftskip=0pt\noindent\small The selected 6-tangle on which $OTS$ will
operate.}
 & (b) \vtop{\hsize=1.25in\leftskip=0pt\noindent\small Midway through the $OTS$ operation.}
 & (c) \vtop{\hsize=1in\leftskip=0pt\noindent\small The completed $OTS$ operation.}
\end{tabular}
\caption{The $OTS$ operator}\label{ots operator}
\end{figure}

This is the first of the two operators that transform an alternating prime knot into
another of the same crossing size.  It is also the only operator of the four that does
not involve replacing one (4-)tangle by another. A 6-tangle with 3 crossings is called
an $OTS$ 6-tangle, and the $OTS$ operator acts on an $OTS$ 6-tangle, replacing it by another
$OTS$ 6-tangle as illustrated in Figure \ref{ots operator}. Since an $OTS$ 6-tangle has 3
vertices, all of degree 4 in the knot configuration, the sum of the vertex degrees in the
knot configuration is 12. However, by definition, there are 6 edges incident to the 6-tangle,
so the sum of the vertex degrees in the subgraph which is the 6-tangle itself is 6. Thus
there are 3 edges in the 6-tangle, and so an $OTS$ 6-tangle is a triangle, with each crossing
being an endpoint for two  of the 6 incident edges.  In (a), we have shown an $OTS$ 6-tangle,
and (b) and (c) illustrate how an $OTS$ operation is performed on such a 6-tangle.
The operator can be considered to consist of two stages: in the first stage, one of the
three strands $ad$, $be$, or $cf$ is chosen. The chosen strand is cut, moved to the other
side of the crossing formed by the other two strands, and rejoined so as to preserve the
over/under pattern on the strand that has been moved. This stage has been completed in
Figure \ref{ots operator} (b).  At this point, the crossing formed by the other two strands
is now an over-pass when it should be an under-pass or vice-versa, so to complete the $OTS$
operation, we must apply the unknotting surgery to this crossing.
The completed $OTS$ operation is shown in Figure \ref{ots operator} (c).

We remark that for a given $OTS$ 6-tangle, the $OTS$ operation yields the same configuration
(not just flype equivalent) independently of which of the three strands is selected. It is
also evident that $OTS$ is self-inverse.

Since an $OTS$ operation can be viewed as the cutting and rejoining one strand, followed by
the cutting and rejoining of another strand, it follows that the end result is a configuration
of a knot rather than of a link. It remains only to verify that the knot that results is
prime. Suppose to the contrary that for some configuration $C(K)$ of a prime alternating knot
$K$, there is an $OTS$ 6-tangle $O$ in $C(K)$ such that the knot configuration $C(K')$
that results when the $OTS$ operation is performed on $O$ contains a 2-tangle $T$. Let us
denote the three crossings in $O$ by $x$, $y$ and $z$. Since $T\cap O=\emptyset$ would mean
that $T$ is a 2-tangle in $C(K)$, which is not possible since $K$ is prime, we identify three
cases to consider: $T$ contains exactly one, exactly two, and all three crossings of $O$.

\noindent Case 1: $T$ contains exactly one crossing of $O$, say $x$. Then the two edges
that join $a$ to $y$ and to $z$, respectively, are the two edges incident to $T$. If $T$
contains other crossings in addition to $x$, then the other two edges incident to $x$ must
join $x$ to crossings in $T-\{x\}$, whence they are distinct edges and so $T-\{x\}$ is a
2-tangle in $C(K)$ (since the $OTS$ operator did not affect any arcs or crossings of $C(K)$
other than those in $O$). But $C(K)$ has no 2-tangles, so this situation cannot occur.
Consider the possibility that $T$ has $x$ as its only crossing. But then the other two arcs
incident to $x$ must be one and the same, whence $x$ is a nugatory crossing, as shown
in Figure \ref{nugatory in ots} (a). Upon undoing the $OTS$ operation, we find that the
6-tangle $O$ in $C(K)$ must be of the form shown in Figure \ref{nugatory in ots} (b), so
that $O$ in fact was a (4-)tangle, rather than a 6-tangle, so this situation does not arise
either. This completes the proof that Case 1 cannot occur.

\begin{figure}[ht]
\centering
\begin{tabular}{c@{\hskip40pt}c}
$\vcenter{\xy /r7pt/:,
{\xypolygon3"a"{~>{}~<{}}},
"a1";"a2"**\dir{}?(-.6)="x1"?(-.2)="x2"?(.2)="x3"
    ?(.8)="x4"?(1.2)="x5"?(1.6)="x6",
    "x6";"x3"**\dir{-},
"a2";"a3"**\dir{}?(-.6)="y1"?(-.2)="y2"?(.2)="y3"
    ?(.8)="y4"?(1.2)="y5"?(1.6)="y6",
    "y1";"y2"**\dir{-},
    "y3";"y6"**\dir{-},
"a3";"a1"**\dir{}?(-.6)="z1"?(-.2)="z2"?(.2)="z3"
    ?(.8)="z4"?(1.2)="z5"?(1.6)="z6",
    "z1";"z2"**\dir{-},
    "z3";"z5"**\dir{-},
    "a1"+(.75,0)*=0{\hbox{\footnotesize $x$}},
    "x2";"z5"**\crv{"x2"+(.75,.75) & "a1"+(0,2.75) & "z5"+(-.75,.75)},
    "a1"+(0,1.25)*\xycircle(2,2){{.}},
    "x1"*!(-2.25,-1.25){\hbox{\footnotesize $T$}},
"x1"*!<-5pt,-5pt>{\hbox{\footnotesize $d$}},
"z6"*!<5pt,-5pt>{\hbox{\footnotesize $e$}},
"y1"*!<7pt,0pt>{\hbox{\footnotesize $f$}},
"x6"*!<5pt,5pt>{\hbox{\footnotesize $a$}},
"z1"*!<-5pt,5pt>{\hbox{\footnotesize $b$}},
"y6"*!<-7pt,0pt>{\hbox{\footnotesize $c$}},
\endxy}$
&
$\vcenter{\xy /l7pt/:,
{\xypolygon3"a"{~>{}~<{}}},
"a1";"a2"**\dir{}?(-.6)="x1"?(-.2)="x2"?(.2)="x3"
    ?(.8)="x4"?(1.2)="x5"?(1.6)="x6",
    "x1";"x4"**\dir{-},
"a2";"a3"**\dir{}?(-.6)="y1"?(-.2)="y2"?(.2)="y3"
    ?(.8)="y4"?(1.2)="y5"?(1.6)="y6",
    "y1";"y4"**\dir{-},
    "y5";"y6"**\dir{-},
"a3";"a1"**\dir{}?(-.6)="z1"?(-.2)="z2"?(.2)="z3"
    ?(.8)="z4"?(1.2)="z5"?(1.6)="z6",
    "z3";"z4"**\dir{-},
    "z5";"z6"**\dir{-},
    "z3";"x5"**\crv{ "z3"+(.9,-1.54) & "a2"+(0,-2) & "x5"+(-.2,-.32)},
"x1"*!<5pt,5pt>{\hbox{\footnotesize $a$}},
"z6"*!<-5pt,5pt>{\hbox{\footnotesize $b$}},
"y1"*!<-7pt,0pt>{\hbox{\footnotesize $c$}},
"x6"*!<-5pt,-6pt>{\hbox{\footnotesize $d$}},
"z1"*!<5pt,-5pt>{\hbox{\footnotesize $e$}},
"y6"*!<7pt,0pt>{\hbox{\footnotesize $f$}},
\endxy}$\\
\noalign{\vskip8pt}
(a) & (b)
\end{tabular}
\caption{Case 1}\label{nugatory in ots}
\end{figure}

\noindent Case 2: $T$ contains exactly two crossings of $O$, say $x$ and $y$. But then the
arc from $x$ to $z$ and the arc from $y$ to $z$ must be the two arcs of $C(K')$ that are
incident to $T$. We may then include $z$ with $T$ to form a larger 2-tangle in $C(K')$ which
contains all three of $x$, $y$ and $z$. This situation is our third case.

\noindent Case 3: $T$ contains all three crossings $x$, $y$ and $z$. In this case, undoing the
$OTS$ operation produces a 2-tangle $T'$ in $C(K)$ which contains the $OTS$ 6-tangle $O$.
But since $K$ is prime, $C(K)$ has no 2-tangles, so we obtain a contradiction in this case
as well.

It follows now that $C(K')$ can't contain any 2-tangles, and therefore $K'$ is a prime knot.

\subsection{The $T$ operator}
 Given a configuration $C(K)$ of a prime alternating knot $K$, a turnable-subgroup
 of $C(K)$ is a subgroup of any size greater than or equal to 2 of a positive
 group of $C(K)$, or a subgroup of odd size greater than or equal to 3 of any
 negative group of $C(K)$. The $T$ operator acts on a turnable-subgroup of $C(K)$
 via a two stage process. First, label the arcs incident to the turnable-subgroup
 by the following device: draw a simple closed curve surrounding the subgroup,
 then follow this curve in the clockwise direction until the first of a pair
 of edges at one end of the subgroup is encountered. Place the label 1 on
 this edge, and then continue following this curve, placing the labels 2, 3
 and 4 on the remaining three arcs incident to the subgroup as they are
 encountered. Now at each label $c$, place a second label $c'$
 on the arc situated between the position of the label $c$ and the crossing
 of the subgroup that the arc meets. Cut each arc between the labels
 $c$ and $c'$, then reconnect the cut strands so as to have the arc labelled
 1 connect to the arc labelled $2'$, the arc labelled 2 connect to the arc
 labelled $3'$, the arc labelled 3 connect to the arc labelled $4'$ and
 the arc labelled 4 connect to the arc labelled $1'$. The effect is to
 have turned the subgroup one-quarter turn in the counter-clockwise
 direction. However, arcs 1 and 2 originally met at a crossing at one
 end of the subgroup, so one was an overpass while the other was an
 underpass. After the turn, arc 1 will have the overpass/underpass status
 that was formerly held by arc 2, which means that the resulting
 configuration is not alternating. Observe that the arcs labelled 1 and
 3, respectively, have the same overpass/underpass status at the crossings
 of the subgroup that they are incident to, as do arcs labelled 2 and 4.
 Thus, after the turn, arc 2 has the overpass/underpass status formerly
 held by arc 3, which is the same as that of the original arc 1 and therefore
 opposite to that of the original arc 2. The same observation holds for
 all four arcs. Thus we can make an alternating configuration by cutting
 one strand at each crossing of the subgroup, moving the strand across the
 uncut strand and rejoin the cut ends so as to change the cut strand from
 an underpassing strand to an overpassing strand and vice-versa. If the
 configuration after the turn was a knot rather than a link, then after
 completing the second phase of the operation, the result is a configuration
 of an alternating knot. By Proposition 1, if we apply $T$ to a turnable
 subgroup of some prime alternating knot, the result is a configuration
 of a prime alternating knot.

 As it turns out, the $T$ operator could be applied to an even subgroup of
 a negative group, but the result would be a link. We establish now that
 when $T$ is applied to a turnable-subgroup, the result is a configuration of
 an alternating knot, hence by Proposition \ref{tangle replacement}, a configuration
 of a prime alternating knot.

\noindent Case 1: the turnable-subgroup is positive. Let $k$
 denote the number of crossings in the turnable-subgroup.

 \noindent Case 1 (a): $k$ is even. Then arcs $1'$ and $3'$ are on the same strand of the
 subgroup (after the four incident arcs have been cut), and arcs $2'$ and
 $4'$ are on the other strand. Since the subgroup is positive, when the
 subgroup has been cut from the configuration, the remainder of the knot
 decomposes into two strands, and the strand which contains arc 3 is the
 same strand which contains arc 2, while the strand that contains arc 4
 also contains arc 1. When the turnable-subgroup is turned and reattached, a
 traversal of the configuration beginning at arc 1 in the direction of
 the reattached subgroup follows arc $2'$ to arc $4'$, which is connected
 to arc 3, so we continue on to arc 2, which is connected to arc $3'$.
 Continue the traversal to arc $1'$, which is attached to arc 4, so we are
 led back to arc 1. Every arc has been traversed, so the result of
 applying $T$ to an even subgroup of a positive group is a knot configuration.

\begin{figure}[ht]
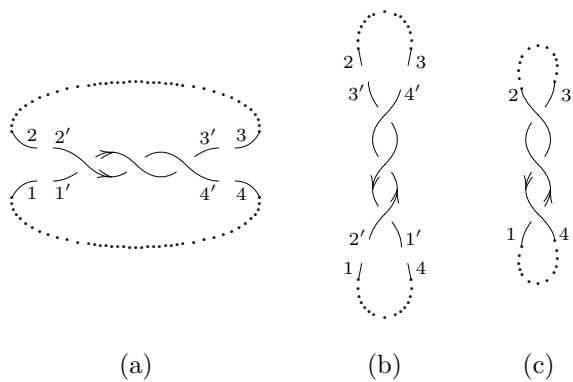

\centering
\begin{tabular}{c@{\hskip30pt}c@{\hskip30pt}c}
 $\vcenter{\xy /r25pt/:,
 (1,0)="c5",(1.75,0)="c6",
 (.5,.25)="b1",(.5,-.25)="b2",(3,-.25)="b4",
 (3,.25)="b3",(2.5,0)="c7",
 "b1";"c5"**\crv{(.75,.25)},
 "c5";"c6"*+{\hbox to 1pt {\hfil}}**\crv{(1.25,-.25) & (1.5,-.25)}?(.5)*\dir{}="c"-<8pt,-2pt>;"c"**\dir{}*\dir{>},
 "c6"*+{\hbox to 1pt {\hfil}};"c7"**\crv{(2,.25) & (2.25,.25)},
 "c7";"b4"**\crv{(2.75,-.25)},
 "b2";"c5"*+{\hbox to 1pt {\hfil}}**\crv{(.75,-.25)},
 "c5"*+{\hbox to 1pt {\hfil}};"c6"**\crv{(1.25,.25) & (1.5,.25)}?(.55)*\dir{}="c"-<8pt,1.25pt>;"c"**\dir{}*\dir{>},
 "c6";"c7"*+{\hbox to 1pt {\hfil}}**\crv{(2,-.25) & (2.25,-.25)},
 "c7"*+{\hbox to 1pt {\hfil}};"b3"**\crv{(2.75,.25)},
 (.3,0)="q1",(3.2,0)="q3",
 (1.75,.75)="q2",
 (1.75,-.75)="q4",
 (-.125,.5)="t1",(3.625,.5)="t3",
 (.25,.25)="t2",
 (3.25,.25)="t4",
 (-.125,-.5)="t5",(.25,-.25)="t6",
 (3.25,-.25)="t7",
 (3.625,-.5)="t8",
 "t1";"t2"**\crv{(0,.25)},
 "t4";"t3"**\crv{(3.5,.25)},
 "t5";"t6"**\crv{(0,-.25)},
 "t7";"t8"**\crv{(3.5,-.25)},
 "t1";"t3"**\crv{~*=<\jot>{.}(-.25,1.5) & (3.75,1.5)},
 "t5";"t8"**\crv{~*=<\jot>{.}(-.25,-1.5) & (3.75,-1.5)},
 (.18,-.43)*{\hbox{\scriptsize 1}},
 (.18,.45)*{\hbox{$\ssize 2$}},
 (3.35,.45)*{\hbox{$\ssize 3$}},
 (3.35,-.44)*{\hbox{$\ssize 4$}},
 (.65,-.41)*{\hbox{$\ssize 1'$}},
 (.65,.44)*{\hbox{$\ssize 2'$}},
 (2.85,.425)*{\hbox{$\ssize 3'$}},
 (2.85,-.425)*{\hbox{$\ssize 4'$}},
 \endxy}$
 &
 $\vcenter{\xy /r25pt/:,
 (0,1)="c5",(0,1.75)="c6",
 (.25,.5)="b1",(-.25,.5)="b2",(-.25,3)="b4",
 (.25,3)="b3",(0,2.5)="c7",
 "b1";"c5"*+{\hbox to 1pt {\hfil}}**\crv{(.25,.75)},
 "c5"*+{\hbox to 1pt {\hfil}};"c6"**\crv{(-.25,1.25) & (-.25,1.5)}?(.5)*\dir{}="c"+<2pt,10pt>;"c"**\dir{}*\dir{>},
 "c6";"c7"*+{\hbox to 1pt {\hfil}}**\crv{(.25,2) & (.25,2.25)},
 "c7"*+{\hbox to 1pt {\hfil}};"b4"**\crv{(-.25,2.75)},
 "b2"*+{\hbox to 1pt {\hfil}};"c5"**\crv{(-.25,.75)},
 "c5";"c6"*+{\hbox to 1pt {\hfil}}**\crv{(.25,1.25) & (.25,1.5)}?(.5)*\dir{}="c"-<2pt,10pt>;"c"**\dir{}*\dir{>},
 "c6"*+{\hbox to 1pt {\hfil}};"c7"**\crv{(-.25,2) & (-.25,2.25)},
 "c7";"b3"*+{\hbox to 1pt {\hfil}}**\crv{(.25,2.75)},
 (0,.25)="q1",(0,3.25)="q3",
 (.75,1.75)="q2",
 (-.75,1.75)="q4",
 (.4,0)="t1",(.4,3.5)="t3",
 (.35,.25)="t2",
 (.35,3.25)="t4",
 (-.4,0)="t5",(-.35,.25)="t6",
 (-.35,3.25)="t7",
 (-.4,3.5)="t8",
 "t1";"t2"**\dir{-},
 "t4";"t3"**\dir{-},
 "t5";"t6"**\dir{-},
 "t7";"t8"**\dir{-},
 "t1";"t5"**\crv{~*=<\jot>{.}(.5,-.75) & (-.5,-.75)},
 "t3";"t8"**\crv{~*=<\jot>{.}(.5,4.25) & (-.5,4.25)},
 (-.55,.18)*{\hbox{$\ssize 1$}},
 (.55,.18)*{\hbox{$\ssize 4$}},
 (.55,3.3)*{\hbox{$\ssize 3$}},
 (-.55,3.3)*{\hbox{$\ssize 2$}},
 (-.41,.65)*{\hbox{$\ssize 2'$}},
 (.44,.65)*{\hbox{$\ssize 1'$}},
 (.425,2.85)*{\hbox{$\ssize 4'$}},
 (-.425,2.85)*{\hbox{$\ssize 3'$}},
 \endxy}$
 &
 $\vcenter{\xy /r25pt/:,
 (0,1)="c5",(0,1.75)="c6",
 (.25,.6)="b1",(-.25,.5)="b2",(-.25,2.9)="b4",
 (.25,3)="b3",(0,2.5)="c7",
 "b1";"c5"**\crv{(.25,.75)},
 "c5";"c6"*+{\hbox to 1pt {\hfil}}**\crv{(-.25,1.25) & (-.25,1.5)}?(.5)*\dir{}="c"+<2pt,10pt>;"c"**\dir{}*\dir{>},
 "c6"*+{\hbox to 1pt {\hfil}};"c7"**\crv{(.25,2) & (.25,2.25)},
 "c7";"b4"**\crv{(-.25,2.75)},
 "b2";"c5"*+{\hbox to 1pt {\hfil}}**\crv{(-.25,.75)},
 "c5"*+{\hbox to 1pt {\hfil}};"c6"**\crv{(.25,1.25) & (.25,1.5)}?(.5)*\dir{}="c"-<2pt,10pt>;"c"**\dir{}*\dir{>},
 "c6";"c7"*+{\hbox to 1pt {\hfil}}**\crv{(-.25,2) & (-.25,2.25)},
 "c7"*+{\hbox to 1pt {\hfil}};"b3"**\crv{(.25,2.75)},
 "b1";"b2"**\crv{~*=<\jot>{.}(.5,-.25) & (-.5,-.25)},
 "b3";"b4"**\crv{~*=<\jot>{.}(.5,3.75) & (-.5,3.75)},
 (-.41,.7)*{\hbox{$\ssize 1$}},
 (.39,.7)*{\hbox{$\ssize 4$}},
 (.425,2.8)*{\hbox{$\ssize 3$}},
 (-.375,2.8)*{\hbox{$\ssize 2$}},
 \endxy}$
\\
\noalign{\vskip8pt}
(a) & (b) & (c)
\end{tabular}
\caption{Case 1 (b)}\label{turning odd subgroup}
\end{figure}

 \noindent Case 1 (b): $k$ is odd. Again, when the turnable-subgroup has been cut
 from the knot configuration, the resulting four strands have
 arcs $1'$ and $3'$ on one, arcs $2'$ and $4'$ on another, arcs 1 and 4
 on a third, and arcs 2 and 3 on the last strand. When the subgroup is
 turned and reattached, we begin a traversal at arc 1 and move towards
 the subgroup. We encounter in the following order arcs $2'$, $4'$, 3,
 2, $3'$, $1'$, 4 and finally we arrive back at 1. Again, we have
 traversed every arc, so the result is a knot configuration. An example
 of this case is shown in Figure \ref{turning odd subgroup}.

\noindent Case 2: the turnable-subgroup is negative. By definition, the number
$k$ of crossings is odd and $k\ge3$. When the turnable-subgroup is cut out
from the knot configuration, arcs $1'$ and $3'$ are on one strand, arcs
$2'$ and $4'$ are on a second strand, while arcs 3 and 4 are on another strand
and arcs 1 and 2 are on the last strand. When the turnable-subgroup is
turned and reattached, a traversal beginning at arc 1 and proceeding in
the direction of the subgroup will encounter arcs in the order $2'$,
$4'$, 3, 4, $1'$, $3'$, 2 and back to 1, so the entire configuration has
been traversed.  The result in this case as well is a knot configuration.

Finally, observe that when $T$ is applied to a subgroup of odd size,
the turned subgroup has the opposite sign to that of the original (that
is, a positive subgroup becomes negative, and a negative subgroup becomes
positive) while turning a even sized positive subgroup results in a
positive subgroup. Consequently, the result of turning a turnable-subgroup is
still a turnable-subgroup. If $T$ is then applied to this resulting
turnable-subgroup, the original knot configuration is produced, whence $T$
is its own inverse.

\par
\section{The Construction Makes Them All}

In this section, we shall prove that if all possible applications of both
$D$ and $ROTS$ are carried out on the complete set of prime alternating 
knots of $n$ crossings, and then all possible applications of both $OTS$ 
and $T$ are carried out on the set of resulting knots of $n+1$ crossings,
then again to the knots that result from this first application of $OTS$
and $T$, and so on until no further knots are produced, the result is 
the set of all prime alternating knots of $n+1$ crossings. In fact, we 
shall prove that the same outcome results if we only apply $D$ to negative 
groups (including loners) and positive 2-groups, $ROTS$ to negative 2-groups
and negative 3-groups, and $T$ to positive 2-groups. There will be no 
restrictions on $OTS$.

It will prove to be convenient to let $\teetwo$ denote the restriction of $T$ to
only positive 2-groups.

To aid us in the proof of this result, we introduce three categories of prime alternating
knots, which have the property that each prime alternating knot belongs to exactly one category.

\subsection{The $K_A$ configurations}

Conceptually, a prime alternating knot configuration is a $K_A$ configuration
if it can be obtained from a prime alternating knot configuration of one lower
crossing size by a single application of either the $D$ or the $ROTS$ operator (under
the restrictions listed at the beginning of this section).

\begin{definition}\label{ka configuration}
 A prime alternating knot configuration $C(K)$ is said to be a $K_A$
 configuration if it contains either a negative group of size 2 or greater,
 or a $rots$-tangle.
\end{definition}

\begin{theorem}\label{characterize ka configurations}
  A prime alternating knot configuration is a $K_A$ configuration
  if it can be obtained from a prime alternating knot configuration of one lower
  crossing size by a single application of either $D$ to a negative group,
  a loner, or a positive 2-group, or the $ROTS$ operator to a negative
  2-group or 3-group.
\end{theorem}

\begin{proof}
Let $C(K)$ be a $K_A$ configuration. Suppose first of all that $C(K)$ contains
a negative group $G$ of size 2 or greater. Then an application of $D^-$ to $G$
results in a configuration $C(K')$ of a prime alternating knot at one lower
crossing size containing a negative group $G'$ with the property that if $D$ is
applied to $G'$, the configuration $C(K)$ is reconstructed. Now suppose that
$C(K)$ contains a $rots$-tangle $T$. Then an application of $ROTS^-$ to $T$
results in a configuration $C(K')$ of a prime alternating knot at one lower
crossing size which contains a negative 2-subgroup $T'$ such that if $ROTS$
is applied to $T'$, $C(K)$ is obtained. If $T'$ is actually a subgroup of
a negative group of size at least 4, in which case $T$ is directly attached
to a negative 2-group $G$, and we have already handled the case when $C(K)$
contains a negative group of size 2 or greater.
\end{proof}

Thus each $K_A$ configuration at a given number of crossings is obtained
when $D$ (applied only to negative groups, loners, and positive 2-groups)
and $ROTS$ (applied only to negative 2 and 3-groups) are applied to the
configurations of the prime alternating knots of one lower crossing size.

\subsection{The $K_B$ configurations}
 To describe the $K_B$ configurations, we require the following concepts.

 \begin{definition}
  Let $C(K)$ be a configuration of a prime alternating knot $K$. An
  {\em interleaved 2-sequence} in $C(K)$ is a pair of positive 2-subgroups
  of $C(K)$ having no crossings in common with the property that, as the knot
  is traversed, one arc of one of the 2-subgroups is traversed first, then one
  arc of the other 2-subgroup is traversed next, then the second arc of the first
  2-subgroup is traversed and finally the
  second arc of the second 2-subgroup is traversed.

 \begin{figure}[ht]
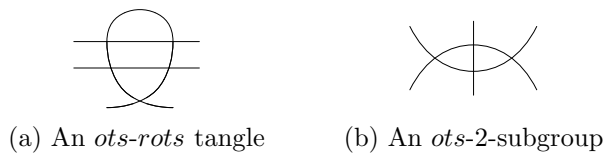

  \centering
    \begin{tabular}{c@{\hskip30pt}c}
    $\vcenter{\setbox0=\hbox{\xy /r25pt/:,
    \vcap,c="a",\vover-,"a",\vunder-,
    (.5,0);(2.4,0)**\dir{-},
    (.5,-.4);(2.4,-.4)**\dir{-},
    \endxy}\dimen0=\ht0\ht0=.5\dimen0\box0}$
    &
    $\vcenter{
    \xy /r20pt/:,
     (0,-.1);(2.4,-.1)**\crv{ (.3,.5) & (1.2, 1) & (2.1,.5)},
     (0,1.1);(2.4,1.1)**\crv{ (.3,.5) & (1.2, 0) & (2.1,.5)},
     (1.2,-.2);(1.2,1.2)**\dir{-},
    \endxy}$\\
    \noalign{\vskip 6pt}
    (a) An $ots$-$rots$ tangle & (b) An $ots$-2-subgroup
   \end{tabular}
   \caption{$ots$ type tangles}
   \label{ots-type 6-tangles}
  \end{figure}

  \noindent
  Furthermore, a 6-tangle of the form shown in Figure
  \ref{ots-type 6-tangles} (a) is called an {\em $ots$-$rots$ tangle},
  while a 6-tangle of the form shown in Figure \ref{ots-type 6-tangles}
  (b) is called an {\em $ots$-2-subgroup}, with sign that of the underlying
  2-subgroup. If the underlying 2-subgroup is actually a group, then the
  6-tangle is called an {\em $ots$-2-group}.
\end{definition}

We observe that any negative group of size 4 or greater provides an instance
of an interleaved 2-sequence.

\begin{definition}\label{kb configuration} 
  A prime alternating knot configuration $C(K)$ is said to be a $K_B$
  configuration if it is not a $K_A$ configuration (that is, contains
  no negative groups of size 2 or greater, and contains no $rots$-tangles),
  but contains one or more of the following structures: a positive group
  of size 3, an interleaved 2-sequence, an $ots$-$rots$ tangle, or a
  negative $ots$-2-group.
\end{definition}

  We shall prove that for each $K_B$ configuration
  $C(K)$, there is a $K_A$ configuration $C(K')$ such that there is either
  a single $\teetwo$ operation or a single $OTS$ operation which transforms
  $C(K')$ into $C(K)$, and that furthermore, if either $\teetwo$ or $OTS$
  is applied to a $K_A$ configuration and the result is not a $K_A$
  configuration, then it is a $K_B$ configuration.

\begin{figure}[ht]
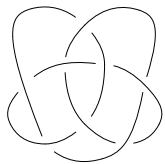

  \centering
   \begin{tabular}{cc}
    $\vcenter{
    \xy /r30pt/:,
     (1,1.75)="a1",(1.75,1)="a2",
     (1.5,.25)="a3",(.5,.25)="a4",(.25,1)="a5",
     (1.25,1.25)="a6",(1,.5)="a7",(.75,1.25)="a8",
     "a1";"a2"*+{\hbox to 1pt {\hfil}}**\crv{(1.25,2) &
                                         (1.75,2) & (2,1.75)},
     "a2"*+{\hbox to 1pt {\hfil}};"a4"*+{\hbox to 1pt {\hfil}}
                       **\crv{(1.625,.25) & (1.25,0) & (.75,0)},
     "a4"*+{\hbox to 1pt {\hfil}};"a1"*+{\hbox to 1pt {\hfil}}
                       **\crv{(.2,1) & (0,1.75) & (.25,2) & (.75,2)},
     "a1"*+{\hbox to 1pt {\hfil}};"a6"**\crv{(1.25,1.5)},
     "a6";"a7"*+{\hbox to 1pt {\hfil}}**\crv{(1.25,.75)},
     "a7"*+{\hbox to 1pt {\hfil}};"a4"**\crv{(.75,.25)},
     "a4";"a5"*+{\hbox to 1pt {\hfil}}**\crv{(.25,.25) & (0,.5) & (0,.75)},
     "a5"*+{\hbox to 1pt {\hfil}};"a6"*+{\hbox to 1pt {\hfil}}
                       **\crv{(.5,1.25) & (.75,1.325)},
     "a6"*+{\hbox to 1pt {\hfil}};"a2"**\crv{(1.5,1.25)},
     "a2";"a3"*+{\hbox to 1pt {\hfil}}**\crv{(2,.75) & (2,.5) & (1.75,.25)},
     "a3"*+{\hbox to 1pt {\hfil}};"a7"**\crv{(1.25,.25)},
     "a7";"a8"*+{\hbox to 1pt {\hfil}}**\crv{(.75,.75)},
     "a8"*+{\hbox to 1pt {\hfil}};"a1"**\crv{(.75,1.5)},
    \endxy}$
    \end{tabular}
   \caption{Knot $K_{8_{16}}$}
   \label{8_16}
  \end{figure}

  By way of illustration, observe that knot $5_1$ is a $K_B$ knot of five
  crossings, since it has both a positive 5-group and an interleaved 2-sequence,
  but no negative groups nor any $rots$-tangles. Knot $8_{16}$ (see Figure \ref{8_16})
  also qualifies as a $K_B$ knot of eight crossings on several counts, since it
  has an interleaved 2-sequence, an $ots$-$rots$ tangle and a negative
  $ots$-2-group. Knot $8_{16}$ does not have any positive groups of size greater
  than 2, and is not a $K_A$ knot since it has no negative groups nor does it
  have any $rots$-tangles.

\begin{proposition}
  Let $C(K)$ be a $K_B$ configuration. Suppose that $C(K)$ contains a
  2-subgroup $G$, and that $C(K')$ denotes the configuration that is
  obtained when $T$ is applied to a 2-subgroup of $G$.  Then the image
  of $G$ will be a 2-group of the resulting configuration.
\end{proposition}  

\begin{proof}
  Suppose that the image $G'$ of $G$ upon this application
  of $T$ is a subgroup of a larger group, as illustrated in Figure \ref{turning
  a 2-subgroup} (a).  Since $G'$ is a subgroup of a larger group in $C(K')$, we
  see that $G$ formed part of a $rots$ tangle, as illustrated in Figure \ref{turning
  a 2-subgroup} (b). But $C(K)$ is not a $K_A$ configuration,
  and so $C(K)$ does not have any $rots$ tangles. Thus $G'$ could not be a
  subgroup of a larger group in $C(K')$.
\end{proof}

  Thus, if the  $T$ operator is applied to any 2-subgroup of a $K_B$ configuration, the
  image of the 2-subgroup will be a 2-group of the resulting configuration, whence
  $T$ can be applied to it to regain $G$.

 \begin{figure}[ht]
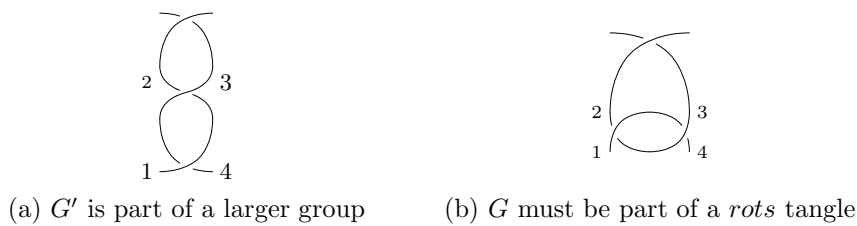

\centering
   \begin{tabular}{c@{\hskip30pt}c}
    $\vcenter{\xy /r20pt/:,
    {\vunder,c="a",\vtwist\vunder-,c="b"},
    "a"*!<5pt,6pt>{\hbox{\scriptsize 2}},
    "a"+(1,0)*!<-5pt,6pt>{\hbox{\small 3}},
    "b"*!<5pt,0pt>{\hbox{\small 1}},    
    "b"+(1,0)*!<-5pt,0pt>{\hbox{\small 4}},
    \endxy}$
    &
    $\vcenter{
    \xy /r30pt/:,
     {\vunder,c="a",\hover[.5],c-(.5,0),\hover[-.5],c="b"},
     "a"*!<5pt,0pt>{\hbox{\scriptsize 2}},
     "a"+(0,-.5)*!<5pt,0pt>{\hbox{\scriptsize 1}},
     "a"+(1,0)*!<-5pt,0pt>{\hbox{\scriptsize 3}},
     "a"+(1,-.5)*!<-5pt,0pt>{\hbox{\scriptsize 4}},  
    \endxy}$\\
   \noalign{\vskip 6pt}
   (a) $G'$ is part of a larger group & (b) $G$ must be part
    of a $rots$ tangle
  \end{tabular}
  \caption{Turning a 2-subgroup}
  \label{turning a 2-subgroup}
 \end{figure}

\begin{theorem}\label{characterize kb configurations}
  A prime alternating knot configuration $C(K)$ is a $K_B$ configuration
  if and only if it is not a $K_A$ configuration, but can be obtained from
  some $K_A$ configuration by a single application of either $\teetwo$ or
  $OTS$.
\end{theorem}

\begin{proof}
  Let $C(K)$ be a $K_B$ configuration, whence by definition, $C(K)$ is not
  a $K_A$ configuration. We must consider four cases.
  
  \noindent Case 1: $C(K)$ contains
  a positive 3-group $G$. Choose a 2-subgroup of $G$ and apply $T$ to the
  selected subgroup. Let $C(K')$ denote the configuration that results, and
  let $H$ denote the positive 2-group which is the image of $T$ applied
  to the selected 2-subgroup of $G$. Then $H$ together with the remaining
  crossing of the original 3-group forms an $ots$-tangle in $C(K')$, whence
  $C(K')$ is a $K_A$ configuration.  If $\teetwo$ is applied to the positive 2-group
  $H$ in $C(K')$, the result is the original configuration $C(K)$, and so in
  this case, $C(K)$ has been obtained by a single application of $\teetwo$ to 
  some $K_A$ configuration.

\noindent Case 2: $C(K)$ contains an interleaved 2-sequence, say $G_1$ and
  $G_2$. Label the arcs incident to the tangle $G_1$ with 1, 2, 3 and 4 in
  preparation for an application of $T$ to $G_1$, and begin a traversal of
  the knot with arc 1 moving towards $G_1$. Label the arcs of $G_1$ and $G_2$
  with $e_1, e_2$ and $f_1,f_2$, respectively, so that the traversal of the knot
  encounters the arcs in the sequence  $1,e_1,4,f_1,2,e_2,3,f_2,1$, so the fact
  that $G_2$ is a positive subgroup is recorded by the data $4,f_1,2$ and $3,f_2,1$
  in the traversal. Consider what happens when $T$ is applied to $G_1$. A
  traversal of the resulting knot that starts on arc 1 will encounter the
  arcs of $G_1$ (turned) and $G_2$ in the sequence $1,e_2,2,f_1,4,e_1,3,f_2,1$,
  so the order of travel through $G_2$ is $2,f_1,4$ and $3,f_2,1$. Thus
  $f_2$ is traversed in the same direction as in the original knot configuration,
  while $f_1$ is now traversed in the opposite direction to that of the original
  knot. Since $G_2$ was originally a positive group, in the configuration $C(K')$
  that results from the application of $T$ to $C_1$, $G_2$ has become a negative
  group, whence $C(K')$ is a $K_A$ configuration. Now, the turned image $G_1'$ of
  $G_1$ is a  positive 2-group of $C(K')$, and so we may apply $\teetwo$ to obtain
  $C(K)$. Thus in this case as well, $C(K)$ is obtained by applying $\teetwo$ to
  some $K_A$ configuration.

  \noindent Case 3: $C(K)$ contains an $ots$-$rots$ tangle $O$. Since
  every $ots$-$rots$ tangle contains an $ots$ tangle, we choose to apply
  $OTS$ to the $ots$ tangle that is contained in $O$. Let $C(K')$ denote
  the resulting configuration. Now, the application of $OTS$ to $O$
  leaves a $rots$ tangle in $C(K')$, so $C(K')$ is a $K_A$ configuration.
  If we apply $OTS$ to the $ots$ tangle of $C(K')$ that is the image of
  $O$ under the initial $OTS$ operation, the result is $C(K)$. Thus every
  $K_B$ configuration that contains an $ots$-$rots$ tangle is obtained by
  a single application of $OTS$ to some $K_A$ configuration.

  \noindent Case 4: $C(K)$ contains a negative $ots$-2-group. Now any
  $ots$ 2-group contains an $ots$ tangle, so we apply $OTS$ to the
  $ots$ tangle $O$ that is contained in the negative $ots$-2-group.
  The result is a configuration $C(K')$ containing a negative 2-group,
  which is therefore a $K_A$ configuration. When $OTS$ is applied to the
  image of $O$ in $C(K')$, the original configuration $C(K)$ is obtained,
  Thus any $K_B$ configuration that contains a negative $ots$-2-group
  is the image under a single $OTS$ operation of some $K_A$ configuration.

  This establishes that for every $K_B$ configuration, there is a $K_A$
  configuration such that a single application of either $\teetwo$ or $OTS$
  to this $K_A$ configuration produces the given $K_B$ configuration.

  Conversely, suppose that $C(K)$ is not a $K_A$ configuration, but that
  there is either a $K_A$ knot configuration $C(K_1)$ such that an
  application of $\teetwo$ to $C(K_1)$ produces $C(K)$, or else there is
  a $K_A$ configuration $C(K_)$ such that an application of $ROTS$ to
  $C(K_1)$ produces $C(K)$. We show that in either case, $C(K)$ is a
  $K_B$ configuration.

  \noindent {\bf Case 1:} an application of $\teetwo$ to some positive
  2-group $G$ in the $K_A$ configuration $C(K_1)$ produces $C(K)$. Since
  $C(K)$ is not a $K_A$ configuration, it has no negative groups or $rots$
  tangles, while $C(K_1)$ is a $K_A$ configuration, and so it must have
  at least one negative group or $rots$ tangle. If $C(K_1)$ contains a
  $rots$ tangle, then $\teetwo$ must have been applied to the 2-group of
  the $rots$ tangle. Since an application of $\teetwo$ to the 2-group of
  a $rots$ tangle results in a positive 3-group, it follows that $C(K)$ is
  a $K_B$ configuration in this case. If $C(K_1)$ does not contain any
  $rots$ tangles, then it must contain a negative group $H$. But $C(K)$
  does not contain any negative groups, so the result of applying $\teetwo$
  to $G$ is to change $G$ into a positive 2-group. But this means that
  $G$ and $H$ constitute an interleaved-2-sequence, and so in this case
  as well, $C(K)$ is a $K_B$ configuration.

  \noindent {\bf Case 2:} an application of $OTS$ to some $ots$ 6-tangle
  $O$ in the $K_A$ configuration $C(K_1)$ produces $C(K)$. Since $C(K_1)$
  is a $K_A$ configuration, it contains either a negative group or a
  $rots$ tangle. If $C(K_1)$ contains a $rots$ tangle, then since $C(K)$
  does not, we see that the $ots$ 6-tangle $O$ must involve crossings of
  the $rots$ tangle. By definition of $OTS$, $O$ cannot contain both
  crossings of the positive 2-group of the $rots$ tangle. If $O$ contained
  one crossing of the positive 2-group of the $rots$ tangle, then $O$
  must contain the third crossing of the $rots$ tangle as well, but in
  such a case, the application of $OTS$ to $O$ would result in a negative
  2-group in $C(K)$. Since this is not possible, $O$ cannot contain either
  crossing of the positive 2-group of the $rots$ tangle. But this means
  that $O$ must contain the third crossing of the $rots$ tangle, and so
  the result of applying $OTS$ to $O$ is an $ots$-$rots$ tangle in $C(K)$,
  whence $C(K)$ is a $K_B$ configuration.

  The last possibility is that $C(K_1)$ contains a negative group $G$.
  Since $C(K)$ contains no negative groups, the $ots$ 6-tangle $O$ must
  involve a crossing at one end of $G$, in which case the application
  of $OTS$ to $O$ will separate this crossing from the rest of $G$. If
  $G$ had been a negative group of size 3 or larger, then $C(K)$ would
  contain a negative group of size 2 or larger, which is not possible.
  Thus $G$ must be a negative 2-group and so $OTS$ results in a negative
  $ots$-2-group. Again, $C(K)$ is a $K_B$ configuration.

  Thus in every case, $C(K)$ is a $K_B$ configuration.
\end{proof}

\subsection{The $K_C$ configurations}
 These configurations make up our last category, and a configuration is
 said to be a $K_C$ configuration if it is neither a $K_A$ configuration
 nor a $K_B$ configuration. Thus a configuration is $K_C$ if all groups
 are of size at most 2, with any 2-groups being positive, and there are no
 interleaved 2-sequences, nor any $rots$ or $ots$-$rots$ tangles nor any
 negative $ots$-2-groups.

 Our objective is to show that for a given $K_C$ configuration, there is a
 $K_B$ configuration and a finite sequence of $\teetwo$ and/or $OTS$
 operations that will transform the $K_B$ configuration into the
 given $K_C$ configuration. In fact, the following result ensures that
 we can always arrange it so that the first operator in the sequence is $OTS$.

\begin{proposition}\label{kb to kc not by t}
 Every $K_C$ configuration that can be obtained from a $K_B$ configuration
 by an application of any finite sequence of $\teetwo$ operations can also
 be obtained from a $K_B$ configuration by a single application of $OTS$. 
\end{proposition}

\begin{proof}
 It suffices to consider only the case where each application of $\teetwo$
 in the sequence results in a $K_C$ configuration.
 
 Let $C(K)$ denote a $K_C$ configuration which has been obtained from a $K_B$
 configuration $C(K_1)$ by a single application of $\teetwo$. Thus $C(K)$
 has all groups of size at most two, with any 2-groups being positive, and
 there are no interleaved 2-sequences, nor any $rots$ or $ots$-$rots$ tangles
 nor any negative $ots$-2-groups. Furthermore, $C(K_1)$ does have a positive
 2-group $G$ that $\teetwo$ turns to produce $C(K)$. Since $C(K_1)$ is a $K_B$
 configuration, it has no negative groups nor any $rots$ tangles, but it must
 have at least one of the following: a positive group of size at least 3, an
 $ots$-$rots$ tangle, a negative $ots$-2-group, or an interleaved-2-sequence.
 Let us examine the possibilities in light of the information that upon
 a single application of $\teetwo$, the resulting configuration is without any
 of these structures. If $C(K_1)$ had a positive group of size at least 3, then
 any finite sequence of $\teetwo$ applications will not affect it, other
 than to possibly change it to a negative group. Since $C(K)$ has neither, we
 see that $C(K_1)$ does not have any positive groups of size greater than or
 equal to 3. Suppose that $C(K_1)$ contains an $ots$-$rots$ tangle. Since $C(K)$
 contains no negative groups and no $ots$-$rots$ tangles, it must be that the
 application of $\teetwo$ was to the positive 2-group of the $ots$-$rots$ tangle.

 \begin{figure}[ht]
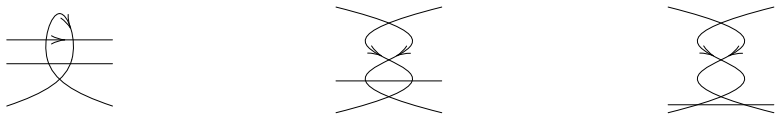

\centering
   \begin{tabular}{c@{\hskip12.5pt}c@{\hskip12.5pt}c}
    $\vcenter{\xy /r20pt/:,
     (0,0);(2,0)**\crv{ (.75,.25) & (1.5,.75) & (1,2.25) & (.5,.75) & (1.25,.25) }
       ?(.38)="x","x"+(-.1,.2);"x"**\dir{}*\dir{>},
     (0,.8);(2,.8)**\dir{-},
     (0,1.25);(2,1.25)**\dir{-}?(.55)*\dir{>},
    \endxy}$
    &
    $\vcenter{\xy /r20pt/:,
    (0,0);(2,2)**\crv{  (2,.5) & (1,1) & (0,1.5)  }?(.55)*{\dir{<}},
    (0,2);(2,0)**\crv{  (2,1.5) & (1,1) & (0,.5)  }?(.45)*{\dir{>}},
    (0,.6);(2,.6)**\dir{-},
    \endxy}$
    &
    $\vcenter{\xy /r20pt/:,
    (0,0);(2,2)**\crv{  (2,.5) & (1,1) & (0,1.5)  }?(.55)*{\dir{<}},
    (0,2);(2,0)**\crv{  (2,1.5) & (1,1) & (0,.5)  }?(.45)*{\dir{>}},
    (0,.15);(2,.15)**\dir{-},
    \endxy}$\\
    \noalign{\vskip 6pt} 
    (a) the $ots$-$rots$ tangle & (b) after turning the 2-group &
      (c) a $K_B$ configuration
\end{tabular}
\caption{Turning the 2-group of an $ots$-$rots$ tangle}
\label{turning the 2-group of an ots-rots tangle}
\end{figure}

\noindent This situation is shown, before and after the turn, respectively,
 in Figure \ref{turning the 2-group of an ots-rots tangle} (a) and (b).
 On the other hand, Figure \ref{turning the 2-group of an ots-rots tangle}
 (c) shows a $K_B$ configuration which will result in the configuration
 illustrated in (b) upon applying the indicated $OTS$ operation. If one
 or more additional $\teetwo$ operations are applied, the end result is
 either to convert the underlying 3-group to a negative group, or to leave
 it unchanged as a positive 3-group. In the first case, the resulting knot
 configuration can also be obtained by a single application of $OTS$ to a
 $K_A$ configuration, which would mean that the resulting knot is $K_B$,
 not $K_C$. Thus this case can't occur. On the other hand, if the underlying
 3-group is positive, then the same argument that we applied for a single
 application of $\teetwo$ will also establish that the configuration that
 results can be obtained by an $OTS$ operation applied to a $K_B$ knot.
 
 We next consider the situation when $K(C_1)$ does not contain any positive
 groups of size greater than or equal to 3 nor any $ots$-$rots$ tangles.
 Suppose that $K(C_1)$ contains a negative $ots$-2-group. Since $K(C)$ does
 not, it must be the case that an application of $\teetwo$ to some positive
 2-group $G$ in $K(C_1)$ will
 cause the underlying 2-group of the $ots$-2-group to become positive.
 However, if we were to apply $OTS$ to $K(C_1)$ to convert the negative
 $ots$-2-group to a negative group $H$ (so the resulting configuration is $K_A$),
 then apply $\teetwo$ to the positive 2-group $G$, the result is a $K_B$
 configuration in which the 2-group $H$ is now positive, and thus $OTS$ may
 be applied to this configuration to obtain $K_C$. Furthermore, if we were
 to apply subsequent $\teetwo$ operations, each time resulting in a
 $K_C$ configuration, we see that none of the subsequent applications
 can change the sign of the underlying 2-group of the positive $ots$-2-group
 back to negative, since then the resulting configuration would be $K_B$.
 Since all positive 2-groups must be non-interleaved, after the $\teetwo$
 operations have been applied, $G$ remains positive, so the same
 argument that applied after the first $\teetwo$ applications is still
 applicable.

 Finally, suppose that $K(C_1)$ does not contain any positive
 groups of size greater than or equal to 3, any $ots$-$rots$ tangles or any
 negative $ots$-2-groups. There must then be an interleaved pair of positive
 2-groups in $K(C_1)$. However, the only way that an application of $T$ to
 a positive 2-group can result in an interleaved pair of positive 2-groups is
 to cause them to change sign. However, $K(C)$ has no negative groups, so this
 situation cannot occur. Obviously, subsequent $\teetwo$ operations which
 are required to produce in a $K_C$ configuration would also not be possible.

 This completes the proof that in every situation where an application of
 any finite sequence of $\teetwo$ operations to a $K_B$ configuration
 results in a $K_C$ configuration $C(K)$, there is a
 corresponding $K_B$ configuration on which $OTS$ can be applied, resulting in
 $C(K)$.
\end{proof}

 Our approach requires that we establish
 some properties of two specific types of subgraphs that are found in
 the plane graph that is obtained when a knot configuration is formed.
 
\begin{definition}
 A plane graph $G$ whose edges are piecewise smooth curves is called a 2-region,
 respectively a minimal loop, if the following conditions are satisfied:
 \begin{alphlist}
  \item $G$ is connected;
  \item $G$ has a face $F$ whose boundary is a cycle $C$;
  \item there are exactly two vertices on $C$ that have degree 2 in $G$
  (called the base vertices of the 2-region), respectively there is exactly
  one vertex on $C$ that has degree 2 in $G$ (called the base vertex of the
  minimal loop);
  \item every non-base vertex on $C$ has degree 3 in $G$;
  \item each vertex that lies in the interior of the region $\mathbb{R}^2-F$
   has degree 4 in $G$.
 \end{alphlist}

 \noindent
  The cycle $C$ is called the {\em boundary} of the 2-region or minimal loop,
  respectively, while the interior of the region $\mathbb{R}^2-F$ is
  called the {\em interior} of the 2-region, respectively minimal loop (note
  that in the event that all vertices of $G$ lie on $C$, then $F$ is not
  uniquely determined--in such a case, let $F$ be the bounded region).
  In the case of a 2-region with boundary cycle $C$ and base vertices $p$ and
  $q$, the two paths between $p$ to $q$ that $C$ determines are called the
  boundary paths of $G$.
  
  More generally, if $G$ is a plane graph with vertex set $V$, then
  a subgraph $G'$ of $G$ with vertex set $V'\subseteq V$ is said to be
  a 2-region (respectively, minimal loop) of $G$ if $G'$ is
  a 2-region (respectively, minimal loop) and every vertex of $G$ that lies in
  the interior of $G'$ belongs to $V'$, and every edge of $G$ that meets the
  interior of $G'$ belongs to $V'$. If $G'$ is a 2-region or minimal loop
  of $G$, then $G$ is said to contain the 2-region or minimal loop $G'$.

  A 2-region $G$ is said to be minimal if $H$ is a 2-region of $G$ implies that
  $H=G$, and a 2-region that has no vertices in its interior is called a 2-group.
  A minimal loop with a single vertex is said to be trivial.
\end{definition}

 Note that the interior of a 2-region or minimal loop could be the unbounded
 region determined by the boundary of the 2-region, respectively minimal loop.
 Further note that a 2-group has exactly two vertices, which are multiply
 connected by two edges.

 Next, observe that by the handshake lemma, there must be an even number of
 non-base vertices on the boundary of a 2-region or minimal loop. There is a
 natural way to pair up the non-base vertices of a 2-region or minimal loop by
 means of what we shall call arcs in the 2-region or minimal loop. 

\begin{definition}
 Let $G$ be a 2-region or minimal loop. An arc in $G$ is a path of maximal
 length which satisfies the following conditions:
 \begin{alphlist}
  \item the initial vertex of the path is a non-base boundary vertex;
  \item the initial edge of the path is the unique non-boundary edge
   of $G$ that is incident to the initial vertex;
  \item each interior (i.e. non-endpoint) vertex of the path is in the
   interior of the 2-region, respectively minimal loop;
  \item at each interior vertex $u$ of the path, the two path edges incident to
   $u$ must be non-adjacent edges at $u$ in $G$;
  \item any non-base boundary vertex on the path is an endpoint of the path. 
  \end{alphlist} 
\end{definition} 
  
\begin{proposition}\label{arcs in 2-region}
 Let $G$ be a 2-region or a minimal loop. Then the endpoints of any arc in $G$
 are distinct and both are non-base boundary vertices of $G$. Furthermore,
 the relation on the set of all non-base boundary vertices of $G$ (if this
 set is nonempty) given by $u$ is related to $v$ if and only if $u=v$ or $u$
 and $v$ are endpoints of an arc of $G$ is an equivalence relation, with
 each equivalence class of size two. Finally, each vertex of $G$ that lies
 in the interior of $G$ is on exactly two arcs in $G$.
\end{proposition}

\begin{proof}
 Suppose that $G$ is a 2-region or minimal loop with nonempty set of non-base
 boundary vertices. Let $P$ be an arc in $G$ with initial vertex $v$ and initial
 edge $e$. We first establish that $P$ is not a closed path. For suppose it is.
 Then there are two edges incident to $v$ that lie in the interior of $G$. Since
 $v$ is a non-base boundary vertex, there are two edges of the boundary cycle
 incident to $v$, whence $v$ has degree at least 4 in $G$. Since this is not
 the case, $P$ is not a closed path. Consequently, the endpoints of $P$ are
 distinct. Let $w$ denote the other endpoint of $P$. We must prove now that $w$
 is a non-base boundary vertex of $G$. Suppose not. Since there are only two
 edges incident to a base vertex, and each such edge has as its other endpoint
 either a base vertex or a non-base boundary, the only way a path could reach
 a base vertex would be to pass through a non-base boundary vertex first. But
 no arc continues on after reaching a non-base boundary vertex, and so $w$ must
 be a vertex in the interior of $G$. Let $f$ denote the last edge of $P$, so
 that $f$ is incident to $w$. Now, $w$ has degree 4, and so there is exactly
 one edge incident to $w$ that is not adjacent to $f$ in $G$. Denote this edge
 by $f'$. If $f'$ were in $P$, then $w$ would not be an endpoint of $P$. Thus
 $f'$ is not in $P$ and so we may extend $P$ by following $f'$, which contradicts
 the maximality of $P$. Thus $w$ is a non-base boundary vertex of $G$.

 The preceding discussion actually establishes that for each non-base boundary
 vertex $v$, there is a unique arc of $G$ with endpoint $v$ (since there
 is only one edge available to start the arc, and only one choice of edge to
 take at each internal vertex of the path). Thus the relation on the set of
 non-base boundary vertices that is defined by $u=v$ or $u$ and $v$ are endpoints
 of an arc of $G$ is an equivalence relation, partitioning the set of
 non-base boundary vertices into cells of size two.

 Furthermore, the discussion of the first paragraph of the proof can also be
 used to establish that each pair of non-adjacent edges incident to a vertex
 in the interior of $G$ lie on exactly one arc in $G$, so there are exactly
 two arcs in $G$ that pass through a given vertex in the interior of $G$.
\end{proof}

\begin{corollary}\label{arc cuts off a 2-region}
 Let $G$ be a 2-region for which some arc $l_1$ in $G$ has both endpoints $x$ and
 $y$ on the same boundary path $L$ of $G$. Let $l_2$ denote the restriction of $L$
 to that part which lies between $x$ and $y$, and let $C$ denote the cycle obtained
 by following $l_1$ by $l_2$. Let $F$ denote the open region of the plane determined
 by $C$ which contains the base points of $G$. Then the graph $G'$ whose vertices
 and edges are those of $G$ which are in $R=\mathbb{R}^2-F$ is a 2-region of $G$.
\end{corollary}

\begin{proof}
  Observe that $C$ lies in the closure of the interior of $G$, but the
  base vertices of $G$ are not on $C$. Thus $R$ is contained
  in the interior of $G$. The boundary edge of $G$ that is incident to $x$
  but which does not belong to $l_2$ lies in the closure of $F$ and is
  therefore not in $G'$, whence $x$ has degree 2 in $G'$. Similarly, $y$
  has degree 2 in $G'$. As well, each vertex of $l_2$ other than $x$ and $y$
  still has degree 3 in $G'$. Each vertex of $l_1$ other than $x$ and $y$
  has two edges of $G$ incident to it, one lying in the closure of $F$, and
  one lying in $R$. Thus each vertex of $l_1$ other than $x$ and $y$ has degree
  3 in $G'$.  Each vertex $v$ in the interior of $R$ is in the interior of $G$,
  and all four edges of $G$ incident to $v$ are in $G'$, whence
  $v$ has degree 4 in $G'$. Furthermore, $v$ lies on an arc
  of $G$, which either has both endpoints on $l_2$ or else must leave $R$ and thus
  meet $l_1$. In either case, each vertex and edge on such an arc between $v$
  and the point of intersection of the arc with either $l_1$ or $l_2$ belongs
  to $G'$, whence there is a walk in $G'$ from $v$ to $x$. Thus $G'$ is connected.
  It follows that $G'$ is a 2-region with interior the interior of $R$, whence
  each vertex of $G$ that is in the interior of $G'$ belongs to $G'$, and so
  $G'$ is a 2-region of $G$.
\end{proof}

 \begin{proposition}\label{existence of 2-region}
  Every non-trivial minimal loop contains a 2-region, and every 2-region contains
  a minimal 2-region.
 \end{proposition}

 \begin{proof}
  Let $G$ be a non-trivial minimal loop with boundary cycle $C$ and
  base vertex $u$. Since $G$ is non-trivial, it must be at least one
  additional vertex. If $G$ has vertices in its interior, then by
  Proposition \ref{arcs in 2-region}, there are arcs in $G$, and so,
  again by Proposition \ref{arcs in 2-region}, there are non-base
  boundary vertices in $G$. Let $v$ be a non-base boundary vertex, and
  let $L_1$ denote the unique arc in $G$ that has $v$ as initial vertex, and
  let $w$ denote the other end point of this arc. Then $v$ and $w$ are
  distinct non-base boundary vertices of $G$. Let $L_2$ denote the path
  obtained by following $C$ from $v$ to $w$ in the direction which does
  not pass through $u$. Then $L_1$ followed by $L_2$ is a cycle $C'$ in $G$.
  Now $C'$ divides the plane into two regions, one of which contains $u$.
  Let $F$ denote this region. Let $H$ be the
  graph whose vertex set consists of all vertices of $G$ that lie either
  on $C'$ or else in the interior of the region $R=\mathbb{R}^2-F$, and whose
  edges are those edges of $G$ that lie in the region $R$ (which includes
  the edges of $C'$). It is straightforward now to verify that $v$ and $w$
  have degree 2 in $H$, any vertex on $C'$ other than $v$ and $w$ has degree
  3 in $H$, every vertex in the interior of $R$ has degree 4 in $H$, and
  $H$ is connected, whence $H$ is a 2-region. But then by construction, $H$
  is a 2-region of $G$.

  Now suppose that $G$ is a 2-region. Of all 2-regions contained in $G$ (and
  $G$ is such a 2-region), let $H$ be one with fewest vertices. Suppose that
  $H'$ is a 2-region of $H$. Let $C_H$ denote the boundary cycle of $H$.
  If $H=C_H$, then $H'=H$ and there is nothing to show. Suppose then that $H$
  has non-boundary vertices. Now all vertices of $H$ lie on $C_H$ or in one of the
  two regions of the plane that $C_H$ determines, and there are vertices of $H$
  that are not on $C_H$. Let $R_H$ denote the region determined by $C_H$ which
  contains vertices of $H$. Now there is a cycle $C_{H'}$ in $H$ such that all
  vertices of $H'$ lie on $C_{H'}$ or in one of the two regions of the plane that
  $C_{H'}$ determines. If $H'=C_{H'}$, then $H'$ is a 2-group and we let $R_{H'}$
  denote the open region determined by $C_{H'}$ that contains no vertices of $H$.
  Otherwise, let $R_{H'}$ denote the open region determined by $C_{H'}$ which
  contains vertices of $H'$. By definition, every vertex of $G$ that is in the
  interior of $H$ belongs to $H$, and every vertex of $H$ that is in the interior
  of $H'$ belongs to $H'$. Let $x$ be a vertex of $G$ that lies in the interior of
  in $R_H$. If $x$ is in $H$, then $x$ is in $H'$, as required.
  Suppose that $x$ is not in $H$. If $R_{H'}\subseteq R_H$, then $x$ would be in the
  interior of $R_H$, whence $x$ belongs to $H$. Since this is not the case,
  it must be that $R_{H'}$ is not contained in $R_H$, whence there is a vertex on
  the boundary of $R_H$ that is contained in the interior of $R_{H'}$. But
  such a vertex must then belong to $H'$ and have degree 4 in $H'$. However,
  vertices on the boundary of $H$ have degree at most 3 in $H$, so this is not
  possible. This contradiction stemmed from the assumption that $x$ was not in
  $H$, so we conclude that $x$ must be in $H$, from which it follows that $x$
  is in $H'$, as required. Thus $H'$ is a 2-region of $H$ such that
  every vertex of $G$ that is contained in the interior of $H'$ belongs to
  $H'$, whence $H'$ is a 2-region of $G$. By the minimality of $H$, we must have
  $H'=H$, whence $H$ is a minimal 2-region.
 \end{proof}

\begin{definition}
 Let $C(K)$ be a configuration of a prime alternating knot $K$, and
 let $v$ be a crossing of $C(K)$. Choose any of the four arcs incident
 to $v$ and construct a closed walk based at $v$ by following a knot
 traversal in the direction of the chosen arc until $v$ is reached for
 the first time. Such a closed walk is called a knot circuit based
 at $v$.
\end{definition}

\begin{proposition}\label{existence of min loops in ck}
 Let $C(K)$ be a configuration of a prime alternating knot $K$. Then
 \begin{alphlist}
  \item For every crossing $v$ of $C(K)$, and each knot circuit $C$
    based at $v$, there is a minimal loop of $C(K)$ whose boundary
    cycle is a subwalk of $C$;
  \item every minimal loop of $C(K)$ contains a 2-region of $C(K)$, and
  \item every 2-region of $C(K)$ contains a minimal 2-region of $C(K)$.
 \end{alphlist} 
\end{proposition}

\begin{proof}
 Let $v$ be a crossing of $C(K)$ and let $C$ denote a knot circuit based at $v$.
 For each crossing $x$ on $C$ (including $x=v$), start at $x$ and follow $C$
 until either $x$ or $v$ is reached, whichever comes first.
 If the result is a circuit based at $x$ (which will be the case if $x\ne v$ and
 $x$ is encountered first), then let $c_x$ denote
 the number of crossings that were encountered (excluding $x$ itself) during
 this traversal, otherwise let $c_x=\infty$. Of all crossings of $C$, choose
 one, say $w$, for which $c_w$ is smallest. Let $C'$ denote the closed walk in
 $C(K)$ that is obtained by beginning at $w$ and following $C$ until $w$ is
 reached again. Then by the minimality of $c_w$, $C'$ must be a cycle.  Let
 $F$ denote the region determined by $C'$ that contains the two edges
 incident to $w$ that are not on $C'$, and let $V'$ be the set consisting of
 all crossings of $C(K)$ that lie on $C'$ or in the interior of the region
 $R=\mathbb{R}^2-F$. Let $H$ be the plane graph whose vertex set is $V'$ and
 whose edge set consists of all edges of $C(K)$ whose endpoints are in $V'$
 and which lie in $R$. We claim that $H$ is a minimal loop in $C(K)$ with
 boundary cycle $C'$.

 First of all, observe that $H$ is a plane graph with a face $F$ whose boundary
 is the cycle $C'$, and that $w$ has degree 2 in $H$. Next, consider a vertex $u$
 on $C'$ different from $w$. Since $C'$ is a cycle, there are two edges of $C'$ that
 are incident to $u$, and since these two edges were successively travelled in
 the knot circuit $C$, they are not adjacent edges at $u$. Consequently,
 one of these two edges, $e_1$ say, lies in $F$ and the other, $e_2$ say, lies in $R$.
 Then $e_2$ belongs to $H$ while $e_1$ does not. It follows that there are exactly
 three edges incident to $u$ in $H$; namely the two on $C'$ and $e_2$, the one which
 lies in $R$. Thus $u$ has degree 3 in $H$. Finally, consider any vertex $y$ of $H$
 which lies in the interior of $R$. Now $y$ has degree 4 in $C(K)$, and all four
 edges that are incident to $y$ have endpoints that are either in the interior of $R$
 or else lie on $C'$, whence all four edges of $C(K)$ that are incident to $y$ belong
 to $H$. Thus $y$ has degree 4 in $H$. Furthermore, any path in $H$ of maximal length
 with initial vertex $y$ must terminate at a vertex of odd degree in $H$. Since the
 only vertices of odd degree in $H$ lie on $C'$, it follows that there is a walk in
 $H$ from $y$ to $w$, whence $H$ is connected. This completes the proof that $H$ is
 a minimal loop with boundary cycle $C'$.

 That $H$ is a minimal loop of $C(K)$ follows now from the fact that the interior
 of $H$ is the interior of $R$, and every vertex of $C(K)$ that is in $R$ was
 put in $V'$, the set of vertices of $H$.

 Thus $H$ is a minimal loop of $C(K)$ with boundary cycle $C'$, which is a closed
 subwalk of $C$.

 Next, suppose that $G$ is a minimal loop of $C(K)$, so there exists a cycle
 $C_G$ in $C(K)$ such that all vertices of $G$ lie in one of the two regions
 of the plane that $C_G$ determines. Let $R_G$ denote this region. By Proposition
 \ref{existence of 2-region}, $G$ contains a 2-region $H$, whence there is
 a cycle $C_H$ in $G$ such that all vertices of $H$ lie in one of the
 two regions of the plane that $C_H$ determines. Let $R_H$ denote this region.
 By definition, every vertex of $C(K)$ that is in the interior of $G$ belongs to $G$, and
 every vertex of $G$ that is in the interior of $H$ belongs to $H$. Let
 $x$ be a vertex of $C(K)$ that lies in the interior of $H$; that is, in the
 interior of $R_H$. If $x$ is in $G$, then $x$ is in $H$, as required.
 Suppose that $x$ is not in $G$. If $R_H\subseteq R_G$, then $x$ would be in the
 interior of $R_G$, whence $x$ belongs to $G$. Since this is not the case,
 it must be that $R_H$ is not contained in $R_G$, whence there is a vertex on
 the boundary of $R_G$ that is contained in the interior of $R_H$. But
 such a verttex must then belong to $H$ and have degree 4 in $H$. However,
 vertices on the boundary of $G$ have degree at most 3 in $G$, so this is not
 possible. This contradiction stemmed from the assumption that $x$ was not in
 $G$, so we conclude that $x$ must be in $G$, from which it follows that $x$
 is in $H$, as required. Thus $H$ is a 2-region contained in $G$ such that
 every vertex of $C(K)$ that is contained in the interior of $H$ belongs to
 $H$, whence $H$ is a 2-region contained in $C(K)$.

 An argument similar to the preceding one will establish that if $G$ is a
 2-region of $C(K)$, and $H$ is a minimal 2-region contained in $G$, then
 every vertex of $C(K)$ that is contained in the interior of $H$ must belong
 to $G$ and thus to $H$, whence $H$ is a minimal 2-region of $C(K)$.
\end{proof} 

Thus every configuration of a prime alternating knot has minimal loops,
2-regions and therefore minimal 2-regions. We remark that an empty 2-region
of a prime alternating knot configuration is a 2-subgroup.

We next observe that the $OTS$ operation on a knot configuration has an
analog for 2-regions and minimal loops. 

\begin{definition}
 Let $G$ be a 2-region or a minimal loop. An $ots$-triangle in $G$ is a
 face of degree 3 none of whose boundary edges belong to the boundary of a
 2-group of $G$, while a face of degree 3 which does have a boundary edge
 in common with a face of degree 2 is called a $rots$-triangle.

 Let $O$ be an $ots$-triangle in $G$, say with boundary edges $B_1$, $B_2$
 and $B_3$. There are three possible situations: none, exactly one, or
 exactly two of the edges $B_1$, $B_2$, $B_3$ is a boundary edge of $G$.
 We define the $ots$ operation in each of these three cases. Let
 $V$ denote the vertex set of $G$ and $E$ denote the edge set of $G$. 
 
 \noindent Case 1: none of $B_1$, $B_2$, $B_3$ is a boundary edge of $C$.
 Now each pair of these edges have exactly one endpoint in common.
 Since a boundary vertex has at most one edge lying in the interior of $G$
 incident to it, we see that none of these three common endpoints is a
 boundary vertex of $G$, and so the compact set $B_1\cup B_2\cup B_3$ is
 contained in the interior of $G$. Furthermore, each of the three common
 endpoints has degree 4, so at each there are two additional incident edges.
 Since the edge curves are smooth, there is an open neighborhood $U$ of $O$
 whose intersection with each of these six additional edges is connected,
 and which does not meet any other edge curve of $G$. Arbitrarily choose
 one pair of boundary edges of $O$, and suppose that the edges were labelled
 so that $B_1$ and $B_2$ are the chosen edges. Let $a$ denote the common
 endpoint of $B_1$ and $B_2$, and let $e$ and $f$ denote the two edges
 that are incident to $a$ in addition to $B_1$ and $B_2$. Further suppose
 that all labelling has been done so that in a clockwise scan at $a$, the edges
 are encountered in the order $B_1$, $e$, $f$ and $B_2$. Choose a point $x\ne a$
 in $e \cap U$, and choose a point $y\ne a$ in $f\cap U$. Let $B_3'$ be a
 smooth curve from $x$ to $y$ within $U$ that does not meet any curve of $G$
 other than $e$ at $x$ and $f$ at $y$. Let $b$ denote the common endpoint of
 $B_2$ and $B_3$, and let $c$ denote the common endpoint of $B_3$ and $B_1$.
 Further, let $g$ and $h$ denote the two edges incident to $b$ other than $B_2$
 and $B_3$, labelled in order in the clockwise direction, and let $i$ and $j$
 denote the two edges incident to $c$ other than $B_3$ and $B_1$, labelled in
 order in the clockwise direction. Let $j'$ denote a smooth curve with endpoint
 $x$ that agrees with $j$ outside of $U$ and within $U$ meets no edge curve of
 $G$ except for $e$ at $x$, and let $g'$ denote a smooth curve with endpoint
 $y$ that agrees with $g$ outside of $U$ and within $U$ meets no edge curve of
 $G$ except for $f$ at $y$. Let the portion of $e$ from $a$ to $x$ be denoted
 by $B_1'$ and denote the remaining portion of $e$ by $e'$. Similarly, let
 the portion of $f$ from $a$ to $y$ be denoted by $B_2'$ and denote the remaining
 portion of $f$ by $f'$. Finally, let $i'=i\cup B_1$ and $h'=h\cup B_2$. Let
 $G'$ denote the plane graph with piecewise smooth edge curves whose vertex
 set is $\left(V-\{\,b,c\,\}\right)\cup\{\,x,y\,\}$ and edge set $\left(E-
 \{\,B_1,B_2,B_3,e,f,g,h,i,j\,\}\right)\cup \{B_1',B_2',B_3',e',f',g',h',i',j'\,\}$.
 Since $G'$ agrees with $G$ outside of $U$, it follows that $G'$ is a 2-region,
 respectively minimal loop, said to be obtained from $G$ by an $ots$ operation
 (on $O$).
 
 \noindent Case 2: exactly one of $B_1$, $B_2$, $B_3$ is a boundary edge of $G$.
 Suppose that the curves were labelled so that $B_3$ is the boundary curve of $G$.
 Label the vertices and edges as in Case 1, with the only differences stemming
 from the fact that the endpoints of $B_3$ are boundary vertices of $G$ (non-base,
 since neither $B_1$ nor $B_2$ is a boundary edge of $G$), so there is no
 edge $h$ or $i$. Carry out the construction as in Case 1 (so there is no
 corresponding $h'$ or $i'$), and let $G'$ denote the plane graph with piecewise
 smooth edge curves whose vertex set is $\left(V-\{\,a,b,c\,\}\right)\cup
 \{\,x,y\,\}$ and edge set $\left(E-\{\,B_1,B_2,B_3,e,f,g,i\,\}\right)\cup
 \{B_3',e',f',g',i'\,\}$. Since $G'$ agrees with $G$ outside of $U$, it follows
 that $G'$ is a 2-region, respectively minimal loop, said to be obtained from $G$
 by an $ots$ operation (on $O$). Note that the number of vertices of $G'$ is one
 less than the number of vertices of $G$. In this case, we say that vertex $a$
 has been $ots$-ed out of $G$.
 
\noindent Case 3: two of the edges $B_1$, $B_2$, $B_3$ are boundary edges of $G$.
 Suppose that the curves have been labelled so that $B_3$ is not a boundary edge
 of $G$. Let the endpoints of $B_3$ be $b$ and $c$. Then $b$ and $c$ are
 non-base boundary vertices, while the common endpoint of $B_1$ and $B_2$
 is a base boundary vertex. Each of $b$ and $c$ have one additional boundary
 edge incident to them. Let $h$ be the additional boundary edge incident to
 $b$ and let $i$ be the additional boundary edge incident to $c$, and set
 $h'=h\cup B_2$ and $i'=i\cup B_1$. Let $G'$ denote the plane graph with
 piecewise smooth edge curves whose vertex set is $V-\{\,b,c\,\}$ and edge
 set $\left(E-\{\,B_1,B_2,B_3,h,i\,\}\right)\cup \{h',i'\,\}$. Since $G'$
 agrees with $G$ outside of $U$, it follows that $G'$ is a 2-region, respectively
 minimal loop, said to be obtained from $G$ by an $ots$ operation (on $O$).

\noindent
 In Cases 2 and 3, we say that the $ots$-triangle is on the boundary of $G$.
\end{definition}

\begin{figure}[ht]
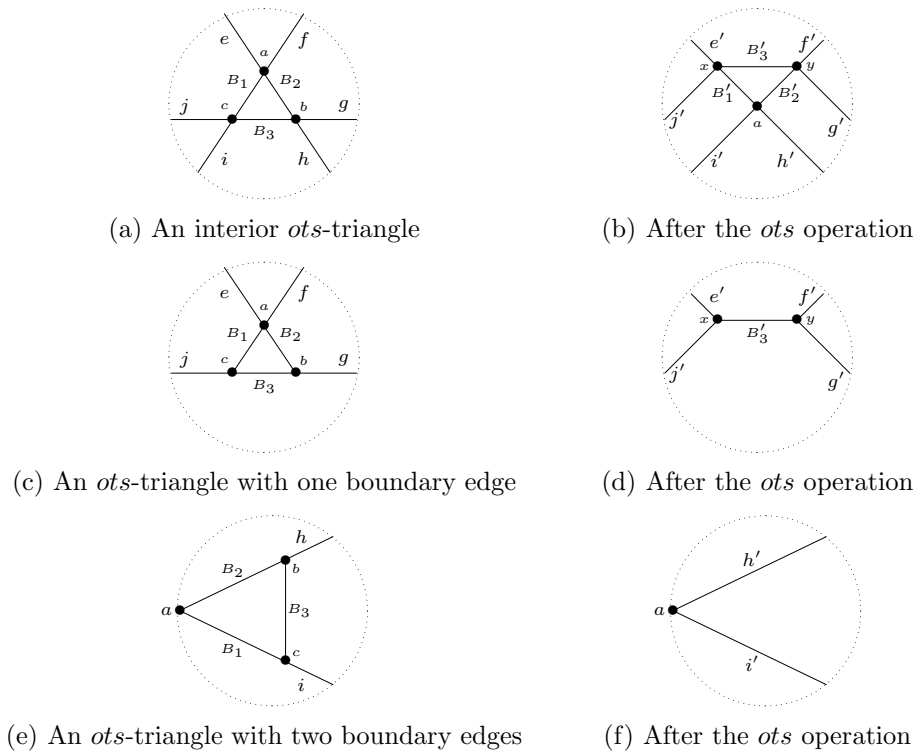

 \centering
   \begin{tabular}{c@{\hskip30pt}c}
$\vcenter{\xy /r40pt/:,
(.625,.65)="a11";"a11";{\ellipse(.9){.}},
(0,0)="a1";(1,1.5)="a2"**\dir{-},
(.25,1.5)="a3";(1.25,0)="a4"**\dir{-},
(-.25,.5)="a5";(1.5,.5)="a6"**\dir{-},
(.25,1.25)*{\hbox{$\ssize e$}},
(1,1.25)*{\hbox{$\ssize f$}},
(1.375,.625)*{\hbox{$\ssize g$}},
(1,.125)*{\hbox{$\ssize h$}},
(.25,.125)*{\hbox{$\ssize i$}},
(-.125,.625)*{\hbox{$\ssize j$}},
(.625,1.125)*{\hbox{$\sssize a$}},
(1,.625)*{\hbox{$\sssize b$}},
(.25,.625)*{\hbox{$\sssize c$}},
(.375,.86)*{\hbox{$\sssize {B_1}$}},
(.875,.86)*{\hbox{$\sssize {B_2}$}},
(.625,.375)*{\hbox{$\sssize {B_3}$}},
(.625,.95)*{\hbox{\small$\bullet$}},
(.325,.5)*{\hbox{\small$\bullet$}},
(.925,.5)*{\hbox{\small$\bullet$}},
\endxy}$
&
$\vcenter{\xy /r40pt/:,
(.625,.65)="a11";"a11";{\ellipse(.9){.}},
(0,0)="a1";(1.25,1.25)="a2"**\dir{-},
(0,1.25)="a3";(1.25,0)="a4"**\dir{-},
(-.25,.5)="a5";(.25,1)="a7"**\dir{-},
"a7";(1,1)="a8"**\dir{-},
"a8";(1.5,.5)="a6"**\dir{-},
(.25,1.25)*{\hbox{$\ssize e'$}},
(1.1,1.225)*{\hbox{$\ssize f'$}},
(1.375,.425)*{\hbox{$\ssize g'$}},
(.9,.125)*{\hbox{$\ssize h'$}},
(.25,.125)*{\hbox{$\ssize i'$}},
(-.125,.5)*{\hbox{$\ssize j'$}},
(.625,.45)*{\hbox{$\sssize a$}},
(.125,1)*{\hbox{$\sssize x$}},
(1.13,1)*{\hbox{$\sssize y$}},
(.3,.75)*{\hbox{$\sssize {B'_1}$}},
(.925,.75)*{\hbox{$\sssize {B'_2}$}},
(.625,1.15)*{\hbox{$\sssize {B'_3}$}},
(.625,.625)*{\hbox{\small$\bullet$}},
(.25,1)*{\hbox{\small$\bullet$}},
(1,1)*{\hbox{\small$\bullet$}},
\endxy}$ \\
\noalign{\vskip 6pt}
(a) An interior $ots$-triangle & (b) After the $ots$ operation \\
\noalign{\vskip 6pt}
$\vcenter{\xy /r40pt/:,
(.625,.65)="a11";"a11";{\ellipse(.9){.}},
(.325,.5)="a1";(1,1.5)="a2"**\dir{-},
(.25,1.5)="a3";(.925,.5)="a4"**\dir{-},
(-.25,.5)="a5";(1.5,.5)="a6"**\dir{-},
(.25,1.25)*{\hbox{$\ssize e$}},
(1,1.25)*{\hbox{$\ssize f$}},
(1.375,.625)*{\hbox{$\ssize g$}},
(-.125,.625)*{\hbox{$\ssize j$}},
(.625,1.125)*{\hbox{$\sssize a$}},
(1,.625)*{\hbox{$\sssize b$}},
(.25,.625)*{\hbox{$\sssize c$}},
(.375,.86)*{\hbox{$\sssize {B_1}$}},
(.875,.86)*{\hbox{$\sssize {B_2}$}},
(.625,.375)*{\hbox{$\sssize {B_3}$}},
(.625,.95)*{\hbox{\small$\bullet$}},
(.325,.5)*{\hbox{\small$\bullet$}},
(.925,.5)*{\hbox{\small$\bullet$}},
\endxy}$
&
$\vcenter{\xy /r40pt/:,
(.625,.65)="a11";"a11";{\ellipse(.9){.}},
(1,1)="a1";(1.25,1.25)="a2"**\dir{-},
(0,1.25)="a3";(.25,1)="a4"**\dir{-},
(-.25,.5)="a5";(.25,1)="a7"**\dir{-},
"a7";(1,1)="a8"**\dir{-},
"a8";(1.5,.5)="a6"**\dir{-},
(.25,1.25)*{\hbox{$\ssize e'$}},
(1.1,1.225)*{\hbox{$\ssize f'$}},
(1.375,.425)*{\hbox{$\ssize g'$}},
(-.125,.5)*{\hbox{$\ssize j'$}},
(.125,1)*{\hbox{$\sssize x$}},
(1.13,1)*{\hbox{$\sssize y$}},
(.625,.875)*{\hbox{$\sssize {B'_3}$}},
(.25,1)*{\hbox{\small$\bullet$}},
(1,1)*{\hbox{\small$\bullet$}},
\endxy}$\\
\noalign{\vskip 6pt}
(c) An $ots$-triangle with one boundary edge & (d) After the $ots$ operation \\
\noalign{\vskip 6pt}
$\vcenter{\xy /r40pt/:,
(.625,.75)="a11";"a11";{\ellipse(.9){.}},
(-.25,.75)="a1";(1.2,1.45)="a2"**\dir{-},
"a1";(1.2,.05)="a3"**\dir{-},
(.75,1.225)="a4";(.75,.275)="a5"**\dir{-},
(-.25,.75)*{\hbox{\small$\bullet$}},
(.75,1.225)*{\hbox{\small$\bullet$}},
(.75,.275)*{\hbox{\small$\bullet$}},
(-.375,.75)*{\hbox{$\ssize a$}},
(.85,1.15)*{\hbox{$\sssize b$}},
(.85,.35)*{\hbox{$\sssize c$}},
(.9,1.45)*{\hbox{$\ssize h$}},
(.9,.05)*{\hbox{$\ssize i$}},
(.25,.375)*{\hbox{$\sssize {B_1}$}},
(.25,1.13)*{\hbox{$\sssize {B_2}$}},
(.88,.75)*{\hbox{$\sssize {B_3}$}},
\endxy}$
&
$\vcenter{\xy /r40pt/:,
(.625,.75)="a11";"a11";{\ellipse(.9){.}},
(-.25,.75)="a1";(1.2,1.45)="a2"**\dir{-},
"a1";(1.2,.05)="a3"**\dir{-},
(-.25,.75)*{\hbox{\small$\bullet$}},
(.5,1.25)*{\hbox{$\ssize h'$}},
(.5,.25)*{\hbox{$\ssize i'$}},
(-.375,.75)*{\hbox{$\ssize a$}},
\endxy}$\\
\noalign{\vskip 6pt}
(e) An $ots$-triangle with two boundary edges & (f) After the $ots$ operation
\end{tabular}
\caption{$ots$ operations in a 2-region or a minimal loop}
\label{ots operations on 2-region}
\end{figure}

For example, in Figure \ref{ots operations on 2-region} (a) and (b), an $ots$-triangle
in the interior of a 2-region or minimal loop and the result of applying an
$ots$ operation to the $ots$-triangle are shown, while in
(c) and (d), an $ots$-triangle with exactly one edge on the boundary of the 2-region
or minimal loop, together with the outcome of an $ots$ operation to this $ots$-triangle
are shown. Finally, in (e) and (f), an $ots$-triangle with two edges on the
boundary of the 2-region or minimal loop and the effect of applying an $ots$
operation to the $ots$-triangle are shown.

We remark that the plane graphs that result from each of the three possible $ots$
operations that may be performed on an $ots$-triangle that is contained in the
interior of a 2-region or minimal loop are isomorphic via an isomorphism that is
the restriction of a homeomorphism of the plane onto itself that fixes all points
in an open neighborhood of the compact set that consists of the $ots$-triangle together
with the six additional edges incident to the vertices of the $ots$-triangle.

The introduction of the notion of $ots$ operation on a 2-region or a minimal loop
is in recognition of the local nature of the $OTS$ operation on a prime knot
configuration. In order to establish that each $K_C$ configuration can be
obtained by a finite sequence of $OTS$ and/or $\teetwo$ operations applied to a
$K_B$ configuration, we shall require the following result about minimal
2-regions.

\begin{theorem}\label{empty min 2-region}
 Given a minimal 2-region, there is a finite sequence of $ots$ operations which 
 will transform the minimal 2-region into an empty 2-region.
\end{theorem}

\begin{proof}
 The proof will be by induction on the number of vertices in the interior of the
 minimal 2-region. Suppose now that we have a
 minimal 2-region $G$ which contains no vertices in its interior.
 If $G$ consists only of the two boundary edges, then $G$ is empty and we are
 done. Suppose that $G$ has edges other than boundary edges. For such an
 edge, both endpoints must be boundary vertices since $G$ has no vertices
 in its interior. Let $p$ and $q$ denote the base vertices of $G$, and let
 $L_1$ and $L_2$ denote the two boundary paths from $p$ to $q$. By Corollary
 \ref{arc cuts off a 2-region}, a non-boundary edge of $G$ must have one endpoint
 on $L_1$ and one endpoint on $L_2$. Of all such edges of $G$, choose the one,
 $e$ say, whose endpoint on $L_1$, $v$ say, is closest to $p$. Then there are no
 vertices on $L_1$ between $p$ and $v$. If the endpoint of $e$ that lies on
 $L_2$ is denoted by $w$, then there can be no vertices on $L_2$ between $p$ and
 $w$, since such a vertex would be the endpoint of an arc of $G$ which must,
 by Corollary \ref{arc cuts off a 2-region}, meet $L_1$ at a vertex between
 $p$ and $x$, which is not possible. Thus $e$ and the two boundary edges incident
 to $p$ form an $ots$-triangle. Let $G'$ denote the 2-region that is formed
 by performing an $ots$ operation on this $ots$-triangle. Then $G'$ is
 a minimal 2-region with no vertices in its interior, having one less non-boundary
 edge than $G$. Induction on the number of non-boundary edges establishes
 the base case of the main induction argument.

 In order to prove the inductive step, we must first establish that any minimal
 2-region with at least one non-boundary edge contains an $ots$-triangle on a
 boundary.
 
 This problem can be modelled as follows. Let $l$ be a line in the plane
 through which pass finitely many given lines in the plane, such that no
 three of the lines (including $l$) are coincident. If some
 pair of lines other than $l$ intersect at a point $X$, then there is a
 triangle on $l$ with third vertex lying on the same side of $l$ as does $X$
 and which is such that none of the lines pass through the interior of this
 triangle. Such a triangle can be found by examining the finite set of
 all triangles based on $l$ that are formed by the given lines and choosing one
 of least height.
 
 We now return to the proof of the theorem. Suppose that $k\ge0$ is an integer
 such that for every minimal 2-region with $k$ or fewer interior vertices,
 there is a finite sequence of $ots$ operations that will transform the minimal
 2-region into an empty 2-region. Consider any minimal 2-region $G$ with
 $k+1$ interior vertices. If there are $ots$-triangles for which two of the
 three edges of the $ots$-triangle are boundary edges of $G$, then we may
 perform an $ots$ operation on such an $ots$-triangle to remove $l$ from $G$,
 leaving a minimal 2-region still with $k+1$ interior vertices. Thus we may
 assume that no $ots$-triangle  of $G$ is such that two of the edges are
 boundary edges of $G$. Since $k+1\ge1$, $G$ has at least one non-boundary edge,
 whence there is an $ots$-triangle on a boundary. Choose one $ots$-triangle on
 a boundary and apply the $ots$ operation to this $ots$-triangle. The result is
 a minimal 2-region $G'$ with $k$ interior vertices, so there is a finite sequence
 of $ots$-operations which will transform $G'$ into an empty 2-region. Thus there
 is a finite sequence of $ots$-operations that will transform $G$ into an empty
 2-region. This completes the proof of the inductive step. The theorem follows now
 by induction.
\end{proof}

Our objective is to prove that for any $K_C$ configuration $C(K)$, there is
a $K_B$ configuration $C(K')$ and a finite sequence of $OTS$ and/or $\teetwo$
operations which will transform $C(K')$ into $C(K)$. For this, we require one
additional result, which involves the notion of a $t$ operation on a
2-region or minimal loop.

\begin{definition}\label{t operation}
 Let $G$ be a 2-region or a minimal loop with edge set $E$, and let $G'$
 be a 2-group contained in $G$ whose base vertices $x$ and $y$ are
 interior vertices of $G$. Let $e$ and $f$ denote the two boundary edges
 of $G'$, labelled in that order clockwise at $x$. Further let the
 remaining two edges of $G$ that are incident to $x$ be denoted by $e_1$
 and $e_2$ in clockwise order, and the remaining two edges of $G$ that
 are incident to $y$ be denoted by $f_1$ and $f_2$ in clockwise order.
 Let $U$ be an open neighborhood of the closure of the interior of $G'$
 which meets no other edges of $G$ than $e_1$, $e_2$, $f_1$ and $f_2$
 (in addition to containing $e$ and $f$). Delete that part of $e_2$
 and of $f_2$ which lies in $U$, except for the endpoints $x$ and $y$,
 and extend $e_2-U$ smoothly into $U$ to terminate at $y$ without
 meeting $e_1$, $e$ or $f_1$, and extend $f_2-U$ smoothly into $U$ to
  terminate at $x$ without meeting $f_1$, $f$ or $e_1$. Denote the
 two extensions of $e_2$ and $f_2$ by $e_2'$ and $f_2'$, respectively.
 Then the graph $H$ whose vertex set is that of $G$ but whose edge set
 is $\left(E-\{\, e_2,f_2\,\}\right)\cup\{\,e_2',f_2'\,\}$ is a 2-region
 or minimal loop respectively, said to be obtained from $G$ by turning
 the 2-group $G'$, or by applying the $t$ operation to $G'$.
\end{definition}

\begin{theorem}\label{reduction of minimal loop}
 Let $G$ be a nontrivial minimal loop. Then there is a finite sequence of $ots$ and/or
 $t$ operations that will transform $G$ into a minimal loop with either a 2-group
 one of whose edges is a boundary edge of the minimal loop, or with two 2-groups
 having a vertex in common, or a $rots$-triangle.
\end{theorem}

\begin{proof}
 The proof will be by induction on the number of vertices in the interior of $G$.
 For the base case, consider any nontrivial minimal loop $G$ with no interior vertices.
 By Proposition \ref{existence of 2-region}, $G$ contains a 2-region $H$, which must
 then have no vertices in its interior and so must be a 2-group in $G$. Since $G$
 itself is not a 2-region, not both boundary edges of $H$ are boundary edges of
 $G$. Thus one boundary edge of the 2-group is a boundary edge of $G$, and the other
 is not, whence $H$ is a 2-group one of whose edges is a boundary edge of $G$.

 Now suppose that $k\ge0$ is an integer such that for any minimal loop with at most $k$
 vertices in its interior, there is a finite sequence of $ots$ and/or $t$ operations that
 will transform the minimal loop into one which contains either a 2-group having one
 edge a boundary edge of the minimal loop, or two 2-groups with a vertex in common, or
 a $rots$-triangle. Consider a minimal loop $G$ with $k+1$ vertices in its interior.
 By Proposition \ref{existence of 2-region}, $G$ contains a minimal 2-region. Of all
 minimal 2-regions of $G$, let $H$ denote one with fewest vertices. By Theorem
 \ref{empty min 2-region}, there is a finite sequence of $ots$ operations that will
 transform $H$ into a 2-group. Suppose that $H$ has an $ots$-triangle $O$. Then $H$
 is not a 2-group, and by the minimality of $H$ in $G$, $G$ contains no 2-groups.
 Thus $O$ is an $ots$-triangle of $G$ and so any $ots$ operation on $O$ in $H$ can
 be obtained as the restriction of an $ots$ operation on $O$ in $G$. Furthermore,
 an $ots$ operation on an $ots$-triangle in $H$ will not create a 2-group outside
 of $H$. Thus there is a finite sequence of $ots$ operations on $G$ that will result
 in a minimal loop $G'$ containing exactly one 2-group. If one of the edges of the
 2-group is a boundary edge of $G'$, we are done. Suppose that this is not the case.
 Let $u$ and $v$ denote the vertices of the 2-group. Then both $u$ and $v$ are
 interior vertices of $G'$. Let the four edges of $G'$ that are not edges of the
 2-group but are incident to either $u$ or $v$ be labelled, in clockwise order, as
 $e_1$, $f_1$, $e_2$, and $f_2$. Let $U$ be an open neighborhood containing the two
 edges of the 2-group which meets no other edges of $G$ except for $e_1$, $f_1$, $e_2$
 and $f_2$. Delete all points on edges of $G$ that lie in $U$, and choose a point
 $x$ in  $U$. Extend $e_1-U$, $f_1-U$, $e_2-U$ and $f_2-U$ smoothly into $U$ 
 without meeting, terminating at $x$. Let $e_1'$, $f_1'$, $e_2'$ and $f_2'$ denote
 these extensions of $e_1$, $f_1$, $e_2$ and $f_2$, respectively. Let $G''$ be
 the plane graph whose vertices are $x$ together with the vertices of $G'$
 except for $u$ and $v$, and whose edges are $e_1'$, $f_1'$, $e_2'$ and $f_2'$,
 together with the edges of $G'$ except for the six edges that are incident to
 either $u$ or $v$. Then $G''$ is a minimal loop with at most $k$ interior
 vertices, whence by the induction hypothesis, there is a finite sequence of
 $ots$ and/or $t$ operations on $G''$ that will transform it into a minimal loop
 with either a 2-group having one edge on the boundary, or two 2-groups with
 a common vertex, or a $rots$-triangle. We shall try to lift each of these
 $ots$ and/or $t$ operations on $G''$ to corresponding operations on $G'$.
 It is clear that if an $ots$ or $t$ operation is being applied to an $ots$-triangle
 or 2-group, respectively, to which the vertex $x$ does not belong, the $ots$
 or $t$ operation can be lifted directly to $G'$. On the other hand, suppose
 that there is an $ots$-triangle $O$ in $G''$ which has $x$ as one of its vertices,
 and an $ots$ operation is to be performed on $O$. 

\begin{figure}[ht]
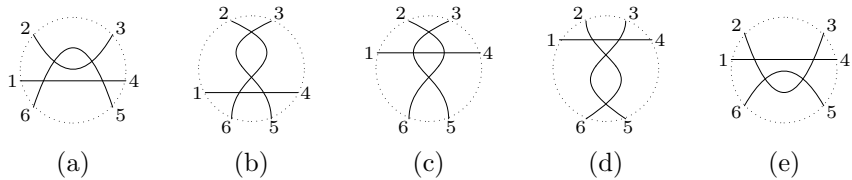

 \centering
   \begin{tabular}{c@{\hskip20pt}c@{\hskip20pt}c@{\hskip20pt}c@{\hskip20pt}c}
$\vcenter{\xy /r20pt/:,
(1,.7)="a11";"a11";{\ellipse(1){.}},
(0,.5)="a1";(2,.5)="a2"**\dir{-},
(.25,1.375)="a3";(1.75,1.375)="a4"**\crv{(.75,.5) & (1.25,.5)},
(.25,0)="a5";(1.75,0)="a6"**\crv{(.75,1.5) & (1.25,1.5)},
(-.15,.5)*{\hbox{$\ssize 1$}},
(.1,1.5)*{\hbox{$\ssize 2$}},
(1.9,1.5)*{\hbox{$\ssize 3$}},
(2.15,.5)*{\hbox{$\ssize 4$}},
(1.9,-.15)*{\hbox{$\ssize 5$}},
(.1,-.15)*{\hbox{$\ssize 6$}},
\endxy}$
&
$\vcenter{\xy /r20pt/:,
(.875,.7)="a11";"a11";{\ellipse(1){.}},
(0,.25)="a1";(1.75,.25)="a2"**\dir{-},
(.5,-.25)="a3";(.5,1.6)="a4"**\crv{(.5,.25) & (1.25,.75) & (1.25,1.25)},
(1.25,-.25)="a5";(1.25,1.6)="a6"**\crv{(1.25,.25) & (.5,.75) & (.5,1.25)},
(-.15,.25)*{\hbox{$\ssize 1$}},
(.35,1.75)*{\hbox{$\ssize 2$}},
(1.4,1.75)*{\hbox{$\ssize 3$}},
(1.9,.25)*{\hbox{$\ssize 4$}},
(1.3,-.4)*{\hbox{$\ssize 5$}},
(.4,-.4)*{\hbox{$\ssize 6$}},
\endxy}$
&
$\vcenter{\xy /r20pt/:,
(.875,.7)="a11";"a11";{\ellipse(1){.}},
(-.1,1)="a1";(1.85,1)="a2"**\dir{-},
(.5,-.25)="a3";(.5,1.6)="a4"**\crv{(.5,.25) & (1.25,.75) & (1.25,1.25)},
(1.25,-.25)="a5";(1.25,1.6)="a6"**\crv{(1.25,.25) & (.5,.75) & (.5,1.25)},
(-.25,1)*{\hbox{$\ssize 1$}},
(.35,1.75)*{\hbox{$\ssize 2$}},
(1.4,1.75)*{\hbox{$\ssize 3$}},
(2,1)*{\hbox{$\ssize 4$}},
(1.3,-.4)*{\hbox{$\ssize 5$}},
(.4,-.4)*{\hbox{$\ssize 6$}},
\endxy}$
&
$\vcenter{\xy /r20pt/:,
(.875,.7)="a11";"a11";{\ellipse(1){.}},
(0,1.25)="a1";(1.75,1.25)="a2"**\dir{-},
(.5,-.25)="a3";(.5,1.6)="a4"**\crv{(1.25,.25) & (1.25,.75) & (.5,1.25)},
(1.25,-.25)="a5";(1.25,1.6)="a6"**\crv{(.5,.25) & (.5,.75) & (1.25,1.25)},
(-.15,1.25)*{\hbox{$\ssize 1$}},
(.35,1.75)*{\hbox{$\ssize 2$}},
(1.4,1.75)*{\hbox{$\ssize 3$}},
(1.9,1.25)*{\hbox{$\ssize 4$}},
(1.3,-.4)*{\hbox{$\ssize 5$}},
(.4,-.4)*{\hbox{$\ssize 6$}},
\endxy}$
&
$\vcenter{\xy /l20pt/:,
(1,.7)="a11";"a11";{\ellipse(1){.}},
(0,.5)="a1";(2,.5)="a2"**\dir{-},
(.25,1.375)="a3";(1.75,1.375)="a4"**\crv{(.75,.5) & (1.25,.5)},
(.25,0)="a5";(1.75,0)="a6"**\crv{(.75,1.5) & (1.25,1.5)},
(-.15,.5)*{\hbox{$\ssize 4$}},
(.1,1.5)*{\hbox{$\ssize 5$}},
(1.9,1.5)*{\hbox{$\ssize 6$}},
(2.15,.5)*{\hbox{$\ssize 1$}},
(1.9,-.15)*{\hbox{$\ssize 2$}},
(.1,-.15)*{\hbox{$\ssize 3$}},
\endxy}$
\\
\noalign{\vskip 6pt}
(a)  & (b) & (c) & (d) & (e)
\end{tabular}
\caption{Lifting an $ots$ operation involving the 2-group vertex}
\label{lifting ots}
\end{figure}

We shall always choose to perform the $ots$ operation that moves the side
opposite the 2-group that is represented by $x$. This might require a $t$
operation if the 2-group is not aligned the right way, followed by two
$ots$ operations to move the edge across both vertices of the 2-group.
Note that since $x$ is an interior vertex of $G''$, after the $ots$
operation on $G''$, $x$ is either still an interior vertex or it has
been $ots$-ed out of $G''$. When we lift the $ots$ operation to $G'$,
the end result is that either both vertices of the 2-group are interior
vertices or else both have been $ots$-ed out of the minimal loop. Thus
every $ots$ operation can be successfully lifted, resulting in possibly
a $t$ operation, but always in two $ots$ operations on $G'$.

Now consider the possibility that a 2-group of $G''$ has $x$ as a vertex,
and the $t$ operator is to be performed on this 2-group. At this stage
in the sequence of $ots$ and/or $t$ operations applied to $G'$, we have
obtained either two 2-groups with a common vertex ($x$ represents one of
the 2-groups), or else we have a $rots$-triangle, depending on the alignment
of the 2-group represented by $x$.

Thus if the sequence of $ots$ and/or $t$ operations that are to be performed
on $G''$ never needs to apply $t$ to a 2-group that has $x$ as one of its
vertices, then the sequence of $ots$ and/or $t$ operations lifts to give a
finite sequence of $ots$ and/or $t$ operations on $G'$ that results in a
minimal loop containing a 2-group with one edge on the boundary of the minimal
loop. On the other hand, if the sequence of $ots$ and/or $t$ operations that are to be performed
on $G''$ does at some point apply $t$ to a 2-group that has $x$ as one of its
vertices, then the sequence of $ots$ and/or $t$ operations lifts to give a
finite sequence of $ots$ and/or $t$ operations on $G'$ that results in a
minimal loop containing either two 2-groups with a common vertex, or a $rots$-triangle.
This completes the proof of the inductive step.
\end{proof}

\begin{definition}
 Let $C(K)$ be a $K_C$ configuration of a prime knot $K$. If $G$ is a
 positive 2-group of $C(K)$ with crossings $u$ and $v$ and edges $e$ and $f$,
 then a knot circuit based at either $u$ or $v$ with initial edge either
 $e$ or $f$ is called a 2-circuit in $C(K)$. The number of crossings on a
 2-circuit is called the weight of the 2-circuit.
\end{definition}

\begin{theorem}\label{effectively emptying a min 2-region}
 Let $C(K)$ be a $K_C$ configuration of a prime alternating knot $K$. Then
 there is a $K_B$ configuration $C(K')$ and a finite sequence of $OTS$ and/or
 $\teetwo$ operations which transforms $C(K')$ into $C(K)$.
\end{theorem}

\begin{proof}
 Let $C(K)$ be a $K_C$ configuration of a prime alternating knot $K$.
 By Proposition \ref{existence of min loops in ck}, $C(K)$ contains a
 minimal 2-region $G$, and then by Theorem \ref{empty min 2-region}, 
 there is a finite sequence of $ots$ operations which will transform $G$ into
 an empty 2-region. Now, an $ots$-triangle which is not on the boundary of a
 minimal 2-region of $C(K)$ is an $OTS$ 6-tangle and the $ots$ operation on
 such an $ots$-triangle is exactly the result of an $OTS$ operation applied to
 the $OTS$ 6-tangle. On the other hand, an $ots$-triangle which is on the
 boundary might fail to be an $OTS$ 6-tangle, but only if it is actually a
 $rots$-tangle. Since $C(K)$ is a $K_C$ configuration, it has no $rots$-tangles,
 so if a $rots$-tangle appears, it was as a result of an $OTS$ operation. Prior
 to such an $OTS$ operation, the $rots$-tangle would have been an $ots$-$rots$
 tangle, so there would be a finite sequence of $OTS$ operations that would
 create an $ots$-$rots$ tangle from crossings of $G$. If no such $ots$-$rots$
 tangle is created as the $OTS$ operations are applied to $G$ to empty it
 out, then we arrive at an empty 2-region, which is then a 2-subgroup. If it
 is negative, then prior to the last $OTS$ operation, we would have
 obtained a negative $ots$-2-subgroup. Otherwise, the finite sequence of $ots$
 operations corresponds to a finite sequence of $OTS$ operations which transform
 the minimal 2-region into a positive 2-subgroup.
 This positive 2-subgroup is either a positive 2-group, or else it is a
 2-subgroup of a positive 3-group (since $C(K)$ had no positive groups of
 size greater than 2). Consequently, from $G$ we determine that either
 there is a finite sequence of $OTS$ operations that transforms $C(K)$ into a $K_B$
 configuration by virtue of creating an $ots$-$rots$ tangle, a negative
 $ots$-2-subgroup, or a positive 3-group, in which case $C(K)$ is obtained by
 the inverse sequence of $OTS$ operations applied to the $K_B$ configuration,
 or else there is a finite sequence of $OTS$ operations that when applied to
 $C(K)$ result in a $K_C$ configuration $C(K')$ containing a positive 2-group.
 The inverse sequence of $OTS$ operations will transform $C(K')$ into $C(K)$.

 It therefore suffices to prove that every $K_C$ configuration that contains
 a positive 2-group is obtained by applying some finite sequence of $OTS$ and/or
 $\teetwo$ operations to some $K_B$ configuration.
 We shall prove by induction on $k$ that
 for every integer $k\ge3$, if $C(K)$ is a $K_C$ configuration of a prime alternating
 knot $K$ with a 2-circuit of weight $k$, then there is a $K_B$ configuration $C(K')$
 and a finite sequence of $OTS$ and/or $\teetwo$ operations that transforms $C(K')$
 into $C(K)$.

 The statement is vacuously true for $k=3$, since
 the existence of a 2-circuit of weight 3 would imply the existence of a 2-tangle,
 and a prime knot has no 2-tangles (we are working with knots of 4 crossings or more).
 Suppose now that $k\ge3$ is an integer such
 that for every $K_C$ configuration $C(K)$ having a 2-circuit of weight at most $k$,
 there is a $K_B$ configuration $C(K')$ and a finite sequence of $OTS$ and/or
 $\teetwo$ operations that will transform $C(K')$ into $C(K)$. Let $C(K)$ be
 a $K_C$ configuration with a 2-circuit $C$ of weight $k+1$. By Proposition
 \ref{existence of min loops in ck}, there is a minimal loop $G$ of $C(K)$ whose
 boundary cycle $C'$ is a subwalk of $C$. By Theorem \ref{reduction of minimal loop},
 there is a finite sequence of $ots$ and/or $t$ operations that will transform $G$
 into a minimal loop with either a 2-group one of whose edges is a boundary edge of
 the minimal loop, or with two 2-groups having a vertex in common, or a
 $rots$-triangle. It is evident that any $t$ operation on a minimal loop or 2-region
 of a knot configuration which contains no negative groups is the restriction of a
 $\teetwo$ operation on the knot configuration, but it is possible that an $ots$-triangle of a minimal loop
 or 2-region of a knot configuration is not an $OTS$ 6-tangle of the knot
 configuration. However, in such a case, the $ots$-triangle would be such that two
 of its vertices were actually vertices of a 2-group in the knot configuration,
 whence the knot configuration would be a $rots$ tangle. Let us start with
 $C(K)$ and follow the sequence of $ots$ and/or $t$ operations which will
 transform $G$ into a minimal loop with either a 2-group one of whose edges
 is a boundary edge of the minimal loop, or with two 2-groups having a vertex in
 common, or a $rots$-triangle. Let $C(K_0)=C(K)$, $C(K_1),\ldots,C(K_n)$ denote
 the configurations formed by the sequence of operations, either carried out to
 completion or until an $ots$ operation is reached where the corresponding
 $ots$-triangle is actually part of a $rots$ tangle in the knot configuration
 produced by the preceding operation. Suppose that $C(K_n)$ is
 not a $K_C$ configuration. Since $C(K_0)$ is a $K_C$ configuration, there is a
 least index $i$ such that $C(K_{i-1}$ is a $K_C$ configuration but $C(K_i)$ is
 not a $K_C$ configuration. Since no single $OTS$ or $\teetwo$ operation will
 convert a $K_C$ configuation into a $K_A$ configuration, it must
 be that $C(K_i)$ is a $K_B$ configuration. Thus the sequence of operations is reversible (it suffices
 to observe that no $\teetwo$ operation resulted in a group of size greater than
 2) and so we have a finite sequence of $OTS$ and/or $\teetwo$ operations which
 will convert the $K_B$ configuration $C(K_i)$ into $C(K)$.

 Otherwise, $C(K_n)$ is a $K_C$ configuration, whence it contains no $rots$ tangles.
 It follows that the sequence of $ots$ and/or $t$ operations has been successfully
 carried out as $OTS$ and/or $\teetwo$ operations on $C(K)$. Moreover, the inverse
 of the sequence of $OTS$ and/or $\teetwo$ operations is a sequence of $OTS$ and/or
 $\teetwo$ operations that will transform $C(K_n)$ into $C(K)$. It therefore will
 suffice to establish that there is a $K_B$ configuration and a finite sequence of
 $OTS$ and/or $\teetwo$ operations that will transform the $K_B$ configuration
 into $C(K_n)$. Now, since $C(K_n)$ is a $K_C$ configuration, it does not contain
 any 3-groups nor $rots$ 6-tangles, nor any negative groups, so we must have obtained
 a minimal loop containing a positive 2-group $H$ with one of its edges on the boundary
 of the minimal loop. Moreover, if the edge of $H$ which is not on the boundary of the
 minimal loop is not on the 2-circuit $C$, then $H$ and the 2-group that gave
 rise to $C$ would constitute an interleaved 2-sequence, whence $C(K_n)$ would
 be a $K_B$ configuration. Since $C(K_n)$ is a $K_C$ configuration, it follows that
 the edge of $H$ which is not on the boundary of the minimal loop is on the 2-circuit
 $C$. Let $v$ denote the first vertex of $H$ to be encountered as $C$ is traversed. 
 Let $C'$ denote the walk obtained by starting at $v$ and following $C$ along the
 first edge incident to $v$ that is traversed by $C$ until $v$ is reached again
 (as it must be, since both edges of $H$ will be traversed by $C$). Then $C'$
 is a 2-circuit of $C(K_n)$. Moreover, since every crossing of $C'$ is a crossing
 of $C$, and the initial vertex of the 2-circuit $C$ does not belong to $C'$, it
 follows that the weight of $C'$ is less than $k+1$, the weight of $C$. But then
 $C(K_n)$ is a $K_C$ configuration with a 2-circuit of weight at most $k$, whence
 by the induction hypothesis, there is a $K_B$ configuration $C(K')$ and a finite
 sequence of $OTS$ and/or $\teetwo$ operations which will transform $C(K')$ into
 $C(K_n)$, as required.

 This completes the proof of the induction step, and the result follow now by
 induction.
\end{proof}

 We remark that for any $K_C$ configuration, it follows from Proposition
 \ref{kb to kc not by t}, that there is a $K_B$ configuration and
 a sequence of $OTS/T$ operations, the first of which is $OTS$, which
 when applied to the $K_B$ configuration results in the given $K_C$
 configuration.

\par

\section{Summary}
  This completes the proof that the algorithm does indeed produce all
  prime alternating knots of minimal crossing size $n+1$ from those
  of minimal crossing size $n$. In the sequel, an implementation of
  the algorithm is presented.

\section{References}

\noindent [1] B. Arnold, M. Au, C. Candy, K. Erdener, J. Fan, R. Flynn,
R. J. Muir, D. Wu, and J. Hoste, Tabulating Alternating Knots through
14 Crossings, J. Knot Theory and its Ram., 3(4), 1994,433--437.

\noindent [2] J. A. Calvo, Knot enumeration through flypes
and twisted splices. J. Knot and its Ram., 6(1997), no. 6, 785--798.

\noindent [3] J. H. Conway, An enumeration of knots and links and
some of their related properties, Computational Problems in Abstract
Algebra (John Leech, ed.), Pergamon Press, Oxford and New York,
1969, 329--358.

\noindent [4] O. T. Dasbach and Stefan Hougardy, Does the Jones Polynomial
Detect Unknottedness?, Experiment. Math. 6(1997), no. 1, 51--56.

\noindent [5] H. de Fraysseix and P. Ossona de Mendez, On a characterization
of Gauss codes, Discrete Comput. Geom. 22 (1999), no. 2, 267--295.

\noindent [6] C. H. Dowker and M. B. Thistlethwaite, Classification of knot
projections, Topology Appl. 16(1983), 19--31.

\noindent [7] C. Ernst and D. W. Sumners, The growth of the number of prime
knots, Math. Proc. Camb. Phil. Soc. 102(1987), 303--315.

\noindent [8] J. Hoste, M. Thistlethwaite, J. Weeks, The First 1,701,936 Knots,
 Math. Intelligencer 20(1998), no. 4, 33--48.

\noindent [9] T. P. Kirkman, The 364 unifilar knots of ten crossings enumerated and
defined, Trans. Roy. Soc. Edinburgh 32(1885), 483-506.

\noindent [10] ---, The enumeration, description and construction of knots of
fewer than ten crossings, Trans. Roy. Soc. Edinburgh 32(1885), 281-309.

\noindent [11] J. C. G\'omez-Larra\~naga, Graphs of tangles. Trans. Amer.
 Math. Soc. 286(1984), no. 2, 817--830.

\noindent [12] W. B. R. Lickorish, Prime Knots and Tangles, Trans. Amer.
Math. Soc. 267(1981), no. 1, 321--332.

\noindent [13] --- Surgery on knots. Proc. Amer. Math. Soc. 60(1976), 296--298 (1977).

\noindent [14] C. N. Little, On knots with a census for order ten, Trans. Connecticut
   Acad. 7(1885), 27--43.

\noindent [15] W. W. Menasco and M. B. Thistlethwaite, The classification of
alternating links, Ann. Math., 1993, 113--171.

\noindent [16]  ---, The Tait flyping
conjecture. Bull. Amer. Math. Soc. (N.S.) 25(1991), no. 2, 403--412.

\noindent [17] K. A. Perko, On the classification of knots, Proc. Amer. Math.
Soc. 45(1974), 262-266.

\noindent [18] ---, On 10-crossing knots, Portugal. Math. 38(1979), 5--9.

\noindent [19] P. G. Tait, On knots I, II, III, Scientific Papers, Vol. I,
  Cambridge Univ. Press, London, 1898, 273--347.
  
\noindent [20] M. B. Thislethwaite, Knot tabulation and related topics,
Aspects of Topology (I. M. James and E. H. Kronheimer, eds.), London Math. Soc.
Lecture Notes Ser., 93, Cambridge Univ. Press, 1985, 1--76.

\noindent [21] Carl Sundberg and Morwen Thistlethwaite, The rate of growth
 of the number of prime alternating links and tangles. Pacific J. Math.
 182(1998), no. 2, 329--358.

\noindent [22] D. J. A. Welsh, On the Number of Knots and Links, Colloquia
Math. Soc. J\'anos Bolyai, Sets, Graphs and Numbers, Budapest, 1991,
713--718.

\end{document}